\documentclass[11pt]{article}              

\usepackage[margin=25mm]{geometry}
\usepackage{graphicx}
\usepackage{colortbl}
\usepackage{subfigure}
 
\usepackage{fancyhdr}
\usepackage{amssymb,amsfonts,amsmath,amsthm}
\usepackage{natbib}
\usepackage{bstnotations}

\usepackage{authblk}

\graphicspath{{./},{./Figures/}}

\newcommand{\hsp}{\hspace{0.3mm}}

\newcommand{\nval}{n_{\rm{val}}}
\newcommand{\vli}{v_l^{(i)}}
\newcommand{\vri}{v_r^{(i)}}
\newcommand{\zkli}{z_{k,l}^{(i)}}
\newcommand{\zkri}{z_{k,r}^{(i)}}
\newcommand{\Pki}{P_k^{(i)}}
\newcommand{\Pkj}{P_k^{(j)}}
\newcommand{\Imax}{I_{\rm max}}
\newcommand{\Derrmin}{\Delta\widehat{err}_{\rm min}}
\newcommand{\Ropt}{R^{\rm opt}}
\newcommand{\RCV}{R^{\rm CV3}}
\newcommand{\RLOO}{R^{\rm LOO}}
\newcommand{\errCV}{\widehat{err}_{\rm CV3}}
\newcommand{\errLOO}{\widehat{err}_{\rm LOO}}
\newcommand{\errG}{\widehat{err}_G}
\newcommand{\LOOPCE}{\widehat{err}_{\rm LOO}^{\rm PCE}}
\newcommand{\CVLRA}{\widehat{err}_{\rm CV3}^{\rm LRA}}
\newcommand{\errGLRA}{\widehat{err}_G^{\rm LRA}}
\newcommand{\errGPCE}{\widehat{err}_G^{\rm PCE}}

\begin{document}
\title{Polynomial meta-models with canonical low-rank approximations: numerical insights and comparison to sparse polynomial chaos expansions} 

\author[1]{K. Konakli} \author[1]{B. Sudret}

\affil[1]{Chair of Risk, Safety and Uncertainty Quantification,
  Department of Civil Engineering, ETH Zurich, Stefano-Franscini-Platz
  5, 8093 Zurich, Switzerland}

\date{}
\maketitle

\abstract{The growing need for uncertainty analysis of complex computational models has led to an expanding use of meta-models across engineering and sciences. The efficiency of meta-modeling techniques relies on their ability to provide statistically-equivalent analytical representations based on relatively few evaluations of the original model. Polynomial chaos expansions (PCE) have proven a powerful tool for developing meta-models in a wide range of applications; the key idea thereof is to expand the model response onto a basis made of multivariate polynomials obtained as tensor products of appropriate univariate polynomials. The classical PCE approach nevertheless faces the ``curse of dimensionality'', namely the exponential increase of the basis size with increasing input dimension. To address this limitation, the sparse PCE technique has been proposed, in which the expansion is carried out on only a few relevant basis terms that are automatically selected by a suitable algorithm. An alternative for developing meta-models with polynomial functions in high-dimensional problems is offered by the newly emerged low-rank approximations (LRA) approach. By exploiting the tensor-product structure of the multivariate basis, LRA can provide polynomial representations in highly compressed formats. Through extensive numerical investigations, we herein first shed light on issues relating to the construction of canonical LRA with a particular greedy algorithm involving a sequential updating of the polynomial coefficients along separate dimensions. Specifically, we examine the selection of optimal rank, stopping criteria in the updating of the polynomial coefficients and error estimation. In the sequel, we confront canonical LRA to sparse PCE in structural-mechanics and heat-conduction applications based on finite-element solutions. Canonical LRA exhibit smaller errors than sparse PCE in cases when the number of available model evaluations is small with respect to the input dimension, a situation that is often encountered in real-life problems. By introducing the conditional generalization error, we further demonstrate that canonical LRA tend to outperform sparse PCE in the prediction of extreme model responses, which is critical in reliability analysis. \\[1em] 

}

\maketitle

\section{Introduction}

It is nowadays common practice to study the behavior of physical and engineering systems through computer simulation. Proper analysis of the system response must account for the prevailing uncertainties in the system model and the underlying phenomena, which requires repeated simulations under varying scenarios for the input parameters. Modern advances in computer science combined with the improved understanding of physical laws are leading to computational models of increasing complexity. Uncertainty propagation through such models may become intractable in cases when a single simulation is computationally demanding. A remedy is to substitute a complex model with a \textit{meta-model} that possesses similar statistical properties, but has a simple functional form.

The focus of the present work is on meta-models that are built with polynomial functions due to the simplicity and versatility they offer. A popular class of meta-models thereof are the so-called polynomial chaos expansions (PCE) \citep{Xiu2002wiener, Ghanem2003stochastic}. The key idea of PCE is to expand the model response onto an appropriate basis made of orthonormal multivariate polynomials, the latter obtained as tensor products of univariate polynomials in each of the input parameters. In non-intrusive approaches that are of interest herein, the coefficients of the expansion are evaluated in terms of the response of the original model at a set of points in the input space, called the experimental design \citep{Choi2004a, Berveiller2006a, Xiu2009a}. Although PCE have proven powerful in a wide range of applications, they face limitations in cases with high-dimensional input. This is because the number of the basis terms, and thus of the unknown expansion coefficients, grows exponentially with the number of input parameters, which is commonly referred to as the ``curse of dimensionality''. As shown in \citet{BlatmanJCP2011, Doostan2011non}, the efficiency of the PCE approach can be significantly improved by using a sparse basis.

A promising alternative for developing meta-models with polynomial functions in high-dimensional spaces is provided by canonical decompositions. In canonical decompositions, also known as separated representations, a tensor is expressed as a sum of rank-one components. This type of representation constitutes a special case of tensor decompositions, which are typically used to compress information or extract a few relevant modes of a tensor; a survey on different types of tensor decompositions can be found in \citet{Kolda2009tensor}. The original idea of canonical decomposition dates back to 1927 \citep{Hitchcock1927}, but became popular in the second half of the 20th century after its introduction to psychometrics \citep{Carroll1970a, Harshman1970foundations}. Since then, it has been used in a broad range of fields, including chemometrics \citep{Appellof1981strategies, Bro1997parafac}, neuroscience \citep{Mocks1988topographic, Andersen2004structure}, fluid mechanics \citep{Felippa1990mixed, Ammar2006new}, signal processing \citep{Sidiropoulos2000parallel,DeLathauwer2007tensor}, image analysis \citep{Shashua2001linear, Furukawa2002appearance} and data mining \cite{Acar2006collective, Beylkin2009multivariate}. More recently, canonical decompositions are attracting an increasing interest in the field of uncertainty quantification \citep{Nouy2010b, Chevreuil2013, Doostan2013, Hadigol2014, Mathelin2014, RaiThesis, Validi2014, Chevreuil2015}.

By exploiting the tensor-product structure of the multivariate polynomial basis, canonical decompositions can provide equivalent to PCE representations in highly-compressed formats. It is emphasized that the number of parameters in canonical decompositions grows only linearly with the input dimension, which, in cases of high-dimensional problems, results in a drastic reduction of the number of unknowns compared to PCE. Naturally, canonical decompositions with a few rank-one components are of interest, thus leading to the name low-rank approximations (LRA). We note that although the present study is constrained to the use of polynomial bases, different basis functions may be considered for the construction of LRA in a general case (see, \eg \citet{Chevreuil2015}).

Recently proposed methods for building canonical LRA meta-models in a non-intrusive manner rely on the sequential updating of the polynomial coefficients along separate dimensions. The underlying algorithms require solving a series of minimization problems of small size, independent of the input dimension, which can be easily handled using standard techniques. However, the LRA construction involves open questions that call for further investigations. In particular, stopping criteria in the sequential updating of the polynomial coefficients as well as criteria for selection of the optimal rank and polynomial degree are not yet well established. Considering a particular greedy algorithm for building canonical LRA meta-models, we herein shed light on the aforementioned issues through extensive numerical investigations. In the sequel, we assess the comparative accuracy of canonical LRA and sparse PCE in applications involving finite-element models pertinent to structural mechanics and heat conduction. In these applications, sparse PCE are built with a state-of-art method, where a candidate basis is defined by means of a hyperbolic truncation scheme and the final sparse basis is determined using least angle regression. Comparisons between the meta-model errors are carried out for experimental designs of varying sizes drawn with Sobol sequences and Latin hypercube sampling.

The organization of the paper is as follows: In Section 2, we present the mathematical setup of non-intrusive meta-modeling and describe corresponding error measures. Sections 3 and 4 respectively describe the sparse PCE and canonical LRA approaches. After investigating open questions in the construction of LRA in Section 5, we confront the two types of polynomial meta-models in Section 6. The paper concludes with a summary of the main findings and respective outlooks in Section 7.

\section{Non-intrusive meta-modeling}

\subsection{Mathematical setup}
\label {sec:setup}

We consider a physical or engineering system whose behavior is represented by a computational model $\cm$. Let $\ve X=\{X_1 \enum X_M\}$ and $\ve Y=\{Y_1 \enum Y_N\}$ respectively denote the $M$-dimensional input vector and the $N$-dimensional response vector of the model. In order to account for the uncertainty in the input and the resulting uncertainty in the response, the elements of $\ve X$ and $\ve Y$ are described by random variables. For the sake of simplicity, we hereafter restrain our analysis to the case of a scalar model response, \ie  $N=1$. Note that the case of a vector model response can be addressed by separately treating each element of $\ve Y$ as in the case of a scalar response. The above are summarized in the mathematical formalism:
\begin{equation}
\label{eq:model}
\ve X \in \cd_{\ve X} \subset \Rr^M\longmapsto Y=\cm (\ve X) \in \Rr,
\end{equation}
where $\cd_{\ve X}$ denotes the support of $\ve X$. 

In a typical real-life application, the model $\cm$ is not known in a closed analytical form and may represent a complex computational process. A meta-model $\widehat{\cm}$ is an analytical function such that $\widehat{Y}=\widehat{\cm}(\ve X)$ possesses similar statistical properties with $Y$. Non-intrusive methods for building meta-models rely on a series of calls to the original model $\cm$, which may be used as a ``black-box'' without any modification. Building a meta-model in a non-intrusive manner requires an experimental design (ED), \ie a set of realizations of the input vector $\ce=\{\ve \chi^{(1)} \enum \ve \chi^{(N)}\}$, and the corresponding model evaluations $\cy=\{\cm(\ve \chi^{(1)}) \enum \cm(\ve \chi^{(N)})\}$.

\subsection{Error measures}
\label{sec:err}

For a set of realizations of the input vector $\cx=\{\ve x_1 \enum \ve x_n\}\subset \cd_{\ve X}$ and two real-valued functions $a$ and $b$ with common domain $\cd_{\ve X}$, we define the semi-inner product:
\begin{equation}
\label{eq:semi-inner}
\innprod a b {\cx}=\frac{1}{n}\sum_{i=1}^n a(\ve x_i) \hsp b(\ve x_i).
\end{equation}
Eq.~(\ref{eq:semi-inner}) leads to the semi-norm $\norme a {\cx}= <a,\hsp a>_{\cx}^{1/2}$, which is employed in the sequel to describe meta-model error measures.

A good measure of the accuracy of a meta-model is the generalization error $Err_G$, which represents the mean-square of the residual $\vare=Y-\widehat Y$:
\begin{equation}
\label{eq:errG}
Err_G =\Esp{\vare^2}=\Esp{\left(Y-\widehat Y \right)^2}.
\end{equation}
In most practical situations, it is not possible to evaluate the generalization error analytically. An estimator $\widehat{Err}_G$ of this error may be computed via Monte Carlo Simulation (MCS) by using the exact-model and meta-model responses at a sufficiently large set of points in the input space $\cx_{\rm{val}}=\{{\ve x}_1 \enum {\ve x}_{\nval}\}$, called \emph{validation set}:
\begin{equation}
\label{eq:ErrG_hat}
\widehat{Err}_G =\left\|\cm-\widehat \cm\right\|_{\cx_{\rm{val}}}^2.
\end{equation}
The corresponding estimator of the relative generalization error, denoted by $\widehat{err}_G$, is obtained by normalizing $\widehat{Err}_G$ with the empirical variance of $\cy_{\rm{val}}=\{Y({\ve x}_1)\enum Y({\ve x}_{\nval})\}$, the latter representing the set of model responses at the validation set.

Unfortunately, a validation set is not available in typical meta-modeling applications, where a large number of evaluations of the original model is non-affordable. An alternative estimator is the empirical error $\widehat{Err}_E$:
\begin{equation}
\label{eq:ErrE}
\widehat{Err}_E =\left\|\cm-\widehat \cm\right\|_{\ce}^2,
\end{equation}
where the subscript $\ce$ indicates that the semi-norm is evaluated at the points of the ED. The corresponding relative error, denoted by $\widehat{err}_E$, is obtained by normalizing $\widehat{Err}_E$ with the empirical variance of $\cy=\{\cm({\ve \chi}^{(1)})\enum \cm({\ve \chi}^{(N)})\}$, the latter representing the set of model responses at the ED. The empirical error does not require any additional evaluations of the exact model than those already used to build the meta-model. It therefore serves the goal of limiting the number of runs of an expensive computational model to the smallest possible. However, it tends to underestimate the actual generalization error, which might be severe in cases of overfitting.

Relying on the ED only, a fair approximation of the generalization error can be obtained by means of cross-validation (CV) techniques. In $k$-fold CV, the ED is randomly partitioned into $k$ sets of approximately equal size. A meta-model is built considering all but one of the partitions, which comprises the training set, while the excluded or testing set, is used to evaluate the generalization error. By alternating through the $k$ sets, $k$ meta-models are obtained in this way; their average generalization error provides an estimate of the error of the meta-model built with the full ED.

\section{Polynomial Chaos Expansions}

\subsection{Spectral representation}
\label{sec:spectral}

Assuming that $Y$ in Eq.~(\ref{eq:model}) has a finite variance, the following representation is possible \citep{Soize2004}:
\begin{equation}
\label{eq:spectral}
Y = \sum_{j=0}^{\infty}y_j\phi_j,
\end{equation}
in which $\{\phi_j, j \in \Nn\}$ is a set of random variables that forms a Hilbertian basis in the space of second-order variables and $\{y_j, j \in \Nn\}$ are the coordinates of $Y$ in this basis. Eq.~(\ref{eq:spectral}) constitutes a \emph{spectral representation} of the random response $Y$.

Let us consider the Hilbert space $\ch$ of square-integrable real-valued functions of $\ve X$ equipped with the inner product:
\begin{equation}
\label{eq:inn_prod_H}
\innprod uv \ch = \int_{\cd_{\ve X}}{u(\ve x) v(\ve x) f_{\ve X}(\ve x) \di{\ve x}},
\end{equation}
where $f_{\ve X}$  denotes the joint probability density function (PDF) of $\ve X$. Let us also consider the Hilbert space $\ch_i$ of square-integrable real-valued functions of $X_i$ equipped with the inner product:
\begin{equation}
\label{eq:inn_prod_Hi}
\innprod uv {\ch_i} = \int_{\cd_{X_i}}{u(x) v(x) f_{X_i}(x_i) \di{x_i}},
\end{equation}
where $f_{X_i}$ denotes the marginal PDF of $X_i$ defined over the support $\cd_{X_i}$.

Under the assumption that the components of $\ve X$ are independent, it is straightforward to show that the Hilbert space defined by the tensor product  $\bar{\ch}=\otimes_{i=1}^{M} {\ch}_i$ and equipped with the inner product:
\begin{equation}
\label{eq:inn_prod_Hbar}
\innprod uv {\bar{\ch}} = \int_{\cd_{\ve X}}{u(x_1 \enum x_M) v(x_1 \enum x_M) f_{X_1}(x_1)\ldots f_{X_M}(x_M) \di{x_1} \ldots \di{x_M}}
\end{equation}
is isomorphic with $\ch$. Accordingly, if $\{\psi_{\alpha_i}^{(i)},\hsp \alpha_i \in \Nn\}$ is a basis in the Hilbert space $\ch_i$, then a basis in the Hilbert space $\ch$ is defined by the set $\{\Psi_{\ua},\hsp \ua=(\alpha_1 \enum \alpha_M) \in \Nn^M\}$ with:
\begin{equation}
\label{eq:basis}
\Psi_{\ua}(\ve X)=\prod_{i=1}^M\psi_{\alpha_i}^{(i)}(X_i).
\end{equation}
It follows that the problem of specifying a Hilbertian basis in $\ch$ reduces to specifying Hilbertian bases in ${\cal{H}}_i$, $i=1 \enum M$. This proposition will be employed in the sequel for building the bases of polynomial chaos expansions. Although the above analysis has been constrained to the case when the components of $\ve X$ are independent, it will be seen that cases with mutually dependent input variables may be treated similarly after an appropriate isoprobabilistic transformation.

\subsection{Construction of polynomial basis}
\label{sec:PCE_basis}

Polynomial chaos expansions (PCE) are spectral representations in which the basis consists of multivariate polynomials in $\ve X$ that are orthonormal with respect to $f_{\ve X}$. According to Section~\ref{sec:spectral}, these can be obtained as tensor products of univariate polynomials that constitute Hilbertian bases in the spaces of $X_i$, $i=1 \enum M$.

A Hilbertian basis $\{\psi_k^{(i)},\hsp k \in \Nn\}$ of ${\ch}_i$ satisfies the orthonormality condition:
\begin{equation}
\label{eq:orthonorm_cond}
\innprod {\psi_j^{(i)}} {\psi_k^{(i)}} {\ch_i} =\delta_{jk},
\end{equation}
where $\delta_{jk}$ is the Kronecker delta symbol. Classical algebra allows one to build a family of orthogonal polynomials $\{Q_k^{(i)}, k\in\Nn\}$, where $k$ denotes the polynomial degree, so that they satisfy \citep{Abramowitz}:
\begin{equation}
\label{eq:orthog_cond}
\innprod {Q_j^{(i)}} {Q_k^{(i)}} {\ch_i}=c_k^{(i)}\delta_{jk}.
\end{equation}
In the above equation, the constant $c_k^{(i)}$ represents the squared $L_2$ norm of the $k$-th degree polynomial:
\begin{equation}
\label{eq:poly_norm}
c_k^{(i)}=\|Q_k^{(i)}\|_2^2= \innprod {Q_k^{(i)}} {Q_k^{(i)}} {{\ch}_i}.
\end{equation}
The corresponding orthonormal polynomial families are thus obtained through the normalization:
\begin{equation}
\label{eq:orthonorm_poly}
P_k^{(i)}=Q_k^{(i)}/\sqrt{c_k^{(i)}}.
\end{equation}
Once the families of univariate orthonormal polynomials associated with the elements of $\ve X$ have been determined, the multivariate polynomial basis can be obtained through the tensorization shown in Eq.~(\ref{eq:basis}), after substituting $\psi_{\alpha_i}^{(i)}$ with $P_{\alpha_i}^{(i)}$.

For standard distributions, \ie uniform, Gaussian, Gamma, Beta, the associated families of orthogonal polynomials are well-known \citep{Xiu2002wiener}. For instance, a uniform variable with support $[-1,1]$ is associated with the family of Legendre polynomials, whereas a standard normal variable is associated with the family of Hermite polynomials. However, it is common in practical situations that the input variables do not follow standard distributions. In such cases, the random vector $\ve X$ is first transformed into a basic random vector $\ve X'$ (\eg a standard normal or standard uniform random vector) through an isoprobabilistic transformation $\ve X=T^{-1}(\ve X')$ and then, the model response $\cm (T^{-1}(\ve X'))$ is expanded onto the polynomial basis associated with $\ve X'$. This approach also allows dealing with mutually dependent input variables through an isoprobabilistic transformation into a vector of independent variables (\eg Nataf transformation in the case of joint PDF with Gaussian copula).

The exact representation of the random response  requires an infinite number of basis terms. However, in practical implementation of PCE, an approximation containing a finite number is considered:
\begin{equation}
\label{eq:PCE}
\widehat{Y}^{\rm PCE}=\widehat{\cm}^{\rm PCE}(\ve X)=\sum_{\ua \in \ca}{y_{\ua}} \Psi_{\ua}(\ve X),
\end{equation}
where the set of retained multi-indices $\ca$ is determined according to an appropriate truncation scheme. A typical truncation scheme consists in selecting multivariate polynomials up to a total degree $p^t$, \ie $\ca=\{\ua \in \Nn^M: \|\ua\|_1  \leq p^t\}$, with $\|\ua\|_1=\sum_{i=1}^M{\alpha_i}$. The corresponding number of terms in the truncated series is:
\begin{equation}
\label{eq:card} 
{\rm{card}} \ca={M+p^t \choose p^t} =\frac{(M+p^t)!}{M!p^t!}.
\end{equation}
This number increases exponentially with the input dimension $M$ giving rise to the ``curse of dimensionality". To limit the number of basis terms that include interactions between input variables, which are usually less significant, \citep{BlatmanPEM2010} proposed the use of a hyperbolic truncation scheme. In the latter, the set of retained multi-indices is defined as $\ca=\{\ua \in \Nn^M: \|\ua\|_q \leq p^t\}$, with:
\begin{equation}
\label{eq:qnorm} 
\|\ua\|_q=\left(\sum_{i=1}^M{\alpha_i}^q\right)^{1/q}, \hspace{2mm} 0<q\leq1.
\end{equation}
According to Eq.~(\ref{eq:qnorm}), lower values of $q$ correspond to a smaller number of interaction terms in the PCE basis. At the limit $q\rightarrow0$, the expansion becomes additive, \ie a sum of univariate functions in the $X_i$'s.

\subsection{Computation of polynomial coefficients}
\label {sec:PCE_coef}

Next, we briefly review a non-intrusive approach for computing the coefficients $y_{\ua}$ based on least-square analysis, originally introduced by \citet{Choi2004a, Berveiller2006a} under the name \emph{regression method}. In this approach, the exact expansion is viewed as the sum of a truncated series and a residual:
\begin{equation}
\label{eq:PCE_res}
Y=\cm(\ve X)=\sum_{\ua \in \ca}y_{\ua}\Psi_{\ua}(\ve X)+\vare,
\end{equation}
where $\vare$ corresponds to the truncated terms. Then, the set of coefficients $\ve y=\{y_{\ua},\hsp \ua \in \ca\}$ can be obtained by minimizing the mean-square error of the residual over the ED:
\begin{equation}
\label{eq:PCE_coef}
\ve y= \mathrm{arg} \underset{\ve {\upsilon}\in\Rr^{\rm{card}\ca}}{\mathrm{min}}\left\|\cm-\sum_{\ua \in \ca}\upsilon_{\ua}\Psi_{\ua}\right\|_{\ce}^2,
\end{equation}
leading to:
\begin{equation}
\label{eq:OLS_solution}
\ve y= (\ve \Psi^{\rm T} \ve \Psi)^{-1} \ve \Psi^{\rm T}\cy,
\end{equation}
where $\ve \Psi=\{\ve \Psi_{ij}=\Psi_j(\ve \chi^{(i)}),\hsp i=1\enum N, \hsp j=1\enum \rm{card}\ca\}$ and $\cy$ is the set of model responses evaluated at the ED, as defined earlier.

For high-dimensional problems, the number of coefficients to be evaluated can be very large. For instance, in typical engineering problems where $M$ varies in the range $10-50$, by considering only low-degree polynomials with, say, $p^t=3$, Eq.~(\ref{eq:card}) results in $286-23,426$ unknown coefficients. Obviously, the number of model evaluations required to solve Eq.~(\ref{eq:PCE_coef}), which is typically 2-3 times the number of unknowns, becomes prohibitively large in cases with high-dimensional input. This limitation constitutes a bottleneck in the classical PCE approach. 

More efficient schemes for evaluating the PCE coefficients can be devised by considering the respective regularized problem:
\begin{equation}
\label{eq:PCE_coef_reg}
\ve y= \mathrm{arg} \underset{\ve {\upsilon}\in\Rr^{\rm{card}\ca}}{\mathrm{min}}\left\|\cm-\sum_{\ua \in \ca}\upsilon_{\ua}\Psi_{\ua}\right\|_{\ce}^2+\lambda \cp(\ve \upsilon),
\end{equation}
in which  $\cp(\ve \upsilon)$ is an appropriate regularization functional of $\ve \upsilon=\{\upsilon_1 \enum \upsilon_{\rm{card}\ca}\}$. If $\cp(\ve \upsilon)$ is selected as the $L_1$ norm of $\ve \upsilon$, \ie $\cp(\ve \upsilon)=\sum_{i=1}^{\rm{card}\ca} \abs {\upsilon_i}$, insignificant terms may be disregarded from the set of predictors, leading to \emph{sparse} solutions. \citep{BlatmanJCP2011} proposed to use the hybrid least angle regression (LAR) method for building sparse PCE. This method employs the LAR algorithm \citep{Efron2004} to select the best set of predictors and subsequently, estimates the coefficients using ordinary least squares (OLS), as described above in Eq.~(\ref{eq:OLS_solution}).

It will be seen in Section~\ref{sec:LRA} that canonical low-rank approximations offer an alternative approach for dealing with high-dimensional problems. By exploiting the tensor-product form of the polynomial basis, such representations may reduce the number of unknown coefficients by orders of magnitude.

\subsection{Accuracy estimation}
\label{sec:PCE_err}

A good measure of the PCE accuracy is the leave-one-out (LOO) error \citep{Allen1971}, corresponding to the CV error for the case $k=N$ (see Section~\ref{sec:err} for details on the CV technique). Using algebraic manipulations, this can be computed based on a \emph{single} PCE built with the full ED. Let $h(\ve \chi^{(i)})$ denote the $i$-th diagonal term of matrix $\ve \Psi(\ve \Psi^{\rm{T}}\ve \Psi)^{-1}\ve \Psi^{\rm{T}}$. The LOO error can then be computed as \citep{BlatmanThesis}:
\begin{equation}
\label{eq:errLOO_simp}
\widehat{Err}_{\rm LOO}=\left\|\frac{\cm-\widehat{\cm}^{\rm PCE}}{1-h}\right\|_{\ce}^2.
\end{equation}
The relative LOO error is obtained by normalizing $\widehat{Err}_{\rm LOO}$ with the empirical variance of the model responses at the ED, denoted by $\Varhat{\cy}$. Because this error tends to be too optimistic, the following corrected estimate is used instead \citep{Chapelle2002}:
\begin{equation}
\label{eq:errLOO_corr}
\errLOO=\frac{\widehat{Err}_{\rm LOO}} {\Varhat{\cy}} \left(1-\frac{\rm{card}\ca}{N}\right)^{-1}\left(1+\rm{tr}((\ve \Psi^{\rm{T}}\ve \Psi)^{-1})\right).
\end{equation}

\section{Canonical low-rank approximations}
\label{sec:LRA}

\subsection{Formulation with polynomial bases}

We herein consider again the mapping in Eq.~(\ref{eq:model}). A rank-one function of the input vector $\ve X$ has the form:
\begin{equation}
\label{eq:rank1}
w(\ve X)= \prod_{i=1}^M {v^{(i)}(X_i)},
\end{equation}
where $v^{(i)}$ denotes a univariate function of $X_i$. A representation of the model response $Y=\cm (\ve X)$ as a finite sum of rank-one functions constitutes a canonical decomposition with rank equal to the number of rank-one components. A rank-$R$ decomposition of $Y=\cm(\ve X)$ therefore reads:
\begin{equation}
\label{eq:CP}
\widehat{Y}_R = \widehat{\cm}_R(\ve X)=\sum_{l=1}^R b_l\left(\prod_{i=1}^M {\vli(X_i)}\right),
\end{equation}
where $\vli$ denotes a univariate function of $X_i$ in the $l$-th rank-one component and $\{b_l, \hsp l=1 \enum R\}$ are scalars that can be viewed as normalizing constants.

An exact canonical decomposition represents a \emph{rank decomposition}. In general, the rank decomposition of a given tensor is not unique; conditions of uniqueness are discussed in \citet{Kolda2009tensor}. The lowest rank of a rank decomposition of a given tensor, called the \textit{tensor rank} \citep{Hitchcock1927, Kruskal1977three}, can be determined numerically by fitting various canonical decomposition models \citep{Kolda2009tensor}. Naturally, of interest are decompositions where the exact response is approximated with sufficient accuracy by using a relatively small number of terms $R$. Such decompositions are hereafter called low-rank approximations (LRA).

The focus of the present work is on canonical LRA made of polynomial functions due to the combination of simplicity and versatility these offer; in a general case however, the use of polynomial functions is not a constraint (see \eg \citet{Chevreuil2015}). By expanding  $\vli$ onto a polynomial basis that is orthonormal with respect to the marginal distribution $f_{X_i}$, Eq.~(\ref{eq:CP}) takes the form:
\begin{equation}
\label{eq:LRA_pol}
\widehat{Y}_R=\widehat{\cm}_R(\ve X)=\sum_{l=1}^R b_l \left(\prod_{i=1}^M\left(\sum_{k=0}^{p_i} \zkli \hsp \Pki (X_i)\right)\right),
\end{equation}
where $\Pki$ denotes the $k$-th degree univariate polynomial in the $i$-th input variable, $p_i$ is the maximum degree of $\Pki$ and $\zkli$ is the coefficient of $\Pki$ in the $l$-th rank-one component. Appropriate families of univariate polynomials according to the distributions of the respective input variables are determined as discussed in Section~\ref{sec:PCE_basis}. Similarly to PCE, the case of dependent input can be treated through an isoprobabilistic transformation of the input variables $\{X_i, i=1 \enum N\}$ into independent reduced variables.

Disregarding the redundant parameterization arising from the normalizing constants, the number of unknowns in Eq.~(\ref{eq:LRA_pol}) is $R \cdot \sum_{i=1}^M (p_i+1)$, which grows \emph{only linearly} with the input dimension $M$. Thus, a representation of the model response in the form of canonical LRA results in a drastic reduction of the number of unknowns as compared to PCE. To emphasize this, we consider PCE with the candidate basis determined by the truncation scheme $\ca=\{\ua \in \Nn^M: \alpha_i\leq p_i, \hsp i=1 \enum M\}$, so that the expansion relies on the same polynomial functions as those used in Eq.~(\ref{eq:LRA_pol}). For the case when $p_i=p$, $i=1 \enum M$, the resulting number of unknowns is $(p+1)^M$ in PCE versus $(p+1)\cdot M \cdot R$ in LRA. Assuming a typical engineering problem with $M=10$ and low-degree polynomials with $p=3$, these formulas yield $1,048,576$ unknowns in PCE versus $40 R$ unknowns in LRA; for a low rank, say $R\leq 10$, the latter number does not exceed a mere $400$. The reduction in the number of unknowns achieved with the compressed LRA representation becomes even more pronounced in cases with high-degree polynomials.

\subsection{Greedy construction}
\label{sec:LRA_constr}

Different algorithms have been proposed in the literature for building a decomposition in the form of Eq.~(\ref{eq:LRA_pol}) in a non-intrusive manner \citep{Chevreuil2013, Doostan2013, Mathelin2014, RaiThesis, Validi2014, Chevreuil2015}. A common point is that the polynomial coefficients are determined by means of an alternated least-squares (ALS) minimization. The ALS approach consists in sequentially solving a least-squares minimization problem along a single dimension $i\in\{1 \enum M\}$, while ``freezing'' the coefficients in all remaining dimensions. In the present study, we adopt the skeleton of the greedy algorithm proposed in \citet{Chevreuil2013, Chevreuil2015}. This algorithm involves a progressive increase of the rank by successively adding rank-one components up to a prescribed maximal. It thus results in a set of candidate decompositions with varying ranks and requires appropriate criteria for selecting the optimal one. In the following, we describe the construction of the approximation for a prescribed rank and discuss rank selection criteria.

\subsubsection{Approximation for a prescribed rank}
\label{sec:LRA_alg}

Let us denote by $\widehat{\cm}_r$ the rank-$r$ approximation of $\cm$. The corresponding model response is then approximated by:
\begin{equation}
\label{eq:Y_r}
\widehat{Y}_r=\widehat{\cm}_r(\ve X)=\sum_{l=1}^{r}b_l w_l(\ve X),
\end{equation}
where $w_l$ represents the $l$-th rank-one component:
\begin{equation}
\label{eq:w_l}
w_l(\ve X)=\prod_{i=1}^M\left(\sum_{k=0}^{p_i} \zkli \hsp \Pki (X_i)\right).
\end{equation}
The employed algorithm involves a sequence of pairs of a \textit{correction step} and an \textit{updating step}, so that in the $r$-th correction step, the rank-one tensor $w_r$ is built, while in the $r$-th updating step, the set of normalizing coefficients $\{b_1 \enum b_r\}$ is determined. These steps are detailed next.

\textbf{Correction step}: In the $r$-th correction step, the rank-one tensor $w_r$ is obtained as the solution to the minimization problem:
\begin{equation}
\label{eq:solve_wr}
w_r(\ve X)=\mathrm{arg} \underset{\omega \in \cw}{\hsp\mathrm{\min}} \left\|\cm-\widehat{\cm}_{r-1}-\omega\right\|_{\ce}^2,
\end{equation}
where $\cw$ represents the space of rank-one tensors. The sequence is initiated by setting $Y_0=\widehat{\cm}_0(\ve X)=0$. Eq.~(\ref{eq:solve_wr}) is solved by means of an ALS scheme that involves successive minimizations along each dimension $i=1\enum M$. In the minimization along dimension $j$, the polynomial coefficients in all other dimensions are ``frozen'' at their current values and the coefficients ${\ve z}_r^{(j)}=\{z_{1,r}^{(j)} \ldots z_{p_j,r}^{(j)}\}$ are determined as:
\begin{equation}
\label{eq:solve_zr}
\ve z_r^{(j)}= \mathrm{arg}\underset{\ve \zeta \in\Rr^{p_j+1}}{\mathrm{min}}\left\|\cm-\widehat{\cm}_{r-1}-\left(\prod_{i\neq j}\vri \right)\left(\sum_{k=0}^{p_j} \zeta_k \hsp \Pkj\right)\right\|_{\ce}^2,
\end{equation}
where:
\begin{equation}
\label{eq:v_r}
\vri(X_i)=\sum_{k=0}^{p_i} \zkri \hsp \Pki(X_i).
\end{equation}
To initiate the $r$-th correction step, one needs to assign arbitrary values to $\vri$, $i=1\enum M$; in the subsequent example applications, we use $v_r^{(1)}(X_1)=\ldots=v_r^{(M)}(X_M)=1$. We underline that a correction step may involve several iterations over the set of dimensions $\{1 \enum M\}$. This aspect of the algorithm, which is not detailed in some of the aforementioned studies, can be critical for the LRA accuracy, as shown in the numerical investigations in Section \ref{sec:Examples_LRA}. \citet{Chevreuil2015} proposed a stopping criterion that combines thresholds on the number of iterations $I_r$ and on the empirical error:
\begin{equation}
\label{eq:Err_r}
\widehat{Err}_r=\left\|\cm-\widehat{\cm}_{r-1}-w_{r}\right\|_{\ce}^2.
\end{equation}
A threshold on the empirical error was also imposed by \citet{Doostan2013} and \citet{Validi2014}. However, appropriate values for these thresholds as well as their effects on the LRA accuracy were not examined. The stopping criterion employed in the present study involves the number of iterations $I_r$ and the decrease in the relative empirical error in two successive iterations, denoted by $\Delta \widehat{err}_r$. The relative empirical error is computed by normalizing $\widehat{Err}_r$ in Eq.~(\ref{eq:Err_r}) with $\Varhat{\cy}$, \ie with the empirical variance of the model responses at the ED. Accordingly, the algorithm exits the $r$-th correction step if either $I_r$ reaches a maximum allowable value $\Imax$ or $\Delta \widehat{err}_r$ becomes smaller than a prescribed threshold $\Derrmin$. Appropriate values for $\Imax$ and $\Derrmin$ will be later discussed based on the numerical investigations in Section \ref{sec:Examples_LRA}.

\textbf{Updating step}: After the completion of a correction step, the algorithm moves to an updating step, in which the set of coefficients $\ve b=\{b_1 \ldots b_r\}$ is obtained as the solution of the minimization problem:
\begin{equation}
\label{eq:solve_b}
\ve b =\mathrm{arg} \underset{\ve \beta \in\Rr^{r}}{\mathrm{min}}\left\|\cm-\sum_{l=1}^r \beta_l w_l\right\|_{\ce}^2.
\end{equation}
Note that in each updating step, the size of vector $\ve b$ is increased by one. In the $r$-th updating step, the value of the new element $b_r$ is determined for the first time, whereas the values of the existing elements $\{\beta_1 \enum \beta_{r-1}\}$ are updated.

Construction of a rank-$R$ decomposition in the form of Eq.~(\ref{eq:LRA_pol}) requires repeating pairs of a correction and an updating step for $r=1\enum R$. The algorithm is summarized below.\\

\textbf{Algorithm 1} Non-intrusive construction of a polynomial rank-$R$ approximation of $Y=\cm(\ve X)$ with an experimental design $\ce=\{\ve \chi^{(1)} \enum \ve \chi^{(N)}\}$:
\begin{enumerate}
	\item Set $\widehat{\cm}_0(\ve \chi^{(q)})=0$, $q=1 \enum N$.
	\item For $r=1\enum R$, repeat steps (a)-(d):
	\begin{enumerate}
		\item Initialize: $\vri(\chi^{(q)})=1$, $i=1 \enum M$, $q=1 \enum N$; $I_r=0$; $\Delta \widehat{err}_r=\epsilon>\Derrmin$.
		\item 
		While $\Delta \widehat{err}_r>\Derrmin$ and $I_r<I_{\rm max}$, repeat steps i-iv:
		\begin{enumerate}
			\item Set $I_r \leftarrow I_r+1$.
			\item For $i=1 \enum M$, repeat steps A-B:
			\begin{enumerate}
				\item Determine  ${\ve z}_r^{(i)}=\{z_{1,r}^{(i)} \ldots z_{p_i,r}^{(i)}\}$ using Eq.~(\ref{eq:solve_zr}).
				\item Update $\vri$, using Eq.~(\ref{eq:v_r}).
			\end{enumerate}
			\item Update $w_r$ using Eq.~(\ref{eq:w_l}).
			\item Compute $\widehat{Err}_r$ using Eq.~(\ref{eq:Err_r}) and update $\Delta \widehat{err}_r$.
		\end{enumerate}
		\item Determine $\ve b=\{b_1 \ldots b_r\}$ using Eq.~(\ref{eq:solve_b}).
		\item Evaluate $\widehat{\cm}_r(\ve \chi^{(q)})$, $q=1 \enum N$, using Eq.~(\ref{eq:Y_r}).
	\end{enumerate}
\end{enumerate}

We emphasize that the above algorithm relies on the solution of several small-size minimization problems; in particular, one needs to solve $M$ minimization problems of size $\{p_i+1, \hsp i=1 \enum M\}$ in each iteration of a correction step (note that $p_i<20$ in typical applications) and one minimization problem of size $r$ in the $r$-th updating step (recall that small ranks are of interest in LRA). Thus, the LRA construction substitutes the single large-size minimization problem involved in the PCE construction with a series of small-size ones; in high-dimensional applications, this can also offer a significant advantage in terms of required computer memory.

Because of the small number of unknowns involved in Eq.~(\ref{eq:solve_zr}) and Eq.~(\ref{eq:solve_b}), these minimization problems can be efficiently solved with OLS, as shown in the subsequent example applications. An alternative approach employed in \citep{Chevreuil2013, RaiThesis, Chevreuil2015} is to substitute these equations with the respective regularized problems:
\begin{equation}
\label{eq:solve_zr_reg}
\ve z_r^{(j)}= \mathrm{arg}\underset{\ve \zeta \in\Rr^{p_j}}{\mathrm{min}}\left\|\cm-\widehat{\cm}_{r-1}-\left(\prod_{i\neq j}\vri\right)\left(\sum_{k=0}^{p_j} \zeta_k \hsp \Pkj\right)\right\|_{\ce}^2+\lambda \cp (\ve \zeta)
\end{equation}
and
\begin{equation}
\label{eq:solve_b_reg}
\ve b =\mathrm{arg} \underset{\ve \beta \in\Rr^{r}}{\mathrm{min}}\left\|\cm-\sum_{l=1}^r \beta_l w_l\right\|_{\ce}^2+\lambda \cp (\ve \beta),
\end{equation}
where $\cp(\ve u)$ is selected as either the $L_1$ or the $L_2$ norm of $\ve u$.

We conclude this section by briefly referring to alternative algorithms for developing LRA non-intrusively. The algorithms in \citet{Doostan2013,Mathelin2014,Validi2014} involve a progressive increase of the rank as well. In \citet{Doostan2013,Validi2014}, when the $r$-th rank-one component is added, the polynomial coefficients in the $(r-1)$ previously built rank-one terms are also updated. Thus, the $r$-th step requires the solution of minimization problems of size $\{\left(p_i+1\right)r, \hsp i=1 \enum M\}$; in applications with high-dimensional input, the size of these minimization problems remains orders of magnitude smaller compared to those involved in the computation of the PCE coefficients. Conversely, the algorithm in \citet{Mathelin2014} does not require the updating of coefficients in previously added components when a new rank-one term is added. It is noteworthy that in the latter, the most-relevant dimensions to be included in the ALS scheme are determined by means of a LAR-based selection technique and CV error criteria.

\subsubsection{Rank selection}
\label{sec:rank}

In a typical application, the optimal rank $R$ is not known \textit{a priori}. As noted earlier, the progressive construction of LRA described in Section~\ref{sec:LRA_alg} results in a set of decompositions of increasing rank. Thus, one may set $r=1 \enum r_{\rm max}$ in Step 2 of Algorithm 1, where $r_{\rm max}$ is a maximum allowable candidate rank, and at the end, select the optimal-rank decomposition using error-based criteria. Two approaches for rank selection are described below and investigated extensively in Section~\ref{sec:Examples_LRA}.

\citet{Chevreuil2015} proposed to select the optimal rank by means of 3-fold CV (see Section \ref{sec:err}); they did not however examine the accuracy of this approach. In the general case of $k$-fold CV, the procedure requires building LRA of increasing rank $r=1 \enum r_{\rm max}$ for each of the $k$ training sets. For $r=1 \enum r_{\rm max}$, the generalization error of each of the $k$ meta-models is estimated using the respective testing set. The rank $R\in \{1 \enum r_{\rm max}\}$ yielding the smallest average generalization error over the $k$ meta-models is identified as optimal; then, a new decomposition of rank $R$ is built using the full ED. The average generalization error corresponding to the selected rank provides an estimate of the actual generalization error of the final meta-model.

Although the above procedure for rank selection relies on the set of model evaluations at the ED only, it requires repeating Algorithm 1 $(k+1)$ times, namely $k$ times for rank $r_{\rm max}$ and another time for rank $R$. Because it is of interest to limit the computational effort required for rank selection, we consider an alternative approach based on the LOO error in the updating step. As described in Section~\ref{sec:PCE_err}, the LOO error represents the average generalization error in $k$-fold CV with $k=N$, but can be computed \emph{from a single meta-model} built with the full ED.

Let $\ve W$ denote the information matrix in the minimization problem represented by Eq.~(\ref{eq:solve_b}), \ie $\ve W=\{W_{ij}=w_j(\ve \chi^{(i)}), \hsp i=1 \enum N, \hsp j=1 \enum r \}$, and let $h(\ve \chi^{(i)})$ denote the $i$-th diagonal term of matrix $\ve W(\ve W^{\rm{T}}\ve W)^{-1}\ve W^{\rm{T}}$. Then, the LOO error can be computed as:
\begin{equation}
\label{eq:errLOO_LRA}
\widehat{Err}_{\rm LOO}=\left\|\frac{\cm-\widehat{\cm}_r}{1-h}\right\|_{\ce}^2.
\end{equation}
The relative LOO error is obtained by normalizing $\widehat{Err}_{\rm LOO}$ with $\Varhat{\cy}$, \ie with the empirical variance of the model responses at the ED. As noted in Section~\ref{sec:PCE_err}, this error may be too optimistic. For this reason, a corrected error may be used in the case of PCE (see Eq.~(\ref{eq:errLOO_corr})). A similar correction may not be formally applied in the present case, because this is based on the assumption of orthogonality between the regressors, which does not hold for Eq.~(\ref{eq:solve_b}). However, it is of interest to examine the performance of the corrected error estimate, given by:
\begin{equation}
\label{eq:errLOO_corr_LRA}
\widehat{err}_{\rm LOO}=\frac{\widehat{Err}_{\rm LOO}}{\Varhat{\cy}} \left(1-\frac{r}{N}\right)^{-1}\left(1+\rm{tr}((\ve W^{\rm{T}}\ve W)^{-1})\right),
\end{equation}
used in an approximate sense.

\section{Example applications of canonical low-rank approximations}
\label{sec:Examples_LRA}

In the following, we investigate the construction of LRA according to Algorithm 1 in Section~\ref{sec:LRA_alg} considering four models with different characteristics and input dimensionality. The first model is an analytical function of rank-1 structure with a 5-dimensional input, representing the deflection of a simply-supported beam under a static load. The subsequent three examples involve finite-element models with input dimensionality 10, 21 and 53.  In particular, we develop LRA for the deflection of a truss structure subjected to vertical static loads, the top horizontal displacement of a three-span five-story frame subjected to horizontal static loads and the temperature response in stationary heat conduction with thermal conductivity described by a random field. In each application, (i) we assess the performance of the two criteria for rank selection discussed in Section~\ref{sec:rank} and evaluate the accuracy of the respective error measures and (ii) we investigate optimal stopping criteria in the correction step. The minimization problems in both the correction and the updating steps are solved with the OLS method. In all applications, effects of varying ED size are examined. The EDs are herein obtained using Sobol quasi-random sequences \citep{Niederreiter1992} (LHS designs will be considered in the following section). In order to assess the accuracy of the different rank-selection methods and stopping criteria, we use sufficiently large validation sets drawn with MCS. We note however that such large validation sets are not available in real-life problems, in which the analyst typically needs to rely solely on the ED.

\subsection{Beam deflection}
\label{sec:Beam_LRA}

In the first example, we consider a simply-supported beam with a constant rectangular cross-section subjected to a concentrated load at the midpoint of the span. The response quantity of interest is the mid-span deflection given by:
\begin{equation}
\label{eq:beam_u}
u = \frac{P L^3}{4 E b h^3},
\end{equation}
where $b$ and $h$ respectively denote the width and height of the cross-section, $L$ is the length of the beam, $E$ is the Young's modulus and $P$ is the magnitude of the load. The aforementioned parameters are described by independent random variables, thus leading to an uncertainty propagation problem of dimension $M=5$. The distributions of the input parameters are listed in Table~\ref{tab:beam_input}. Employing an isoprobabilistic transformation of $\ve X=\{b,h,L,E,P\}$ into a vector of standard normal variables, we develop LRA representations of the random response $U=\cm(\ve X)$ with basis functions made of Hermite polynomials.

\begin{table} [!ht]
	\centering
	\caption{Beam-deflection problem: Distributions of input random variables.}
	\label{tab:beam_input}
	\begin{tabular}{c c c c}
		\hline Variable & Distribution & Mean & CoV \\
		\hline $b~[\rm m]$ & Lognormal & $0.15$ & $0.05$ \\
		$h~[\rm m]$ & Lognormal & $0.3$  & $0.05$ \\
		$L~[\rm m]$ & Lognormal & $5$    & $0.01$ \\
		$E~[\rm MPa]$ & Lognormal & $30,000$  & $0.15$ \\
		$P~[\rm MN]$ & Lognormal & $0.01$ & $0.20$ \\
		\hline       
	\end{tabular}
\end{table}

\subsubsection{Rank selection and error measures}
\label{sec:Beam_rank}

We begin our analysis by investigating the selection of the optimal rank among a set of candidate values $\{1 \enum 20\}$. After preliminary investigations, we set a common polynomial degree $p_i=p=5$ for $i=1 \enum 5$. The stopping criterion in the correction step is defined by setting $\Imax=50$ and $\Derrmin= 10^{-8}$ (the selection of $\Derrmin$ will be explained in the following subsection). We use EDs of size $N$ varying from $50$ to $5,000$ and a validation set of size $\nval=10^6$. The actual optimal rank, denoted by $\Ropt$, is identified as the one yielding the minimum relative generalization error $\errG$, the latter estimated with the validation set as described in Section~\ref{sec:err}. It is of interest to compare $\Ropt$ with the rank selected with the 3-fold CV approach, denoted by $\RCV$, and the rank selected with the LOO-error criterion, denoted by $\RLOO$. Note that for a given ED, the value of $\RCV$ is not fixed because of the random partition of the ED into training and testing sets. This aspect of randomness is considered in the following analysis.

\begin{figure}[!ht]
	\centering
	\includegraphics[width=0.45\textwidth]{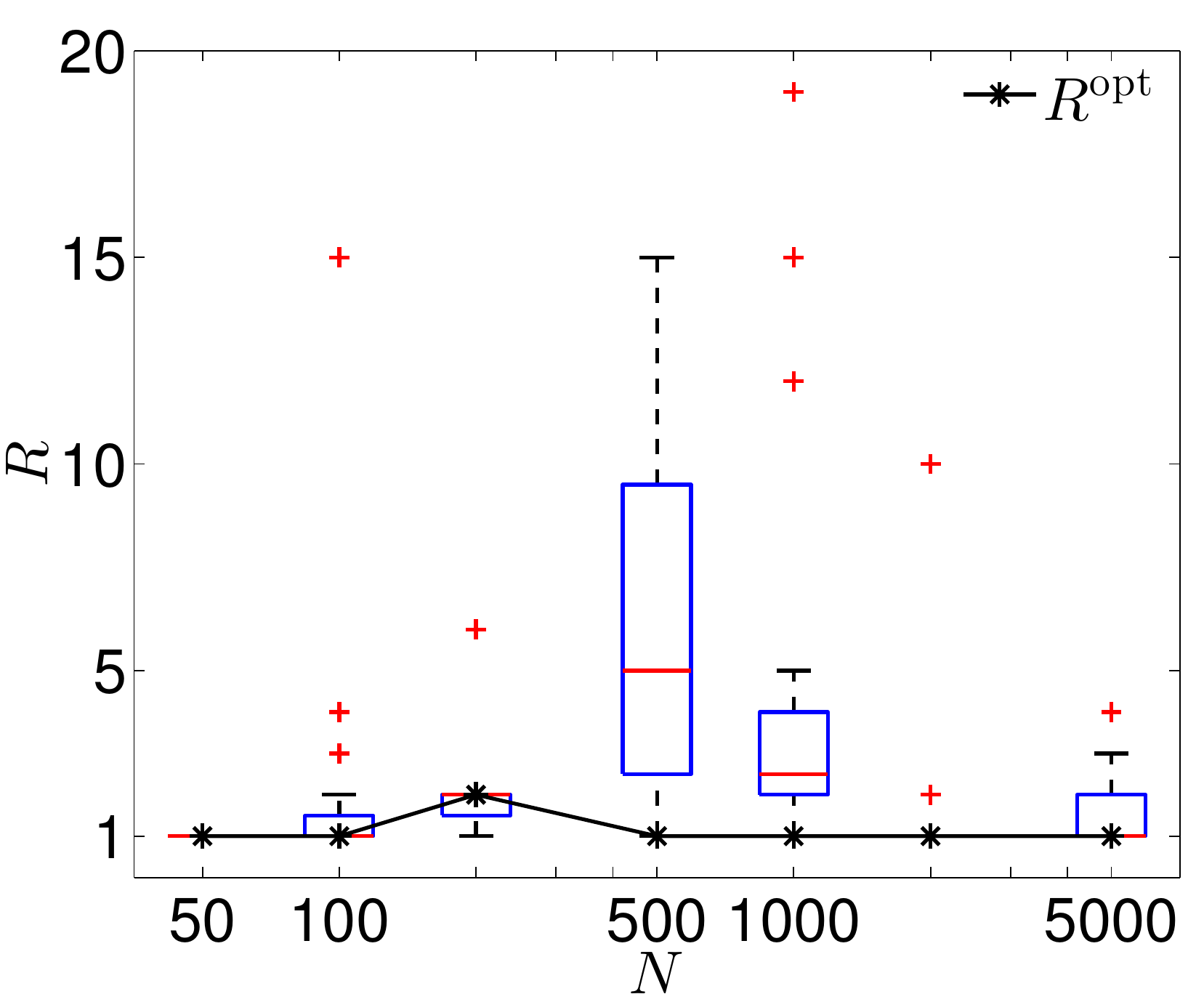}
	\includegraphics[width=0.48\textwidth]{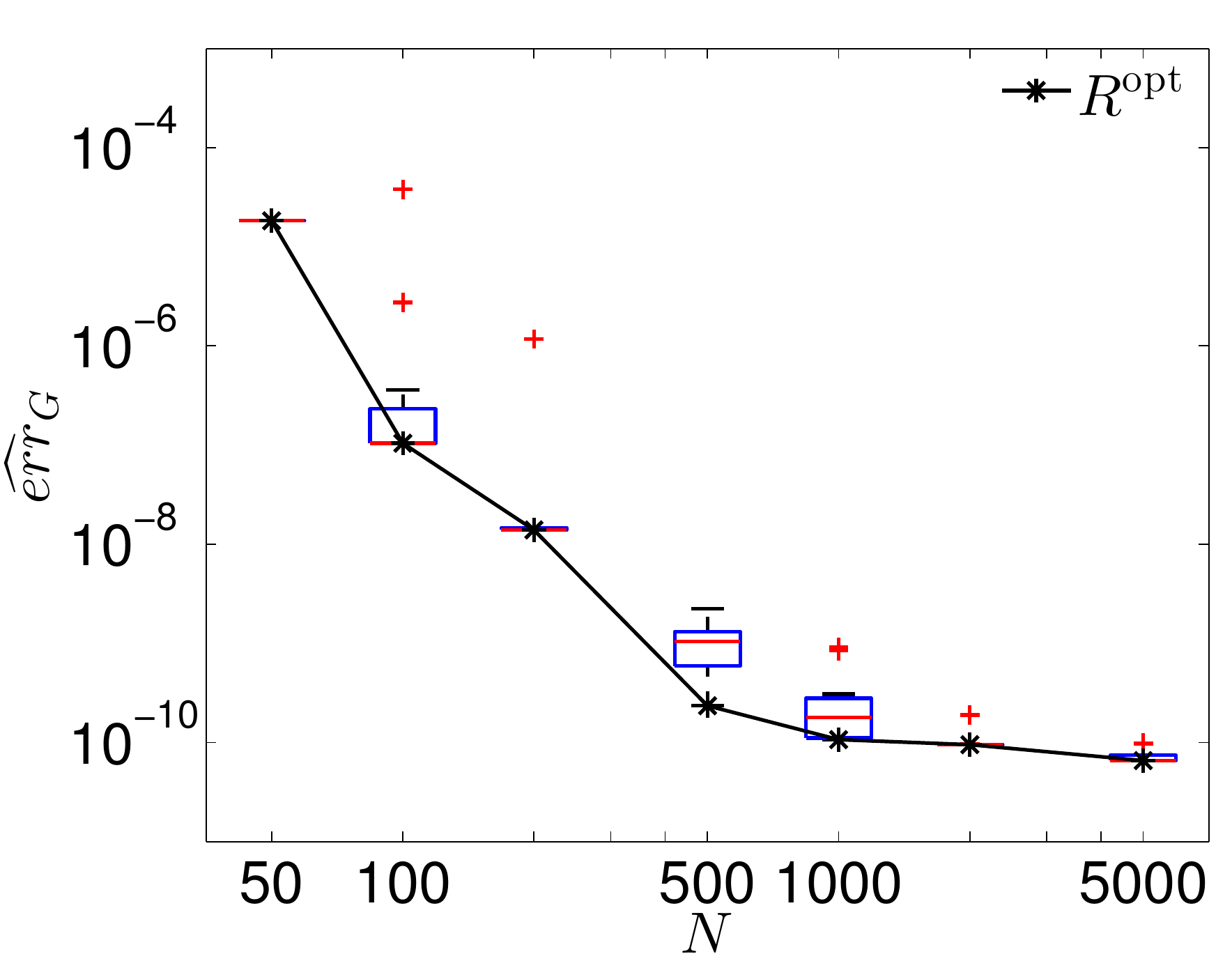}
	\caption{Beam-deflection problem: Comparison of ranks selected with 3-fold CV (20 replications) to optimal ranks (left) and corresponding relative generalization errors (right).}
	\label{fig:Beam_rank}
\end{figure}

In Figure~\ref{fig:Beam_rank}, we investigate rank selection with the 3-fold CV approach. The left graph shows boxplots of $\RCV$ for 20 random partitions of each considered ED into training and testing sets, while the right graph shows boxplots of the relative generalization errors of the resulting meta-models. In order to assess the accuracy of the 3-fold CV approach, we also plot the actual optimal rank $\Ropt$ and corresponding relative generalization error for each ED. The rank-1 structure of the original model (see Eq.~(\ref{eq:beam_u})) matches the structure of the optimal LRA for all $N$ except for $N=100$ ($\Ropt=2$ in this case). The median $\RCV$ coincides with $\Ropt$ for most EDs, while selection of a non-optimal rank appears to have a smaller effect for the larger EDs. We underline that highly accurate meta-models are obtained in all cases, with $\errG$ not exceeding $10^{-4}$. For $N\geq200$, $\errG$ is of the order of $10^{-8}$ or smaller even for cases with non-optimal ranks. The LOO-error criterion yields $\RLOO=1$ for all considered EDs, which coincides with $\Ropt$ except for $N=100$.

\begin{figure}[!ht]
	\centering
	\includegraphics[width=0.48\textwidth]{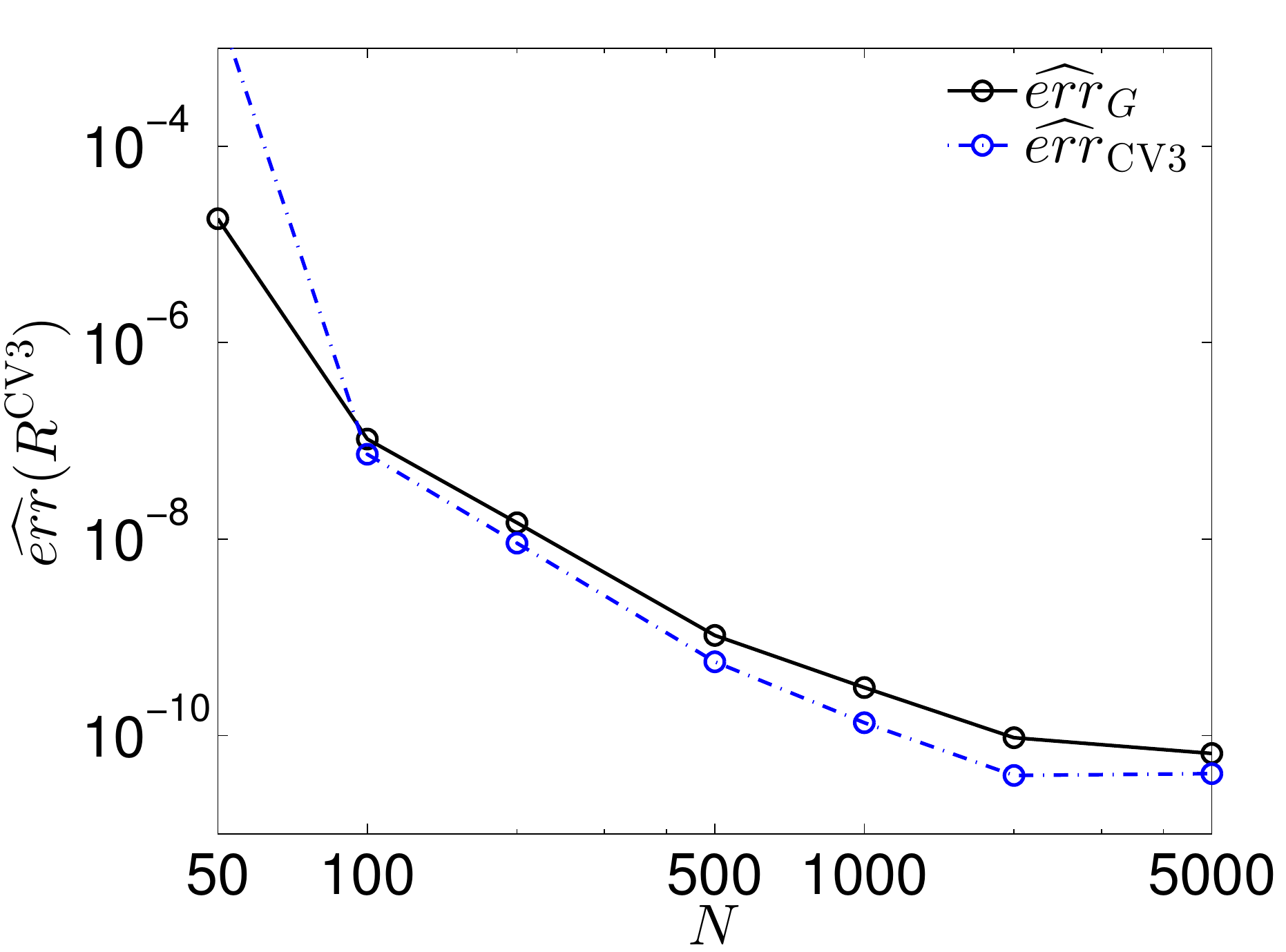}
	\includegraphics[width=0.48\textwidth]{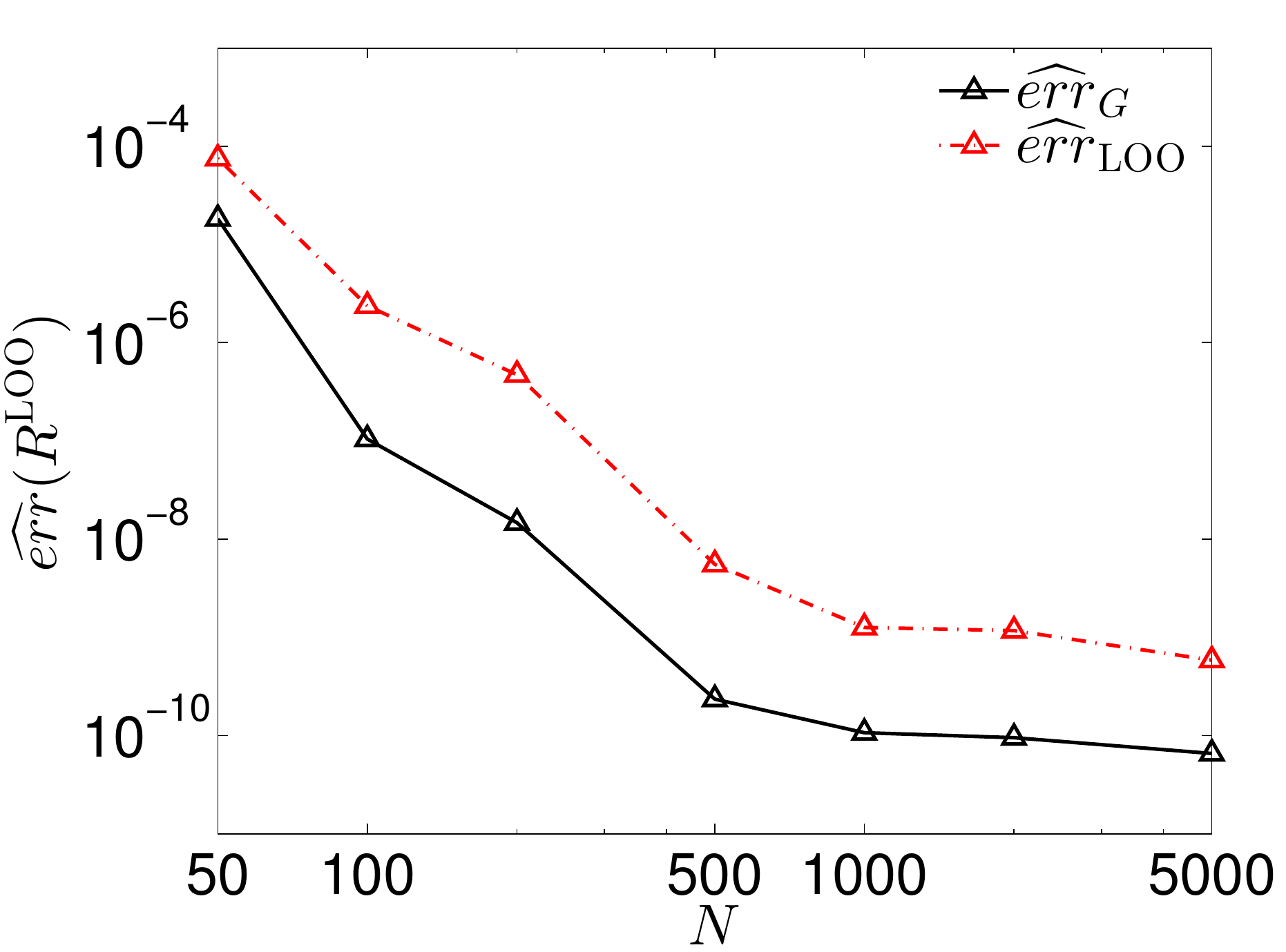}
	\caption{Beam-deflection problem: Comparison of ED-based errors to corresponding relative generalization errors for ranks selected with 3-fold CV (left) and with the LOO-error criterion (right).}
	\label{fig:Beam_rank_2}
\end{figure}

To gain further insight into the criteria for rank selection, we examine the estimation of the generalization error by the ED-based error measures. In the left graph of Figure~\ref{fig:Beam_rank_2}, we consider meta-models with rank $\RCV$ and compare the 3-fold CV error $\errCV$ with the corresponding generalization error $\errG$ for one example partition of each ED. In the right graph of the same figure, we consider meta-models with rank $\RLOO$ and compare the LOO error $\errLOO$ with the corresponding generalization error. In these graphs, the 3-fold CV error appears to be a better estimator of the generalization error compared to the LOO error. It must be noted that contrary to the 3-fold CV error, the estimation of the LOO error largely deteriorates for higher ranks. The lack of orthogonality of the regressors leads to overly high correction factors, whereas the non-corrected errors are too optimistic. Overall, $\errLOO$ is deemed inappropriate for estimating $\errG$ in the present example.

\subsubsection{Stopping criterion in the correction step}
\label{sec:Beam_stopp}

\begin{figure}[!ht]
	\centering
	\includegraphics[width=0.48\textwidth]{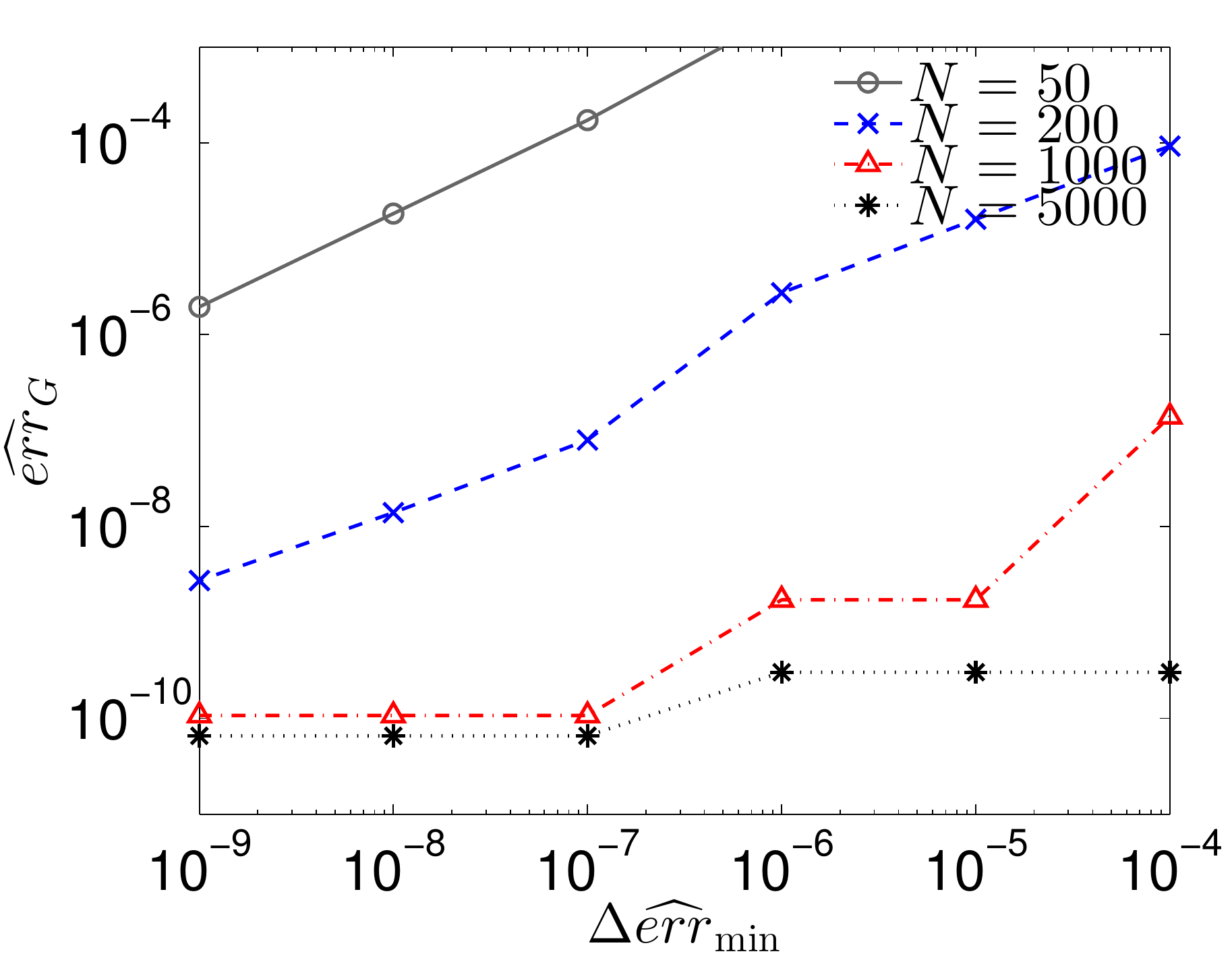}
	\includegraphics[width=0.45\textwidth]{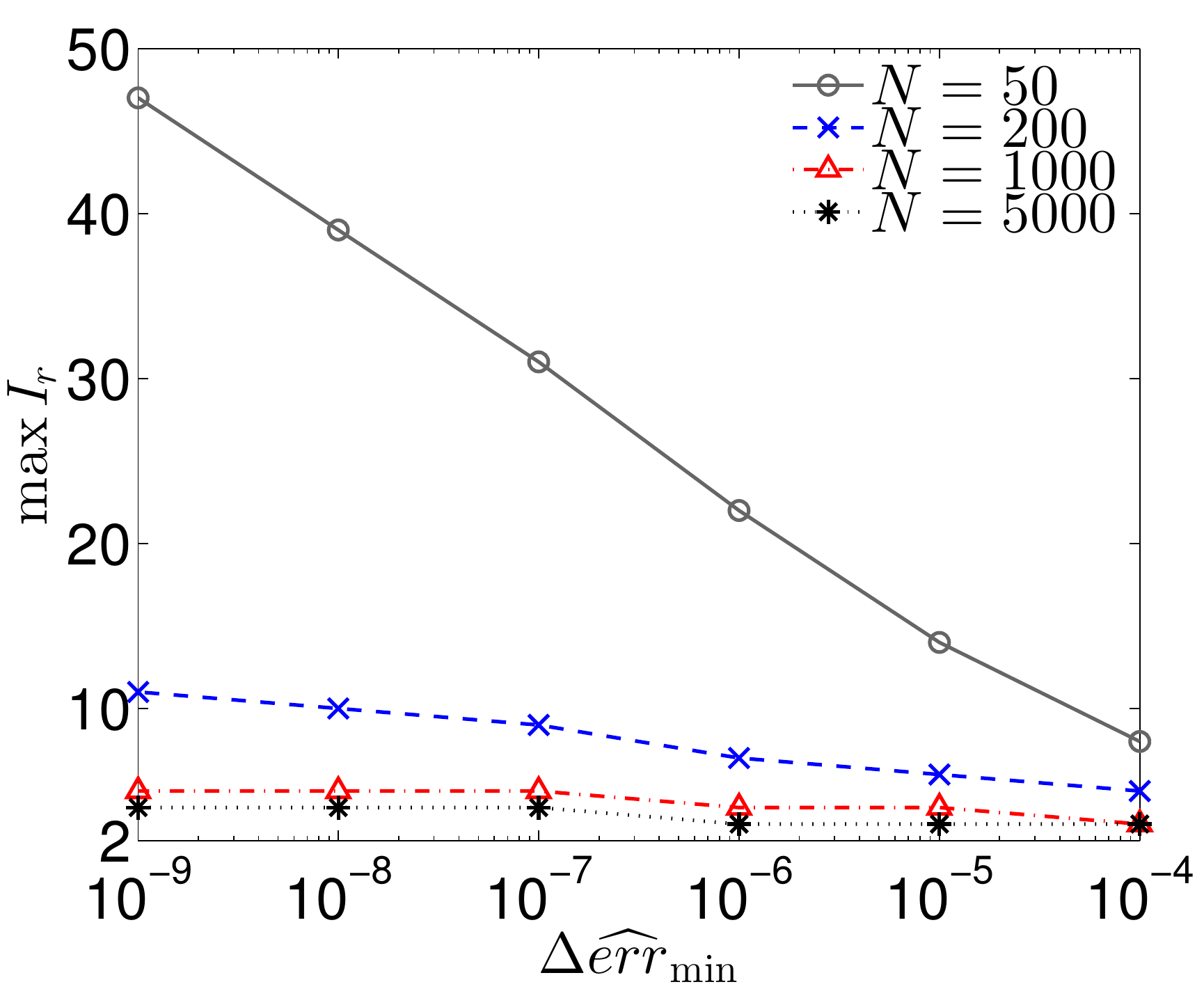}
	\caption{Beam-deflection problem: Relative generalization error (left) and corresponding maximum number of iterations in a correction step (right) versus threshold of differential error in the stopping criterion.}
	\label{fig:Beam_Derr}
\end{figure}

We herein investigate the effects of different stopping criteria in the correction step on the LRA accuracy. For $N\in\{50; 200; 1,000; 5,000\}$, the left graph of Figure~\ref{fig:Beam_Derr} shows the relative generalization errors of the LRA meta-models with optimal rank, while the parameter $\Derrmin$ that drives the loops over the set of dimensions in ALS varies between $10^{-9}$ and $10^{-4}$; other parameters are fixed to their values in Section~\ref{sec:Beam_rank}. The right graph of the same figure indicates the maximum number of iterations $I_r$ performed in a correction step. Except for the largest ED, the accuracy of LRA strongly depends on $\Derrmin$, with decreasing values leading to orders-of-magnitude smaller $\errG$. The right graph indicates that a larger number of iterations is required for the smallest ED.

\subsection{Truss deflection}
\label{sec:Truss_LRA}

\begin{figure}[!ht]
	\centering
	\includegraphics[trim = 15mm 70mm 15mm 70mm, width=0.6\textwidth]{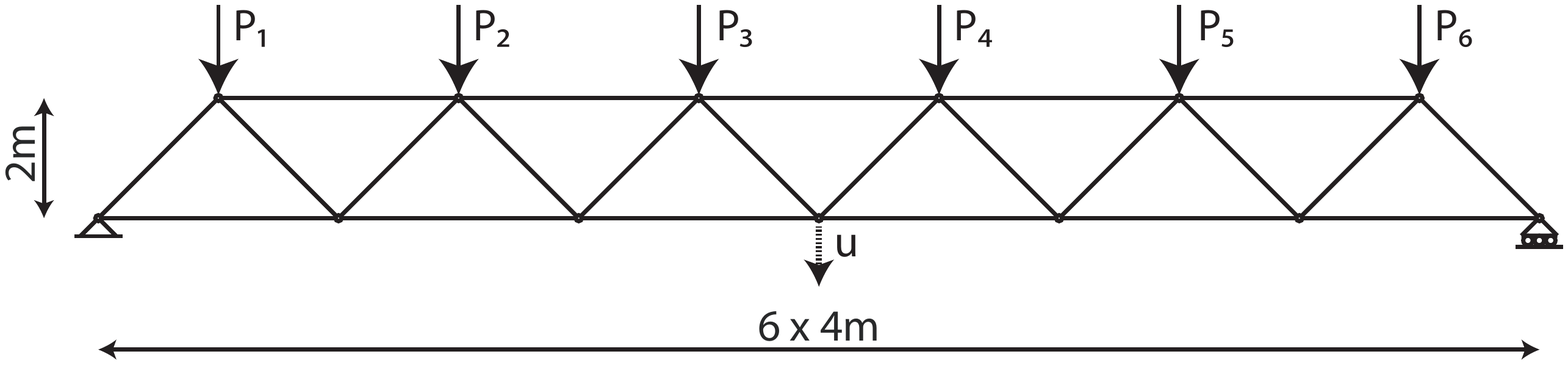}
	\caption{Truss structure.}
	\label{fig:truss}
\end{figure}

The second example also derives from structural mechanics, but involves a finite-element model. In particular, the truss shown in Figure~\ref{fig:truss} is considered (also studied in \citet{BlatmanCras2008}), with the mid-span deflection $u$ representing the response quantity of interest. The random input herein comprises $M=10$ independent variables: the vertical loads $P_1 \enum P_6$, the cross-sectional area and Young's modulus of the horizontal bars, respectively denoted by $A_1$ and $E_1$, and the cross-sectional area and Young's modulus of the diagonal bars, respectively denoted by $A_2$ and $E_2$. The distributions of the input variables are listed in Table~\ref{tab:truss_input}. After employing an isoprobabilistic transformation of $\ve X=\{P_1 \enum P_6,A_1,A_2,E_1,E_2\}$ into a vector of standard normal variables, we develop LRA representation of the random response $U=\cm(\ve X)$ with basis functions made of Hermite polynomials. In the underlying deterministic problem, the mid-span deflection is computed with an in-house finite-element analysis code developed in the Matlab environment.

\begin{table} [!ht]
	\centering
	\caption{Truss-deflection problem: Distributions of input random variables.}
	\label{tab:truss_input}
	\begin{tabular}{c c c c}
		\hline Variable & Distribution & Mean & CoV \\
		\hline $A_1~[\rm m]$ & Lognormal & $0.002$ & $0.10$ \\
		$A_2~[\rm m]$ & Lognormal & $0.001$  & $0.10$ \\
		$E_1, E_2~[\rm MPa]$ & Lognormal & $2.1 \cdot 10^5$   & $0.10$\\
		$P_1\enum P_6~[\rm KN]$ & Gumbel & $50$  & $0.15$ \\
		\hline       
	\end{tabular}
\end{table}

\subsubsection{Rank selection and error measures}

Similarly to the previous example, we investigate rank selection based on the 3-fold CV and LOO errors, considering the candidate values $\{1 \enum 20\}$. After preliminary investigations, the polynomial degree is set to $p_i=p=3$ for $i=1 \enum 10$, while the parameters of the stopping criterion in the correction step are set to $\Imax=50$ and $\Derrmin= 10^{-6}$. Again, we use EDs of size $N$ varying between $50$ and $5,000$ to build the LRA meta-models and a validation set of size $\nval=10^6$ to estimate the relative generalization errors and identify the actual optimal ranks.

\begin{figure}[!ht]
	\centering
	\includegraphics[width=0.46\textwidth]{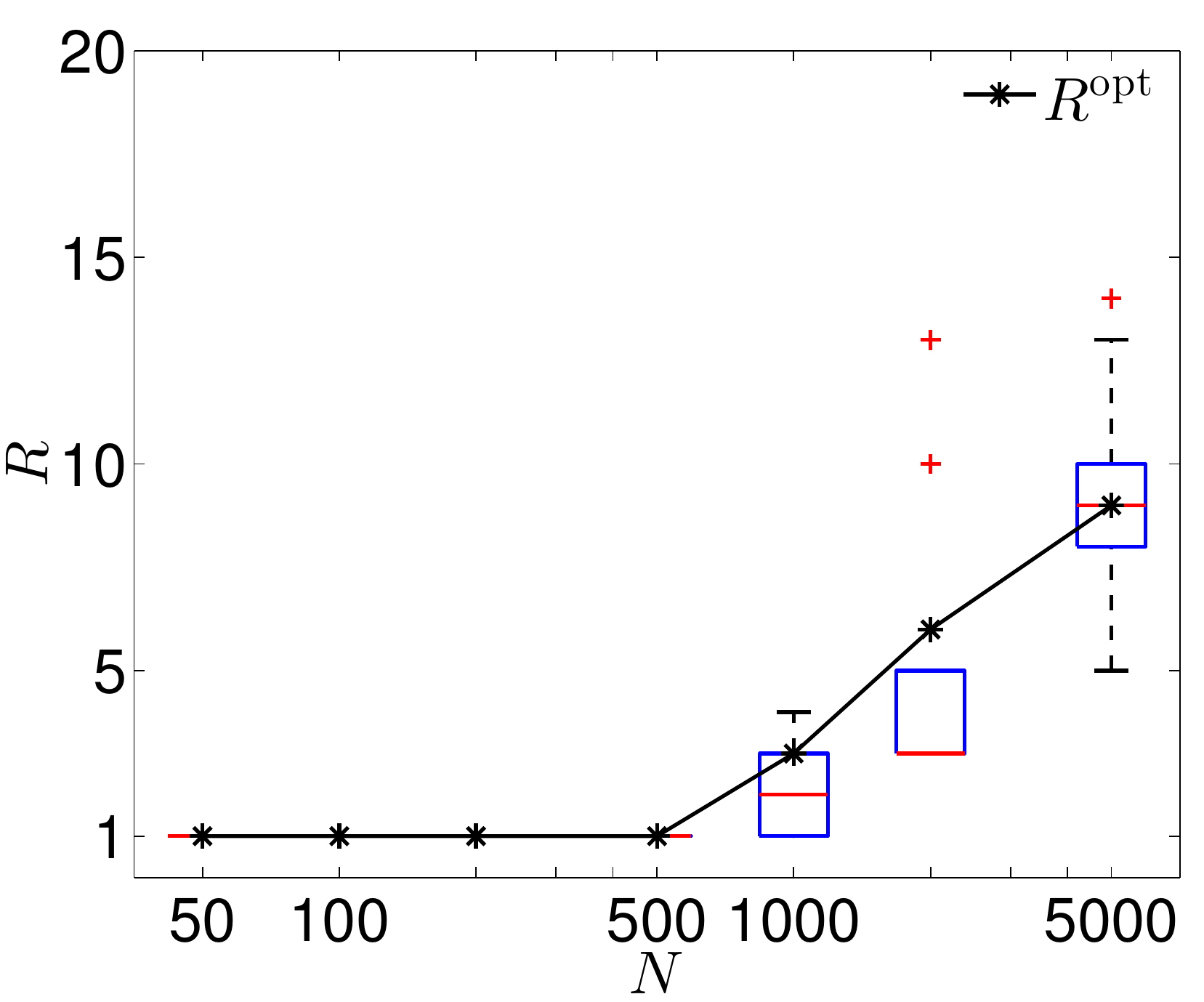}
	\includegraphics[width=0.48\textwidth]{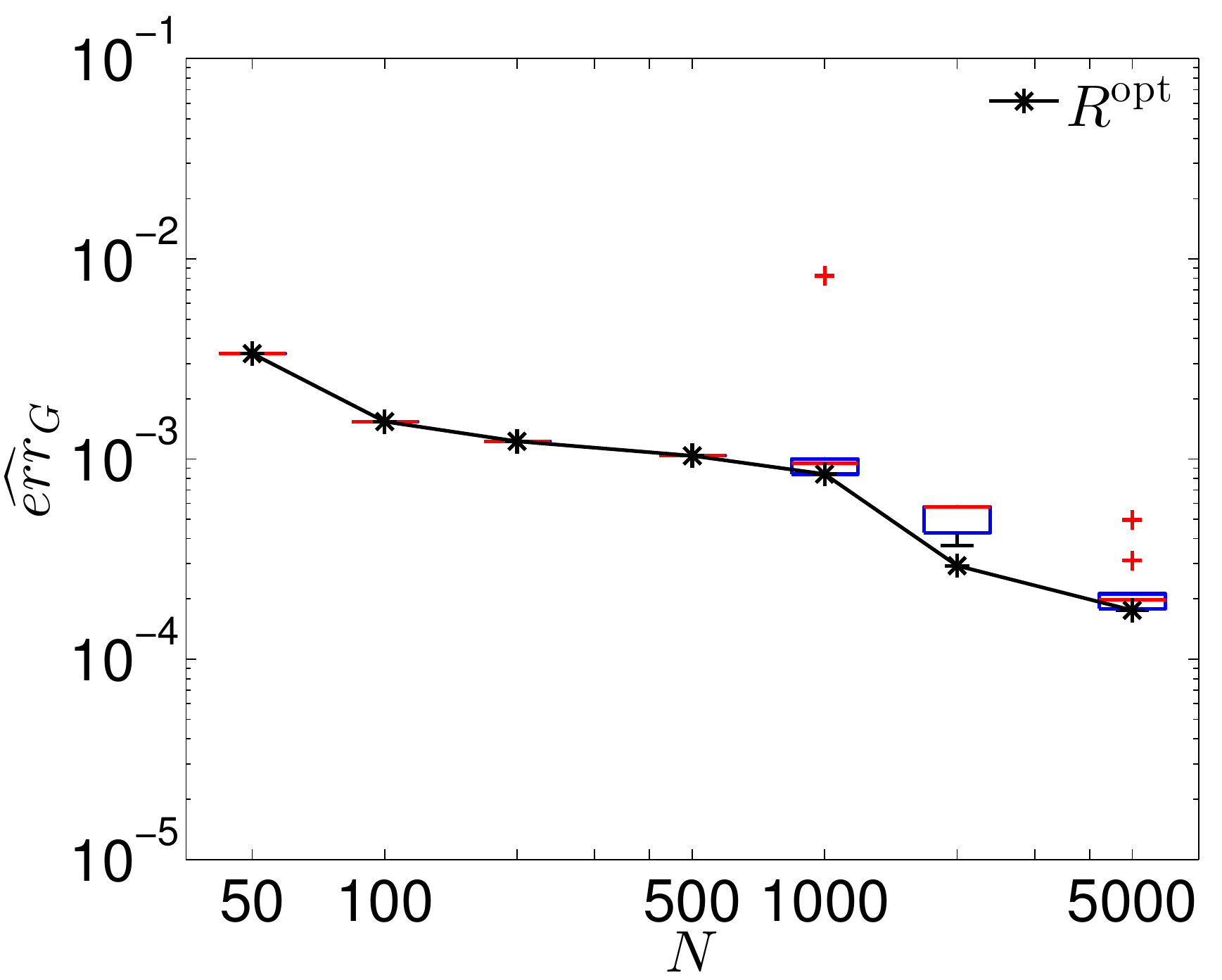}
	\caption{Truss-deflection problem: Comparison of ranks selected with 3-fold CV (20 replications) to optimal ranks (left) and corresponding relative generalization errors (right).}
	\label{fig:Truss_rank}
\end{figure}

The left graph of Figure~\ref{fig:Truss_rank} shows boxplots of $\RCV$, \ie the rank selected with 3-fold CV, for 20 random partitions of each considered ED, together with the corresponding actual optimal rank $\Ropt$. The right graph shows boxplots of the relative generalization errors of the meta-models with rank $\RCV$ as well as the relative generalization errors of the meta-models with rank $\Ropt$. We observe that the actual optimal rank is equal to 1 for $N\leq500$, but attains higher values for larger EDs. This is because the larger amount of information contained in the latter allows the estimation of a larger number of coefficients with higher accuracy.  This trend is captured by the 3-fold CV approach, which yields $\RCV=1$ for all random partitions of the ED when $N\leq500$ and varying higher ranks for the larger EDs. The median of $\RCV$ coincides with $\Ropt$ for $N=5,000$, but is smaller than $\RCV$ for $N=1,000$ and $N=2,000$. As in the beam-deflection problem, effects of non-optimal rank selection on the LRA accuracy tend to be less significant as the ED size increases. The LOO-error criterion yields $\RLOO=1$ for all considered EDs, which coincides with $\Ropt$ only for $N\leq500$.

\begin{figure}[!ht]
	\centering
	\includegraphics[width=0.48\textwidth]{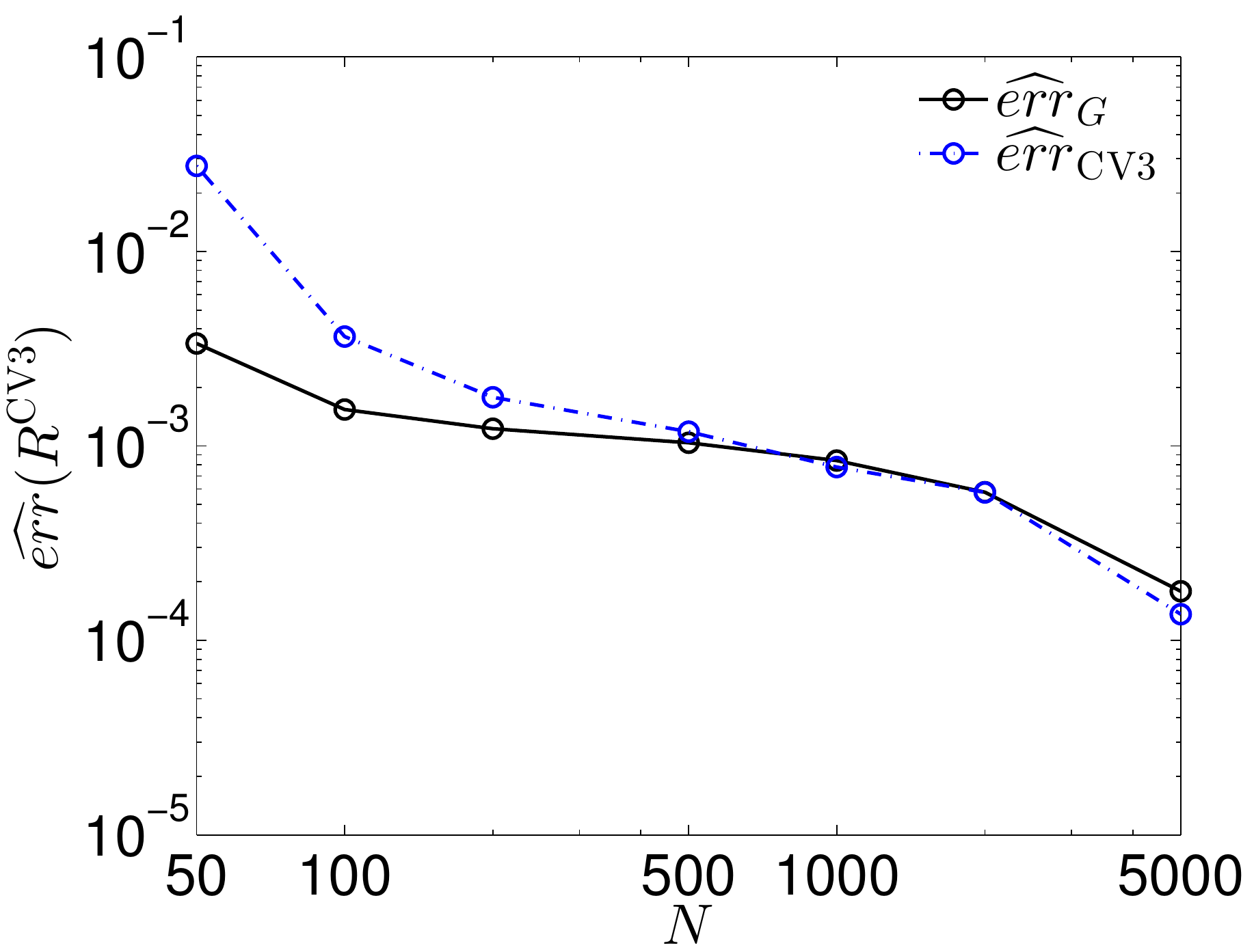}
	\includegraphics[width=0.48\textwidth]{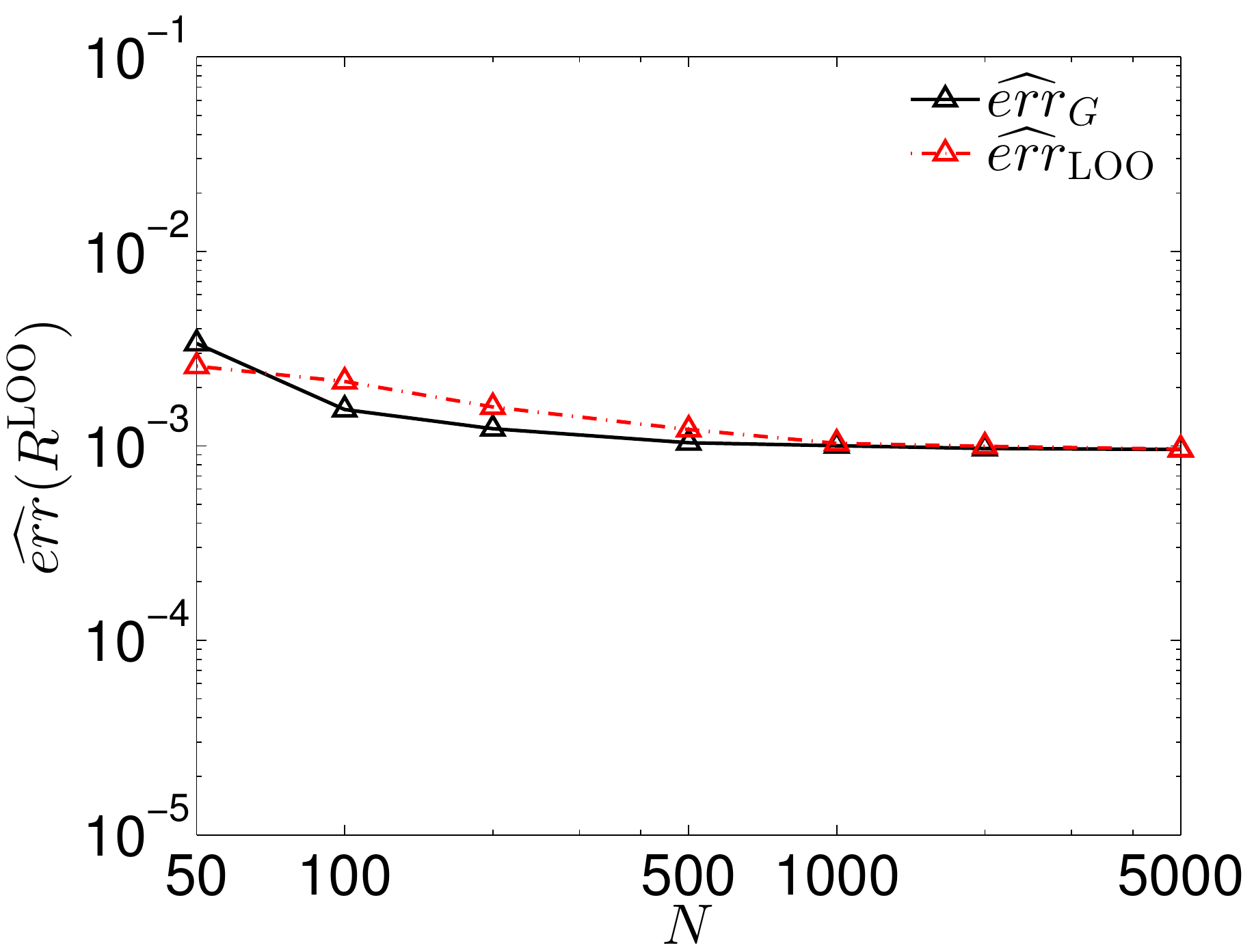}
	\caption{Truss-deflection problem: Comparison of ED-based errors to corresponding relative generalization errors for ranks selected with 3-fold CV (left) and with the LOO-error criterion (right).}
	\label{fig:Truss_rank_2}
\end{figure}

We next assess the estimation of the generalization error by the ED-based error measures. In the left graph of Figure~\ref{fig:Truss_rank_2}, we compare $\errCV$ to $\errG$ for LRA meta-models with rank $\RCV$. In the right graph, we compare $\errLOO$ to $\errG$ for LRA meta-models with rank $\RLOO$. Both ED-based error measures appear to approximate the corresponding generalization errors fairly well, particularly for the larger EDs. However, similarly to the beam-deflection problem, $\errLOO$ is found to be an inappropriate estimator of $\errG$ for higher-rank meta-models. On the other hand, for sufficiently large EDs, $\errCV$ provides fair estimates of $\errG$ for all ranks.

\subsubsection{Stopping criterion in the correction step}

\begin{figure}[!ht]
	\centering
	\includegraphics[width=0.48\textwidth]{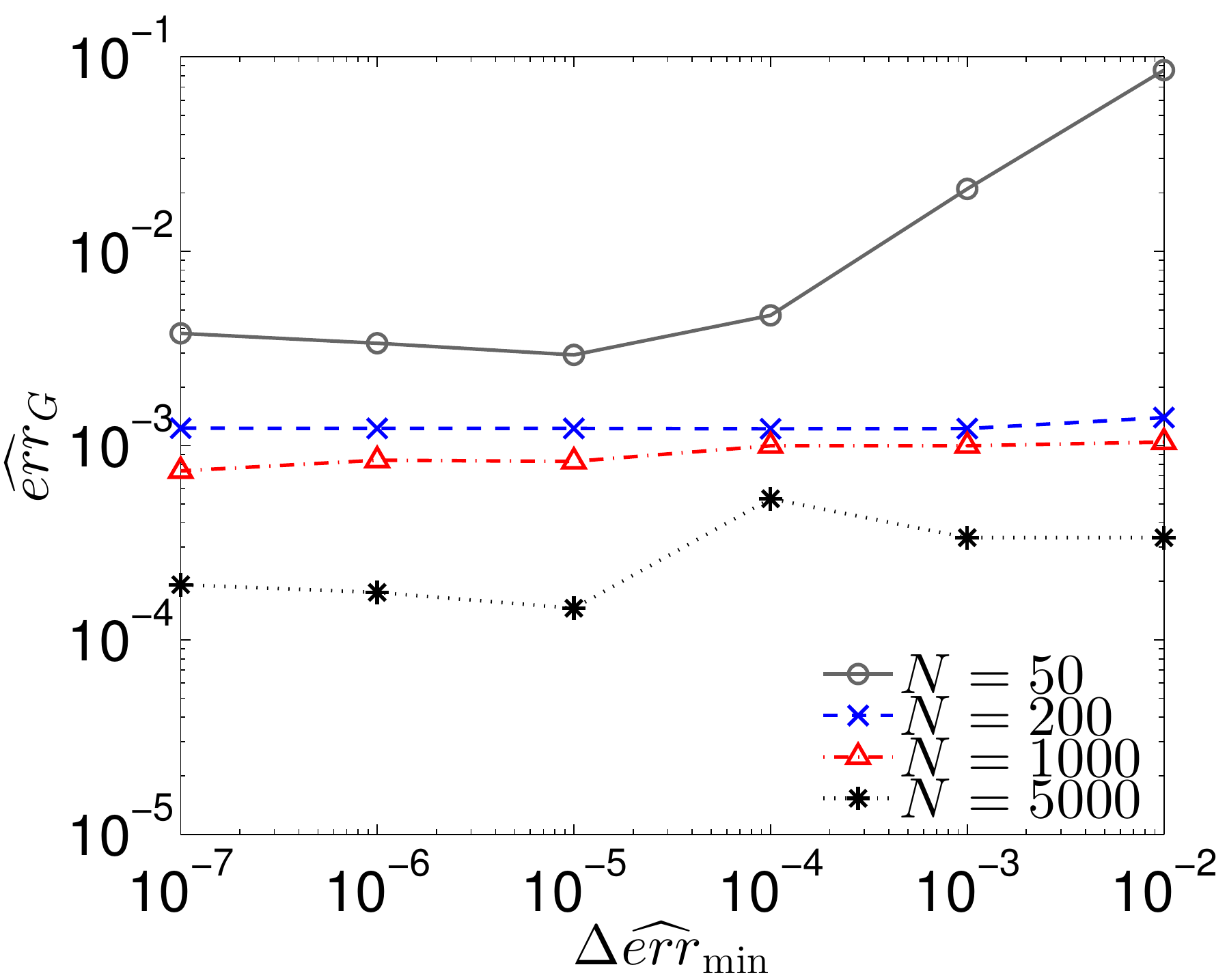}
	\includegraphics[width=0.46\textwidth]{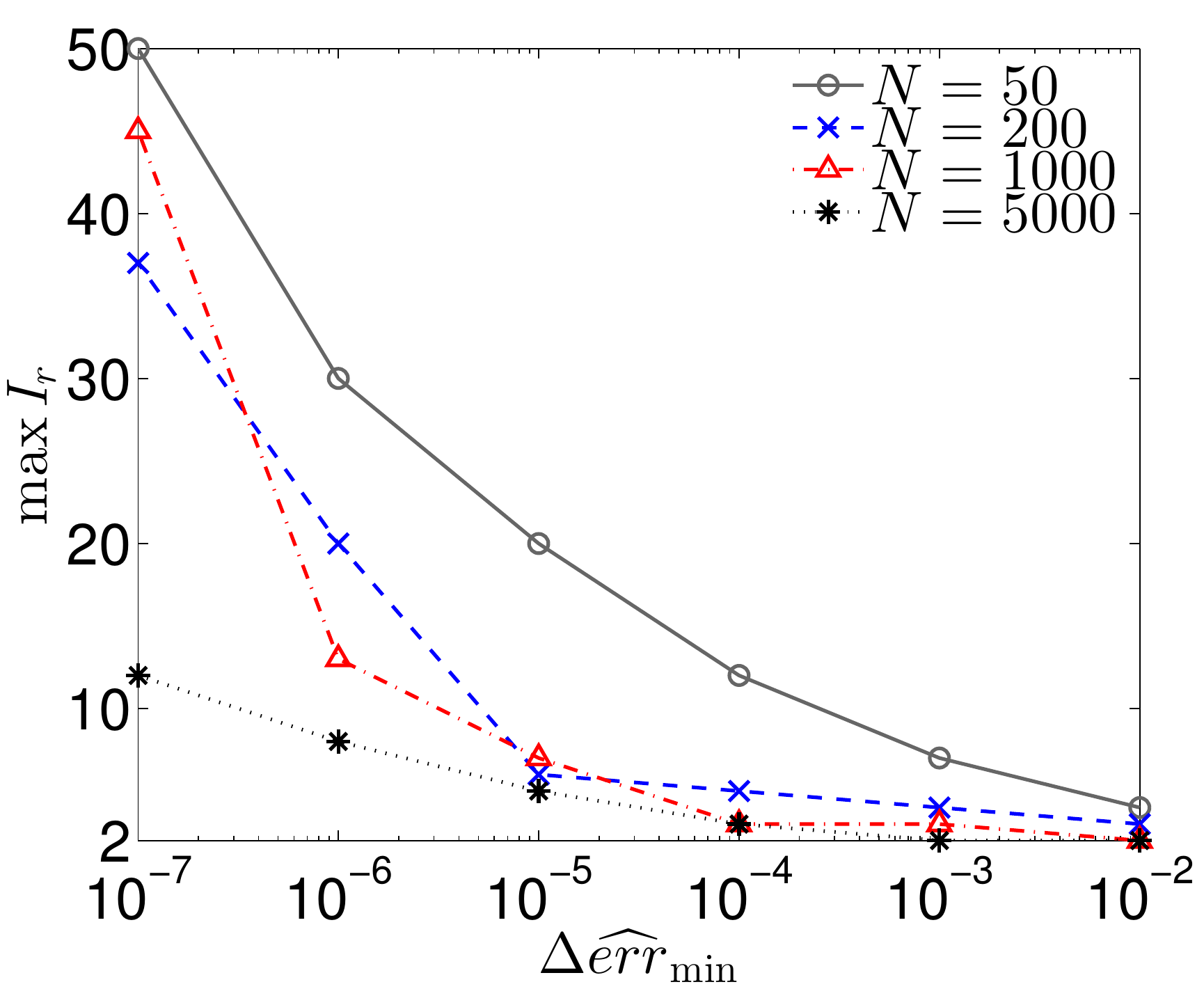}
	\caption{Truss-deflection problem: Relative generalization error (left) and corresponding maximum number of iterations in a correction step (right) versus threshold of differential error in the stopping criterion.}
	\label{fig:Truss_Derr}
\end{figure}

We now examine effects of the differential error threshold in the correction step, while other parameters are fixed to their values above. For $N\in\{50; 200; 1,000; 5,000\}$, the left graph of Figure~\ref{fig:Truss_Derr} shows the relative generalization errors for the LRA meta-models with optimal ranks, while $\Derrmin$ varies from $10^{-7}$ to $10^{-2}$. The right graph of the same figure shows the corresponding maximum number of iterations $I_r$ in a correction step. While selecting a sufficiently small threshold $\Derrmin$ appears critical for the smaller ED, it does not essentially affect the meta-model accuracy when $N=200$ or $N=1,000$. For $N=5,000$, variations of $\errG$ within the same order of magnitude are observed with decreasing $\Derrmin$, following a decreasing but non-monotonic trend. Thus, small values of $\Derrmin$ may lead to unnecessary iterations in the correction step for the larger EDs.

\subsection{Frame displacement}
\label{sec:Frame_LRA}

\begin{figure}[!ht]
	\centering
	\includegraphics[trim = 0mm 0mm 0mm 0mm, width=0.7\textwidth]{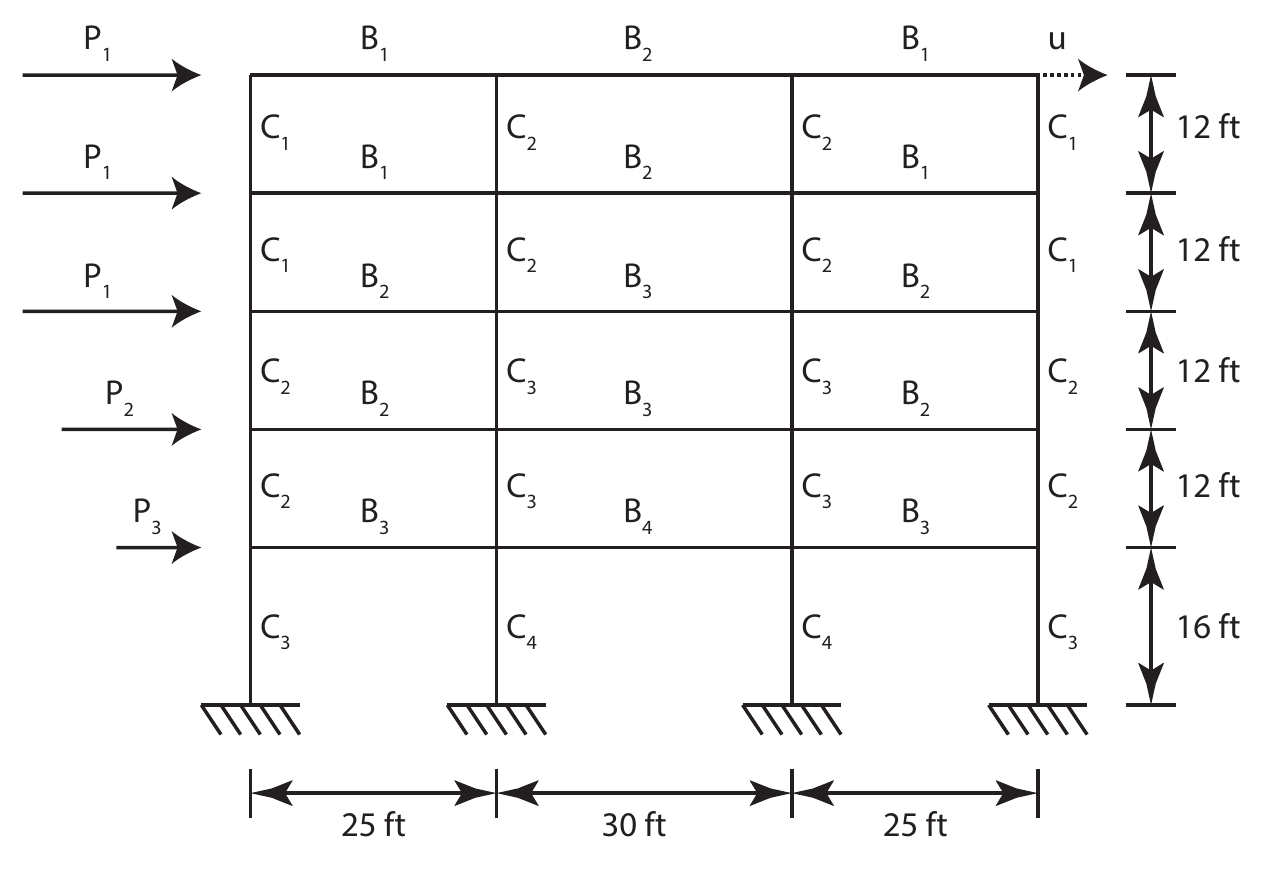}
	\caption{Frame structure.}
	\label{fig:frame}
\end{figure}

The present example involves a finite-element model representing the three-span five-story frame shown in Figure~\ref{fig:frame} (also studied in \citet{LiuPL91,BlatmanPEM2010}). The response quantity of interest is the horizontal displacement $u$ at the top right corner of the top floor, under the depicted horizontal loads acting at the floor levels. Table~\ref{tab:frame_properties} lists the properties (Young's modulus, moment of inertia, cross-sectional area) of the different elements according to their labels in Figure~\ref{fig:frame}. These element properties together with the values of the horizontal loads comprise the random input of the problem of dimension $M=21$. The distributions of the input random variables are listed in Table~\ref{tab:frame_input}. The truncation of the Gaussian distributions used to model the element properties is a modification of the original example in \citet{LiuPL91} that was introduced by \citet{BlatmanPEM2010}. The input variables are correlated, with their dependence structure defined using a Gaussian copula (see \citet{Nelsen2006} for modeling of probabilistic dependence with copulas). The elements of the associated linear correlation matrix are defined as follows:
\begin{itemize}
	\item the correlation coefficient between the cross-sectional area $A_i$ and the moment of inertia $I_i$ of a certain element $i$ is $\rho_{A_i,I_i}=0.95$;
	\item the correlation coefficient between the geometric properties of two distinct elements $i$ and $j$ are $\rho_{A_i,I_j}=\rho_{I_i,I_j}=\rho_{A_i,A_j}=0.13$;
	\item the correlation coefficient  between the two Young's moduli is $\rho_{E_1,E_2}=0.90$;
	\item the remaining correlation coefficients are zero.
\end{itemize}
Note that in the original example, the above values are considered for the corresponding linear correlation coefficients in the standard normal space; however, as noted in \citet{BlatmanPEM2010}, the differences between the two are insignificant. After an isoprobabilistic transformation of $\ve X=\{P_1,P_2,P_3,E_1,E_2,I_1\enum I_8,A_1 \enum A_8\}$ into a vector of independent standard normal variables, we develop LRA representations of the random response $U=\cm(\ve X)$ with the basis functions made of Hermite polynomials.  In the underlying deterministic problem, the displacement of interest is computed with an in-house finite-element analysis code developed in the Matlab environment.

\begin{table} [!ht]
	\centering
	\caption{Frame-displacement problem: Element properties.}
	\label{tab:frame_properties}
	\begin{tabular}{c c c c}
		Element & Young's modulus & Moment of inertia & Cross-sectional area\\
		\hline
		$C_1$ & $E_1$ & $I_1$ & $A_1$ \\
		$C_2$ & $E_1$ & $I_2$ & $A_2$\\
		$C_3$ & $E_1$ & $I_3$ & $A_3$\\
		$C_4$ & $E_1$ & $I_4$ & $A_4$\\
		$B_1$ & $E_2$ & $I_5$ & $A_5$ \\
		$B_2$ & $E_2$ & $I_6$ & $A_6$\\
		$B_3$ & $E_2$ & $I_7$ & $A_7$\\
		$B_4$ & $E_2$ & $I_8$ & $A_8$\\
		\hline       
	\end{tabular}
\end{table}

\begin{table} [!ht]
	\centering
	\caption{Frame-displacement problem: Distributions of input random variables.}
	\label{tab:frame_input}
	\begin{tabular}{c c c c}
		\hline Variable & Distribution & Mean & Standard deviation \\
		\hline
		$P_1~[\rm KN]$ & Lognormal & $133.45$ & $40.04$ \\          
		$P_2~[\rm KN]$ & Lognormal & $88.97$ & $35.59$ \\           
		$P_3~[\rm KN]$ & Lognormal & $71.17$ & $28.47$ \\
		$E_1~[\rm KN/m^2]$ & Truncated Gaussian over $[0,\infty)$ & $2.3796\cdot 10^{7}$ & $1.9152\cdot 10^{6}$ \\   
		$E_2~[\rm KN/m^2]$ & Truncated Gaussian over $[0,\infty)$ & $2.1738\cdot 10^{7}$ & $1.9152\cdot 10^{6}$ \\
		$I_1~[\rm m^4]$ & Truncated Gaussian over $[0,\infty)$ & $8.1344\cdot 10^{-3}$ & $1.0834\cdot 10^{-3}$ \\
		$I_2~[\rm m^4]$ & Truncated Gaussian over $[0,\infty)$ & $1.1509\cdot 10^{-2}$ & $1.2980\cdot 10^{-3}$ \\
		$I_3~[\rm m^4]$ & Truncated Gaussian over $[0,\infty)$ & $2.1375\cdot 10^{-2}$ & $2.5961e-03$ \\
		$I_4~[\rm m^4]$ & Truncated Gaussian over $[0,\infty)$ & $2.5961\cdot 10^{-2}$ &$3.0288\cdot 10^{-3}$ \\
		$I_5~[\rm m^4]$ & Truncated Gaussian over $[0,\infty)$ & $1.0811\cdot 10^{-2}$ & $2.5961\cdot 10^{-3}$ \\
		$I_6~[\rm m^4]$ & Truncated Gaussian over $[0,\infty)$ & $1.4105\cdot 10^{-2}$ & $3.4615\cdot 10^{-3}$ \\
		$I_7~[\rm m^4]$ & Truncated Gaussian over $[0,\infty)$ & $2.3279\cdot 10^{-2}$ & $5.6249\cdot 10^{-3}$ \\
		$I_8~[\rm m^4]$ & Truncated Gaussian over $[0,\infty)$ & $2.5961\cdot 10^{-2}$ & $6.4902\cdot 10^{-3}$ \\
		$A_1~[\rm m^2]$ & Truncated Gaussian over $[0,\infty)$ & $3.1256\cdot 10^{-1}$ & $5.5815\cdot 10^{-2}$ \\
		$A_2~[\rm m^2]$ & Truncated Gaussian over $[0,\infty)$ & $3.7210\cdot 10^{-1}$ & $7.4420\cdot 10^{-2}$ \\
		$A_3~[\rm m^2]$ & Truncated Gaussian over $[0,\infty)$ & $5.0606\cdot 10^{-1}$ & $9.3025\cdot 10^{-2}$ \\
		$A_4~[\rm m^2]$ & Truncated Gaussian over $[0,\infty)$ & $5.5815\cdot 10^{-1}$ & $1.1163\cdot 10^{-1}$ \\
		$A_5~[\rm m^2]$ & Truncated Gaussian over $[0,\infty)$ & $2.5302\cdot 10^{-1}$ & $9.3025\cdot 10^{-2}$ \\
		$A_6~[\rm m^2]$ & Truncated Gaussian over $[0,\infty)$ & $2.9117\cdot 10^{-1}$ & $1.0232\cdot 10^{-1}$ \\
		$A_7~[\rm m^2]$ & Truncated Gaussian over $[0,\infty)$ & $3.7303\cdot 10^{-1}$ & $1.2093\cdot 10^{-1}$ \\
		$A_8~[\rm m^2]$ & Truncated Gaussian over $[0,\infty)$ & $4.1860\cdot 10^{-1}$ & $1.9537\cdot 10^{-1}$ \\
		\hline    
	\end{tabular}
\end{table}

\subsubsection{Rank selection and error measures}

We herein investigate rank selection, while setting the polynomial degree to $p_i=p=3$ for $i=1 \enum 21$. We define the stopping criterion in the correction step by requiring $\Imax=50$ and $\Derrmin= 10^{-6}$. As in the two preceding examples, we use EDs of sizes $N$ varying from $50$ to $5,000$ to build the meta-models and a validation set of size $\nval=10^6$ to identify the actual optimal rank among the candidate values $\{1 \enum 20\}$.

\begin{figure}[!ht]
	\centering
	\includegraphics[width=0.46\textwidth]{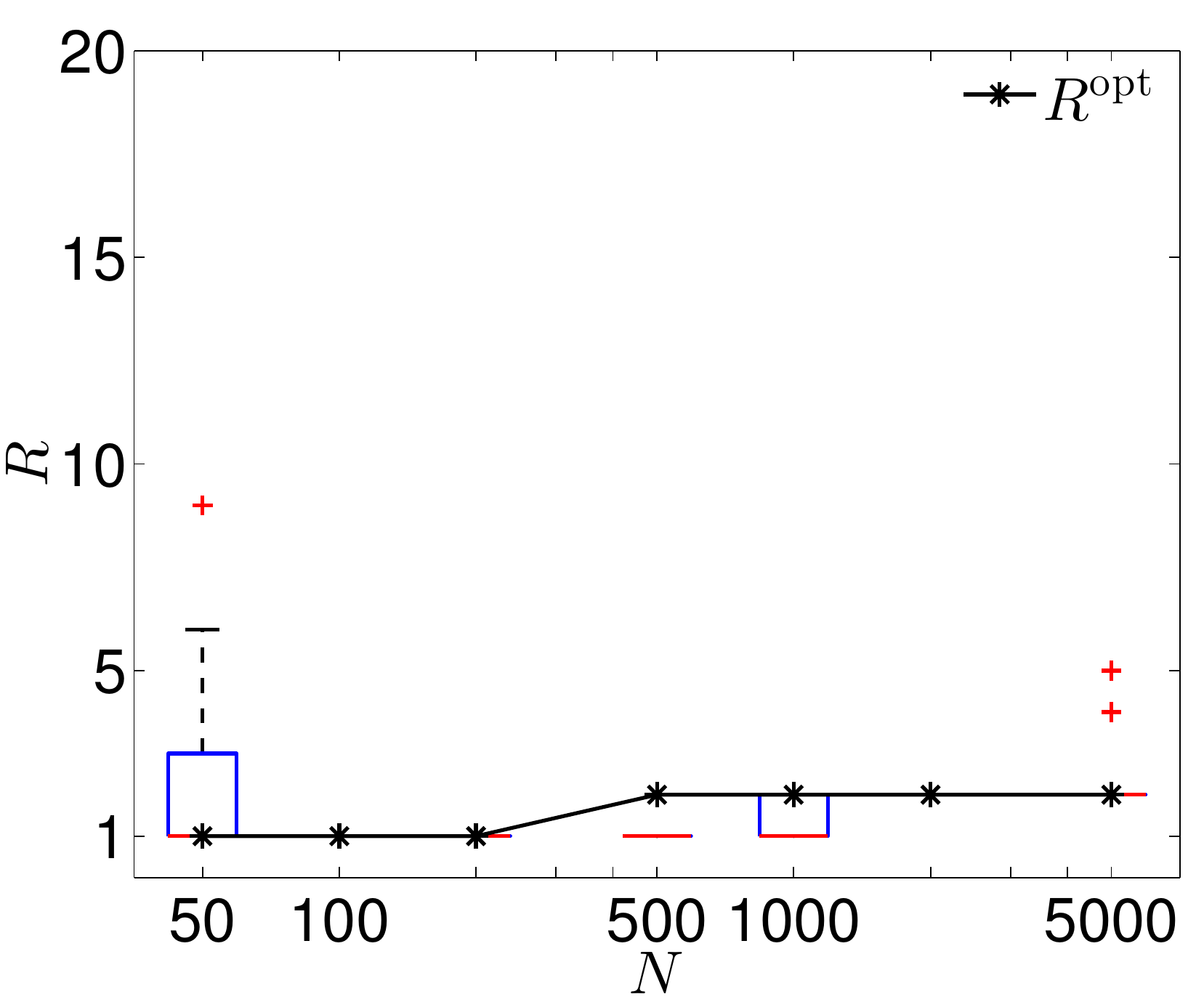}
	\includegraphics[width=0.48\textwidth]{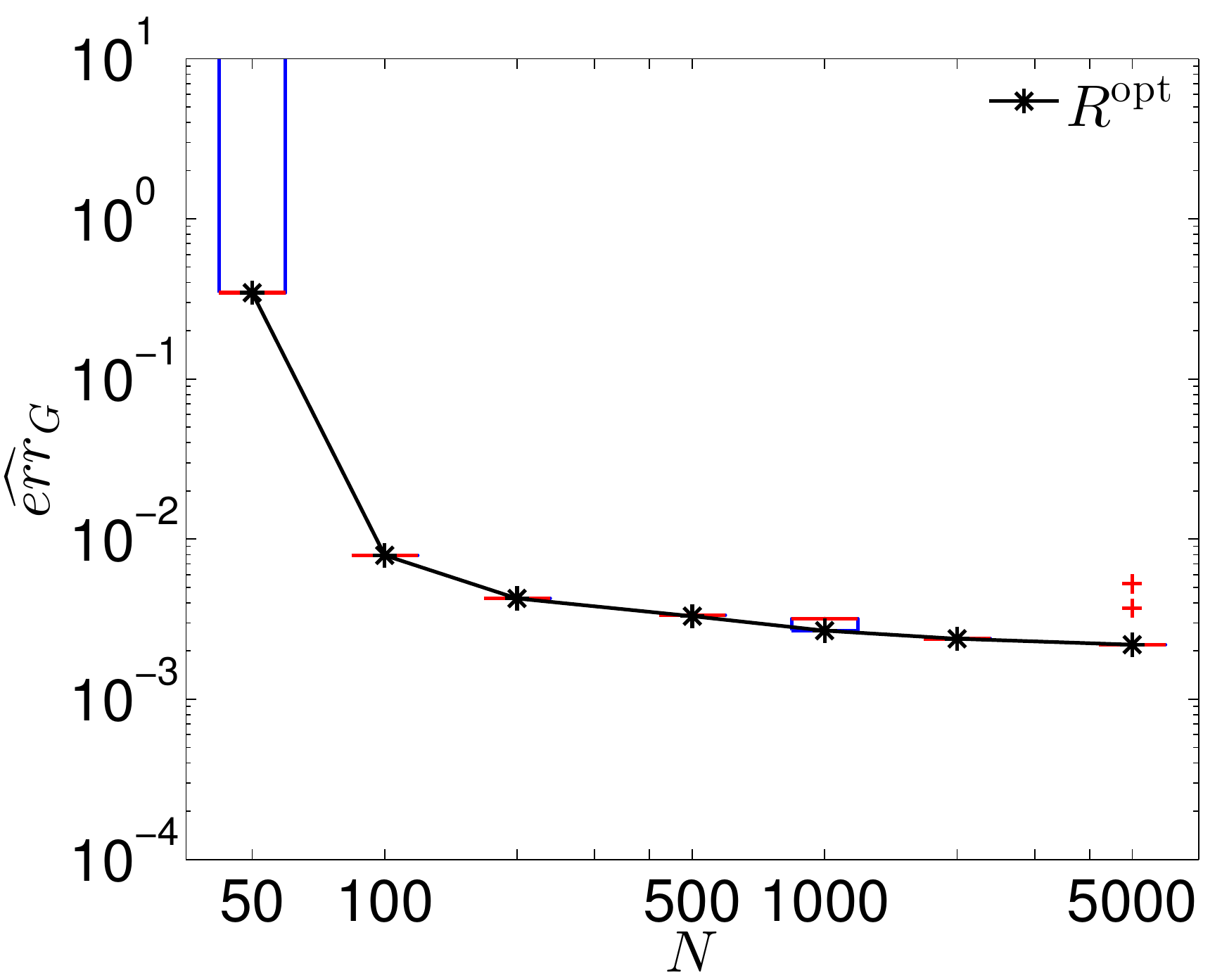}
	\caption{Frame-displacement problem: Comparison of ranks selected with 3-fold CV (20 replications) to optimal ranks (left) and corresponding relative generalization errors (right).}
	\label{fig:Frame_rank}
\end{figure}

Figure~\ref{fig:Frame_rank} presents results of rank selection with the 3-fold CV approach, considering 20 random partitions of each ED. The left and right graphs of this figure respectively show boxplots of $\RCV$ and the relative generalization errors of the corresponding meta-models. The accuracy of the approach is assessed through comparisons to LRA with the actual optimal rank $\Ropt$, identified by means of the validation set. As in the truss-deflection problem, the optimal rank is $\Ropt=1$ for the smaller EDs and increases for larger $N$, without however exceeding $\Ropt=2$ in the present example. Except for $N=50$, the ED-based rank $\RCV$ exhibits small or no dispersion for different random selections of the training and testing sets, while its median value coincides with $\Ropt$ in most cases. Selection of a rank higher than unity may lead to gross generalization errors for $N=50$, because the few points comprising the ED are insufficient to accurately estimate a large number of coefficients. Again, the LOO-error criterion yields $\RLOO=1$ for all considered EDs.

\begin{figure}[!ht]
	\centering
	\includegraphics[width=0.48\textwidth]{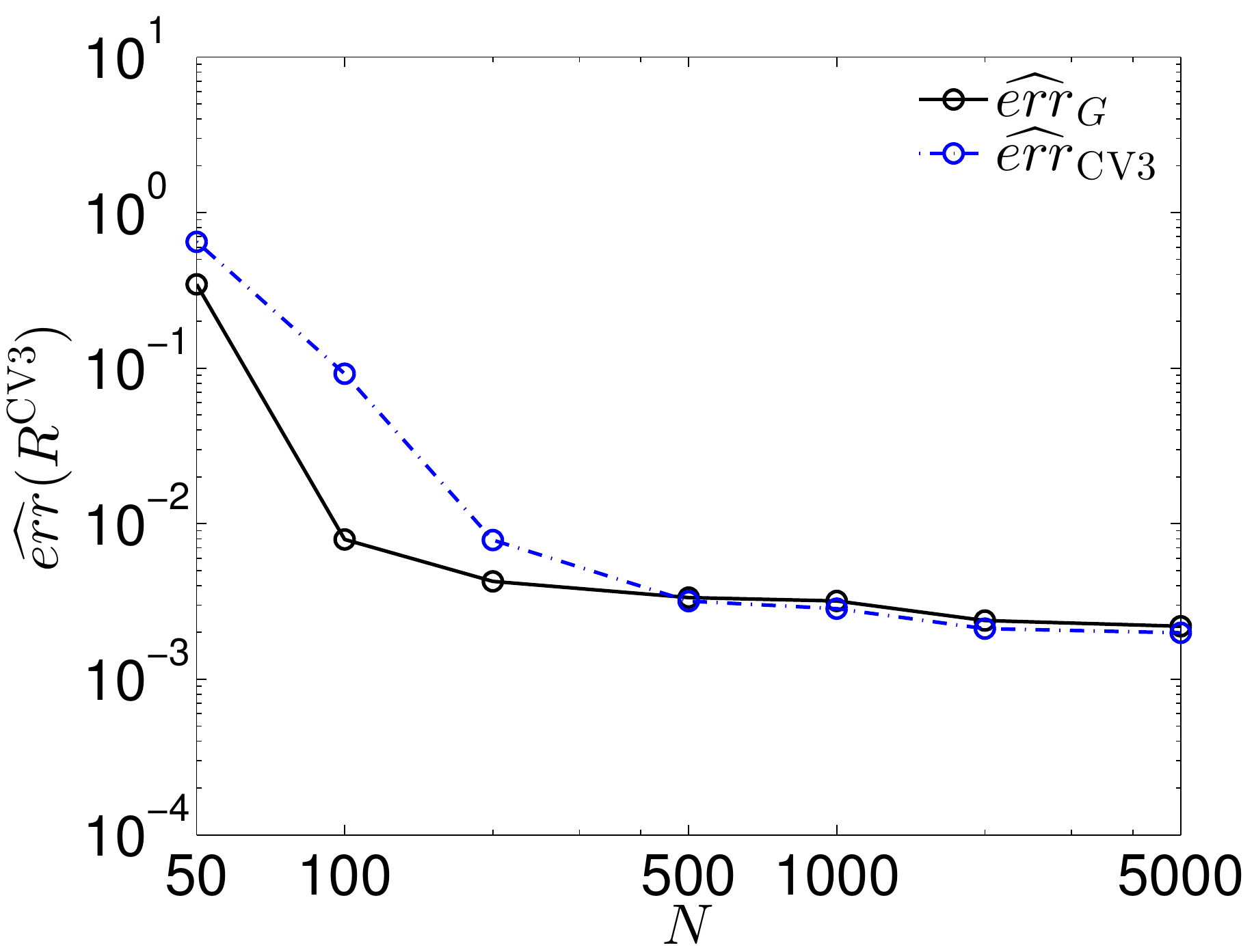}
	\includegraphics[width=0.48\textwidth]{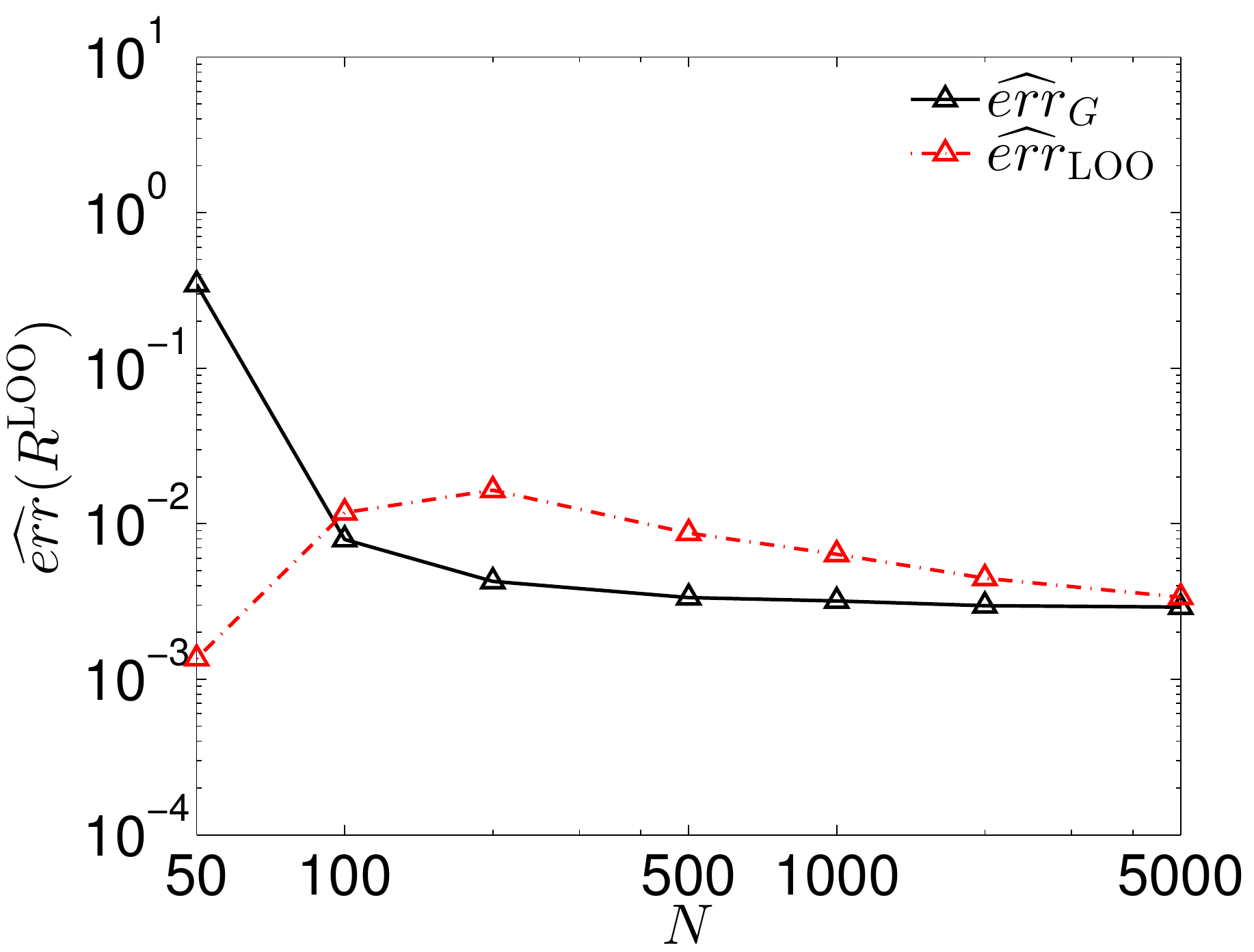}
	\caption{Frame-displacement problem: Comparison of ED-based errors to corresponding relative generalization errors for ranks selected with 3-fold CV (left) and with the LOO-error criterion (right).}
	\label{fig:Frame_rank_2}
\end{figure}

Figure~\ref{fig:Frame_rank_2} provides further insight into rank selection by comparing the ED-based errors $\errCV$ and $\errLOO$ with the corresponding relative generalization errors $\errG$. The left graph plots $\errCV$ together with $\errG$ for LRA with rank $\RCV$, while the right graph plots $\errLOO$ together with $\errG$ for LRA with rank $\RLOO$. The ED-based error measures appear to approach the corresponding generalization errors with increasing ED size, with $\errCV$ providing accurate estimates with only $N=500$. As also observed in the previous examples, $\errLOO$ may largely overestimate the errors of meta-models with higher ranks. 

\subsubsection{Stopping criterion in the correction step}

We next investigate effects of the stopping criteria on the LRA accuracy considering the cases with $N\in\{50; 200; 1,000; 5,000\}$. The left graph of Figure~\ref{fig:Frame_rank_2} shows the relative generalization errors of LRA with optimal ranks while $\Derrmin$ varies between $10^{-6}$ and $10^{-1}$ and other parameters are fixed to their values above. The right graph shows the corresponding maximum number of iterations in a correction step. For $N=200$, the considered variation of $\Derrmin$ corresponds to a variation of $\errG$ of about one order of magnitude. For the other EDs, the value of $\Derrmin$ has a negligible effect on the LRA accuracy and thus, repeated iterations in the correction step add an unnecessary computational effort.

\begin{figure}[!ht]
	\centering
	\includegraphics[width=0.48\textwidth]{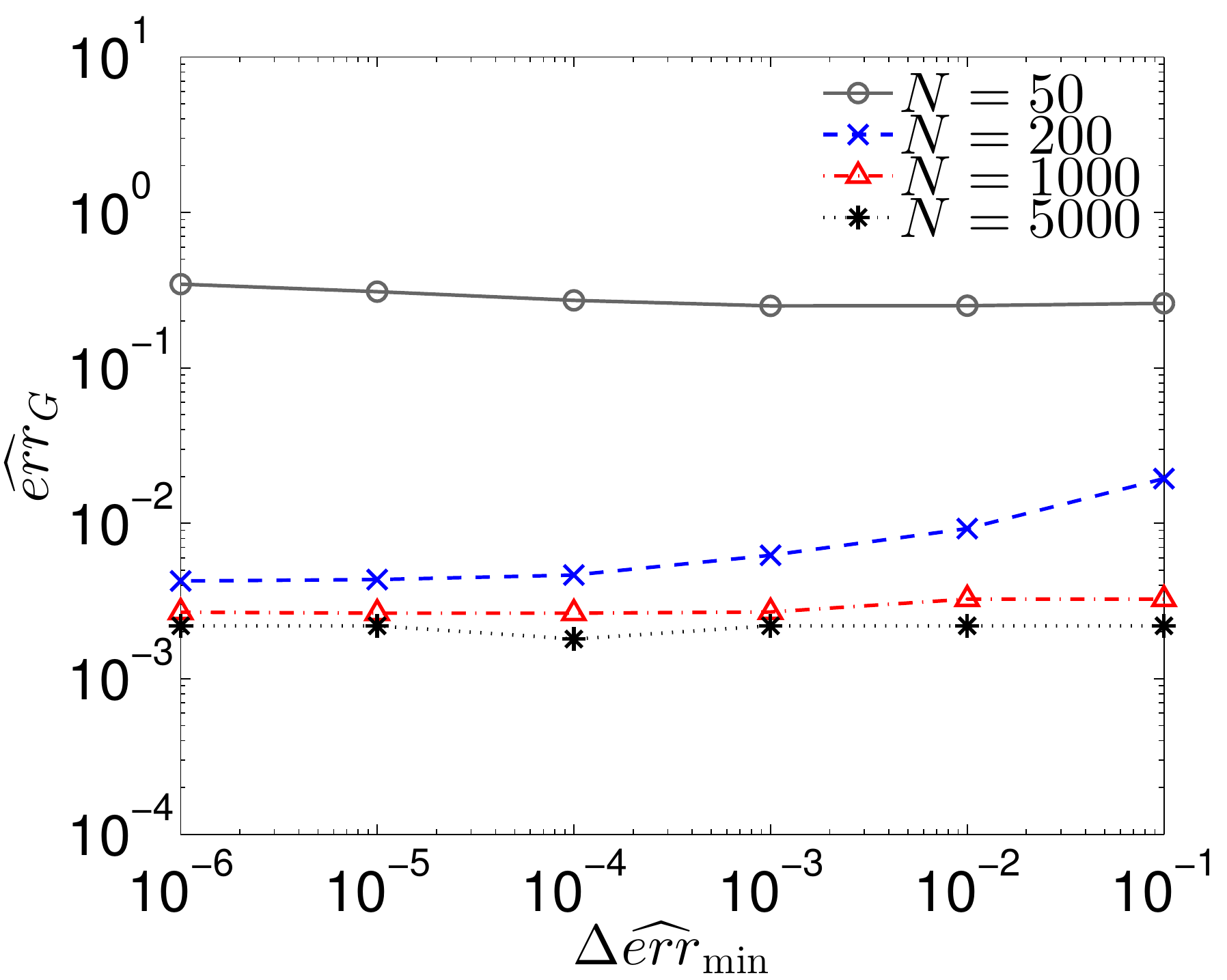}
	\includegraphics[width=0.46\textwidth]{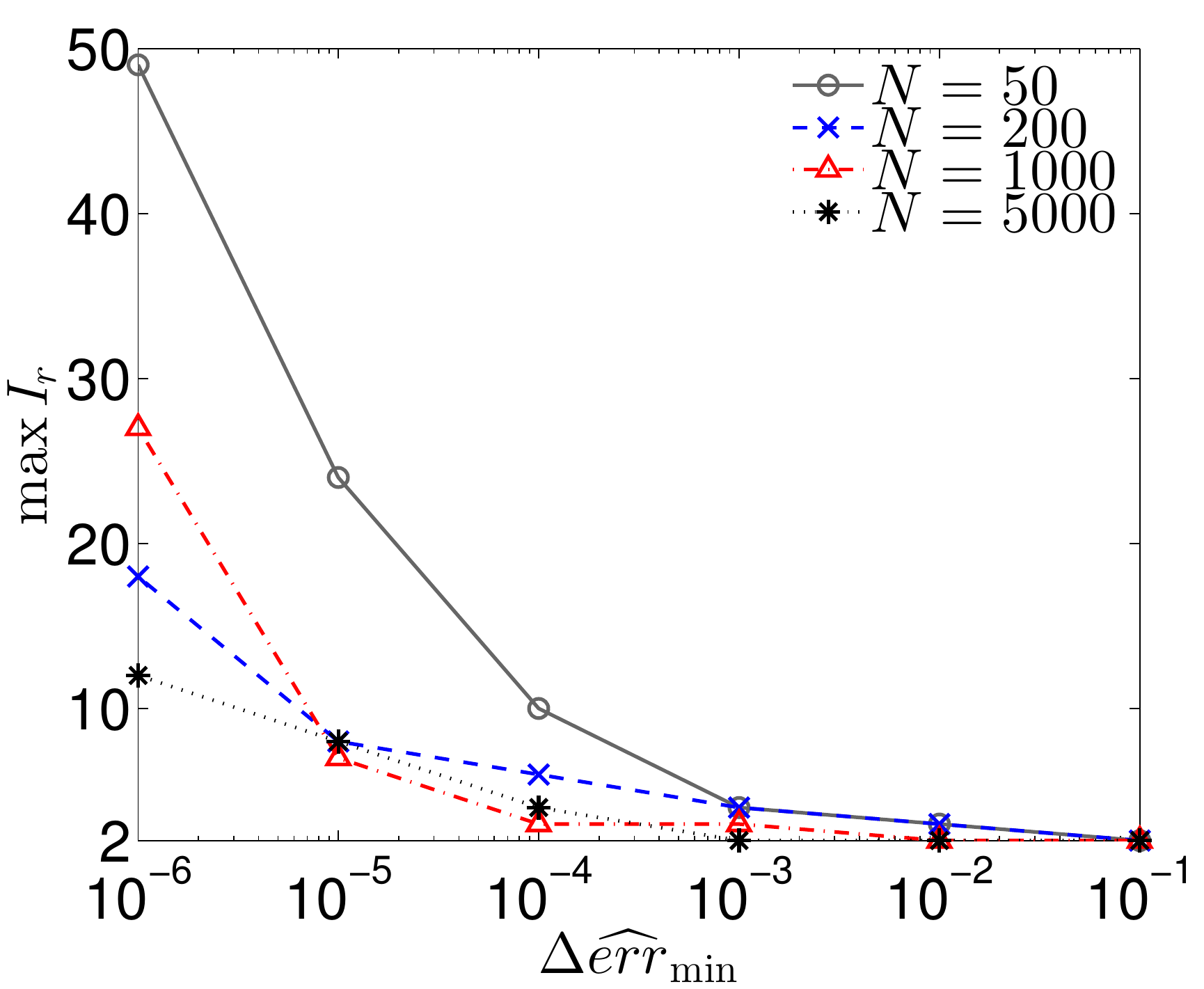}
	\caption{Frame-displacement problem: Relative generalization error (left) and corresponding maximum number of iterations in a correction step (right) versus threshold of differential error in the stopping criterion.}
	\label{fig:Frame_Derr}
\end{figure}

\subsection{Heat conduction}
\label{sec:DiffusionRF_LRA}

\begin{figure}[!ht]
	\centering
	\includegraphics[trim = 0mm 15mm 0mm 10mm, width=0.49\textwidth] {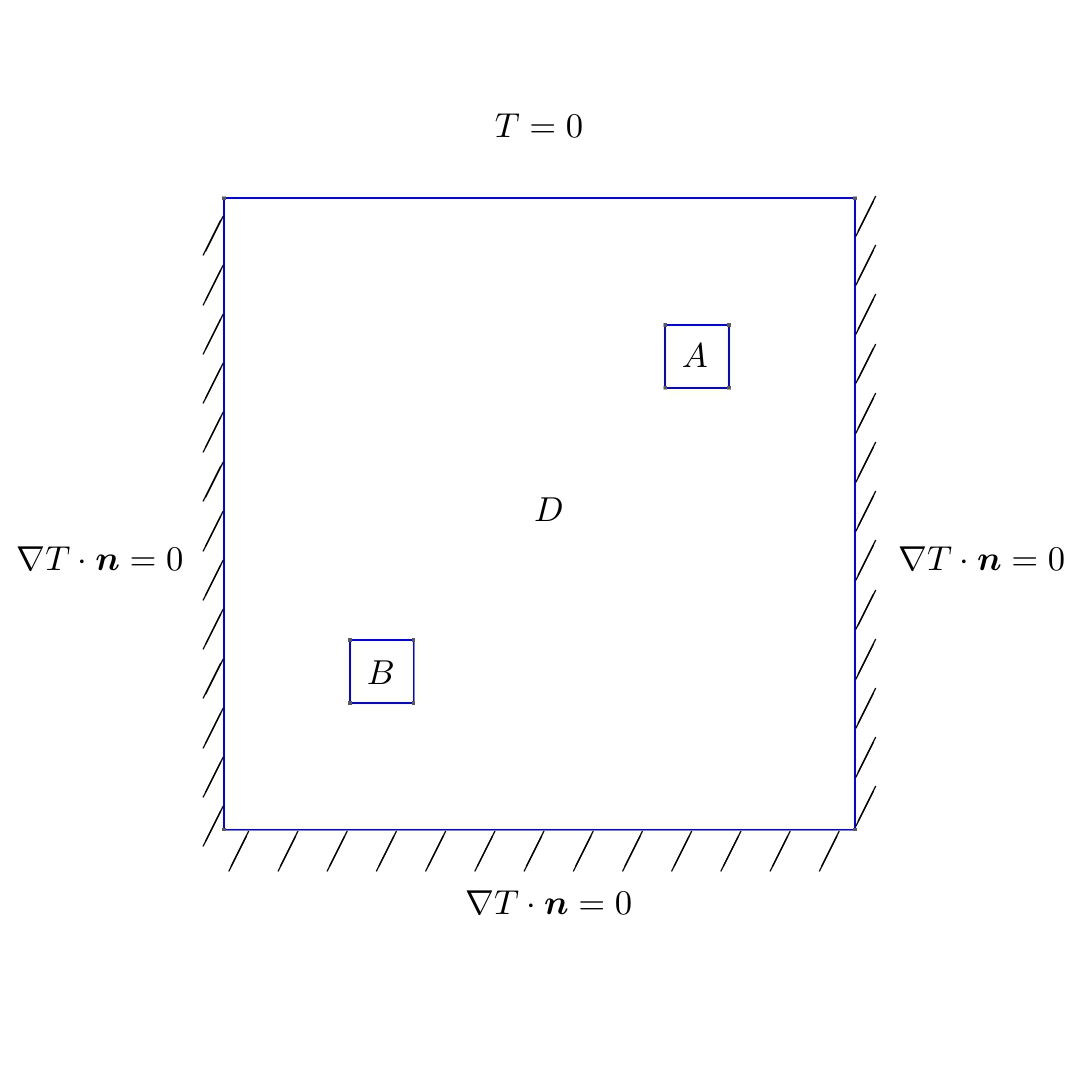}
	\includegraphics[trim = 0mm 15mm 0mm 10mm,width=0.49\textwidth] {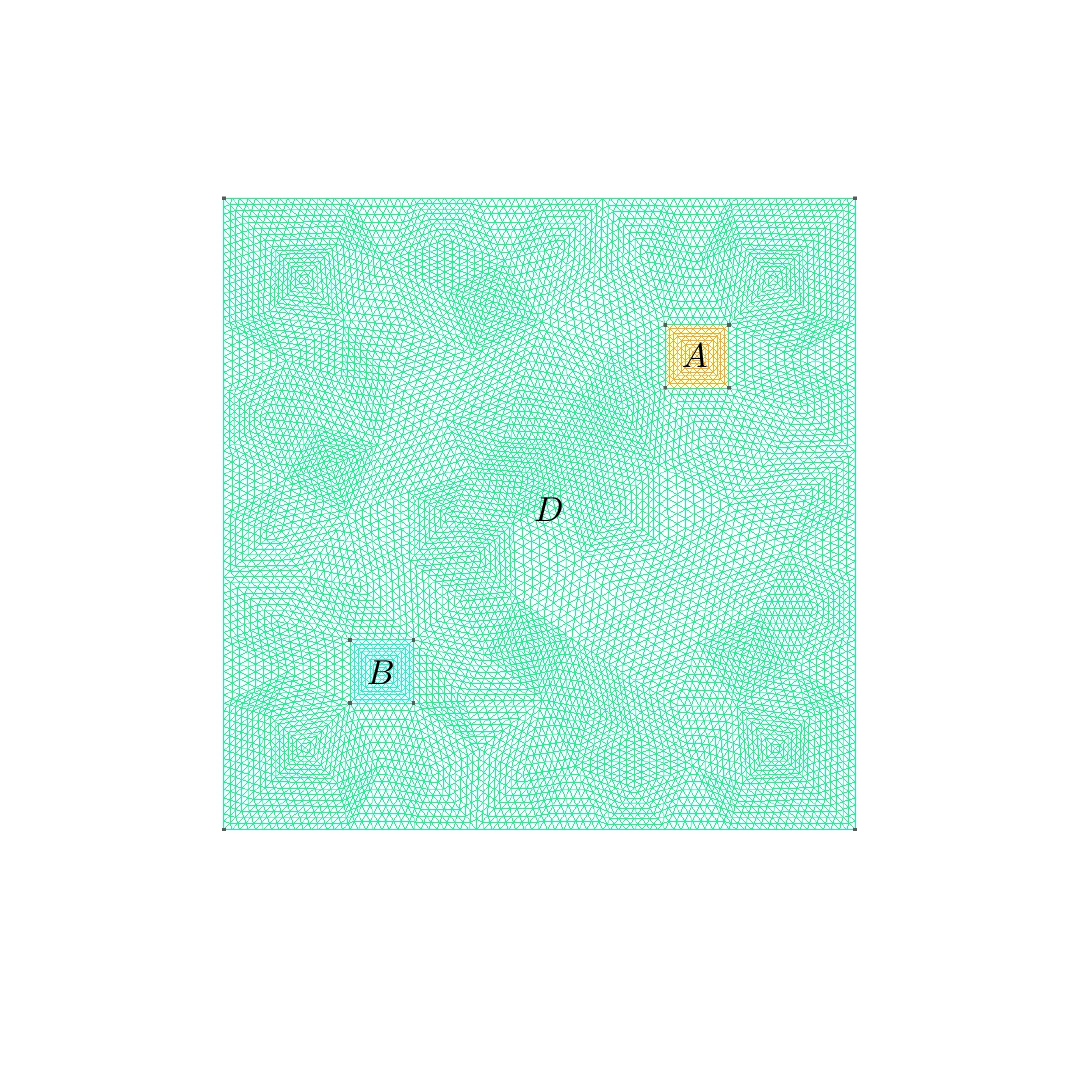}
	\caption{Heat-conduction problem: Domain and boundary conditions (left); finite-element mesh (right).}
	\label{fig:RF_domain}
\end{figure}

The last example application (inspired by a problem studied in \citet{Nouy2010b}) concerns stationary heat conduction in the two-dimensional square domain $D=(-0.5,0.5)\hsp\rm m\times(-0.5,0.5)\hsp\rm m$ shown in the left graph of Figure~\ref{fig:RF_domain}. The temperature field $T(\ve z)$, $\ve z \in D$, is described by the partial differential equation:
\begin{equation}
\label{eq:diffusion_eq}
-\nabla(\kappa(\ve z)\hsp \nabla T(\ve z)) =I_A(\ve z)\hsp Q ,
\end{equation}
with boundary conditions $T=0$ on the top boundary and $\nabla T \cdot \ve n=0$ on the left, right and bottom boundaries, where $\ve n$ denotes the vector normal to the boundary. On the right side of Eq.~(\ref{eq:diffusion_eq}), $Q=2\cdot 10^3\hsp\rm W/m^3$ and $I_A$ is an indicator function equal to unity if $\ve z\in A$, where $A=(0.2,0.3)\hsp\rm m\times(0.2,0.3)\hsp\rm m$ is a square domain within $D$ (see left graph of Figure~\ref{fig:RF_domain}). The thermal conductivity $\kappa(\ve z)$ is a lognormal random field described by:
\begin{equation}
\label{eq:diffusion_coef}
\kappa(\ve z)=\exp[a_{\kappa}+b_{\kappa} \hsp g(\ve z)],
\end{equation}
in which $g(\ve z)$ represents a standard Gaussian random field with a square-exponential autocorrelation function:
\begin{equation}
\label{eq:corr_fcn}
\rho(\ve z,\ve z')=\exp{(-\|\ve z-\ve z'\|^2/\ell^2)}.
\end{equation}
In Eq.~(\ref{eq:diffusion_coef}), the values of $a_{\kappa}$ and $b_{\kappa}$ are such that the mean and standard deviation of $\kappa$ are $\mu_{\kappa}=1\hsp\rm W/°C\cdot m$ and $\sigma_{\kappa}=0.3\hsp\rm W/°C\cdot m$ respectively, while in Eq.~(\ref{eq:corr_fcn}), $\ell=0.2\hsp\rm m$. The response quantity of interest is the average temperature in the square domain $B=(-0.3,-0.2)\hsp\rm m\times(-0.3,-0.2)\hsp\rm m$ within $D$ (see left graph of Figure~\ref{fig:RF_domain}):
\begin{equation}
\label{eq:diffusion_Y}
\widetilde{T}=\frac{1}{|B|}\int_{\ve z \in B}T(\ve z)\hsp d\ve z.
\end{equation}

\begin{figure}[!ht]
	\centering
	\includegraphics[trim = 25mm 25mm 25mm 25mm, width=0.60\textwidth] {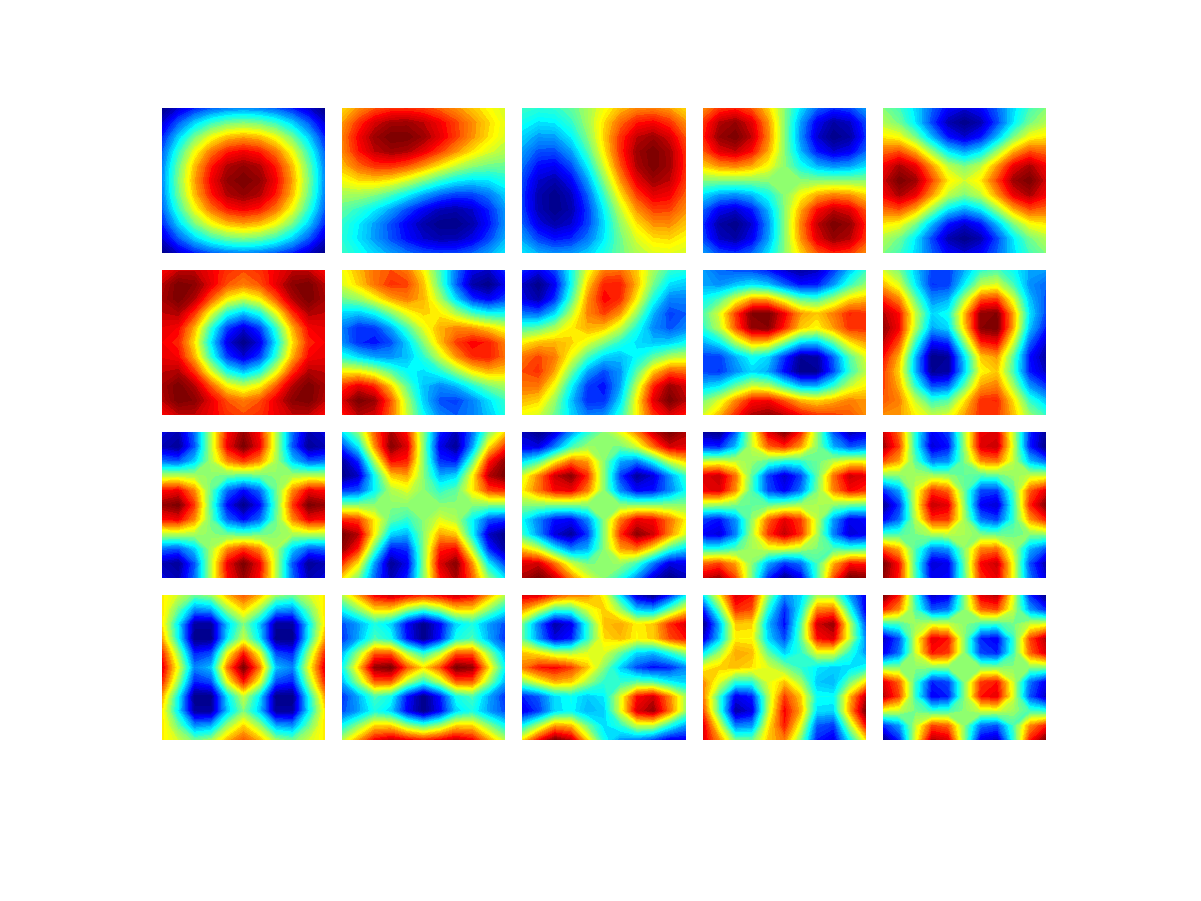}
	\caption{Heat-conduction problem: Shapes of the first 20 basis functions in the EOLE discretization (from left-top to bottom-right row-wise).}
	\label{fig:RF_modes}
\end{figure}

To solve Eq.~(\ref{eq:diffusion_eq}), the Gaussian random field $g(\ve z)$ is first discretized using the expansion optimal linear estimation (EOLE) method \citep{Li1993optimal}, as described next. Let $\{\ve{\zeta}_1 \enum \ve{\zeta}_n\}$ denote the points of an appropriately defined grid in $D$. By retaining the first $M$ terms in the EOLE series, $g(\ve z)$ is approximated by:
\begin{equation}
\label{eq:EOLE}
\widehat{g}(\ve z) = \sum_{i=1}^M \frac{\xi_i}{\sqrt{l_i}}\ve {\phi}_i^{\rm{T}} \ve C_{\ve z \ve \zeta}(\ve z),
\end{equation}
where $\{\xi_1 \enum \xi_M\}$ are independent standard normal variables, $\ve C_{\ve z \ve \zeta}$ is a vector with elements $\ve C_{\ve z \ve \zeta}^{(k)}=\rho(\ve z,\ve \zeta_k)$, for $k=1 \enum n$, and $(l_i,\ve{\phi}_i)$ are the eigenvalues and eigenvectors of the correlation matrix $\ve C_{\ve \zeta \ve \zeta}$ with elements $\ve C_{\ve \zeta \ve \zeta}^{(k,l)}=\rho(\ve{\zeta}_k,\ve{\zeta}_l)$, for $k,l=1 \enum n$. \citet{Sudret2000stochastic} recommend that for a square-exponential autocorrelation function, the size of the element in the EOLE grid must be $1/2-1/3$ of $\ell$. In the present numerical application, we use a square grid with element size $0.1\hsp \rm m$ ($1/2$ of $\ell$), thus comprising $n=121$ points. The number of terms in the EOLE series is determined according to the rule:
\begin{equation}
\label{eq:EOLE_M}
\sum_{i=1}^{M}l_i/\sum_{i=1}^{n} l_i \geq 0.99,
\end{equation}
herein leading to $M=53$. The shapes of the first 20 basis functions $\{{\phi}_i^{\rm{T}} \ve C_{\ve z \ve \zeta}(\ve z), \hsp i=1 \enum 20\}$ are shown in Figure~\ref{fig:RF_modes}.

\begin{figure}[h]
	\centering
	\includegraphics[width=0.46\textwidth] {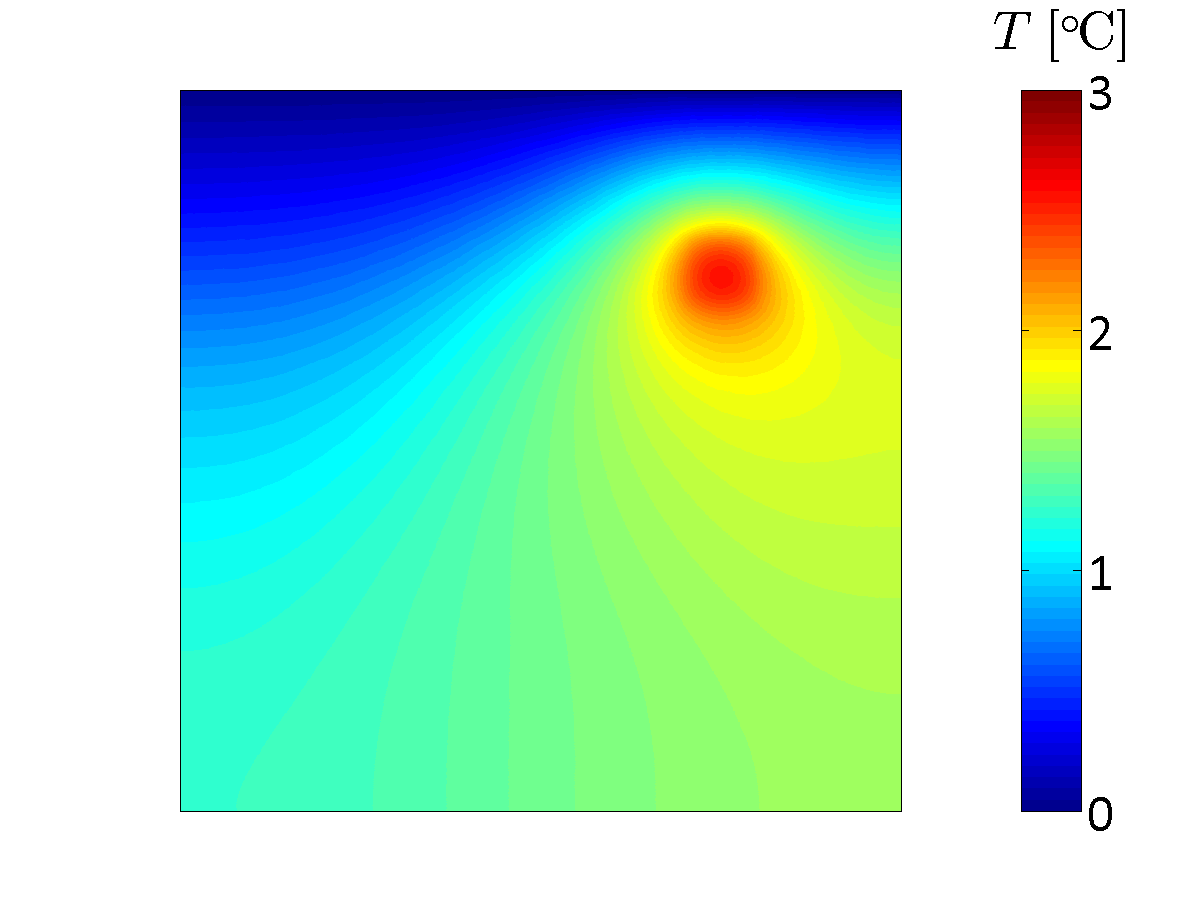}
	\includegraphics[width=0.46\textwidth] {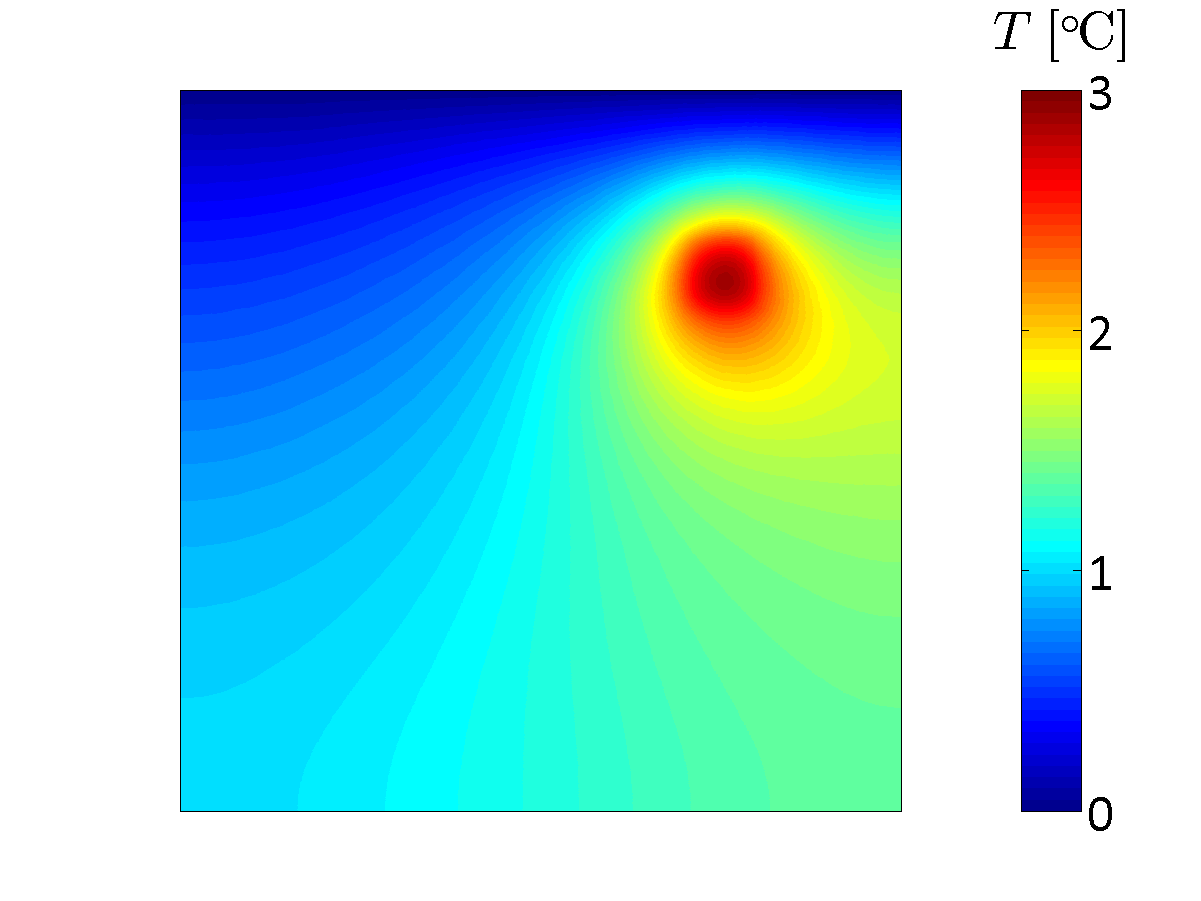}
	\caption{Heat-conduction problem: Example realizations of the temperature field.}
	\label{fig:maps}
\end{figure}

The underlying deterministic problem is solved with an in-house finite-element analysis code developed in the Matlab environment. The domain $D$ is discretized into $16,000$ triangular T3 elements, as shown in the right graph of Figure~\ref{fig:RF_domain}, using the software \emph{Gmsh} \citep{Geuzaine2009gmsh}. The temperature field $T(\ve z)$ for two realizations of the conductivity random field is depicted in  Figure~\ref{fig:maps}. Because the input vector $\ve X=\{\xi_1 \enum \xi_M\}$ consists of independent standard normal variables, we develop LRA representations of the random response $\widetilde{T}=\cm(\ve X)$ with basis functions made of Hermite polynomials.

\subsubsection{Rank selection and error measures}

In examining the accuracy of the two rank-selection criteria and associated error measures, we consider EDs of size varying between $N=50$ and $N=5,000$ and a validation set of size $\nval=10^4$. After preliminary investigations, the maximum polynomial degree is set to $p_i=p=2$ for $\{i=1 \enum 53\}$ and the stopping criterion in the correction step is defined by $I_{\rm max}=50$ and $\Derrmin=10^{-5}$. As in the previous examples, candidate ranks are selected among $\{1 \enum 20\}$.

\begin{figure}[!ht]
	\centering
	\includegraphics[width=0.46\textwidth]{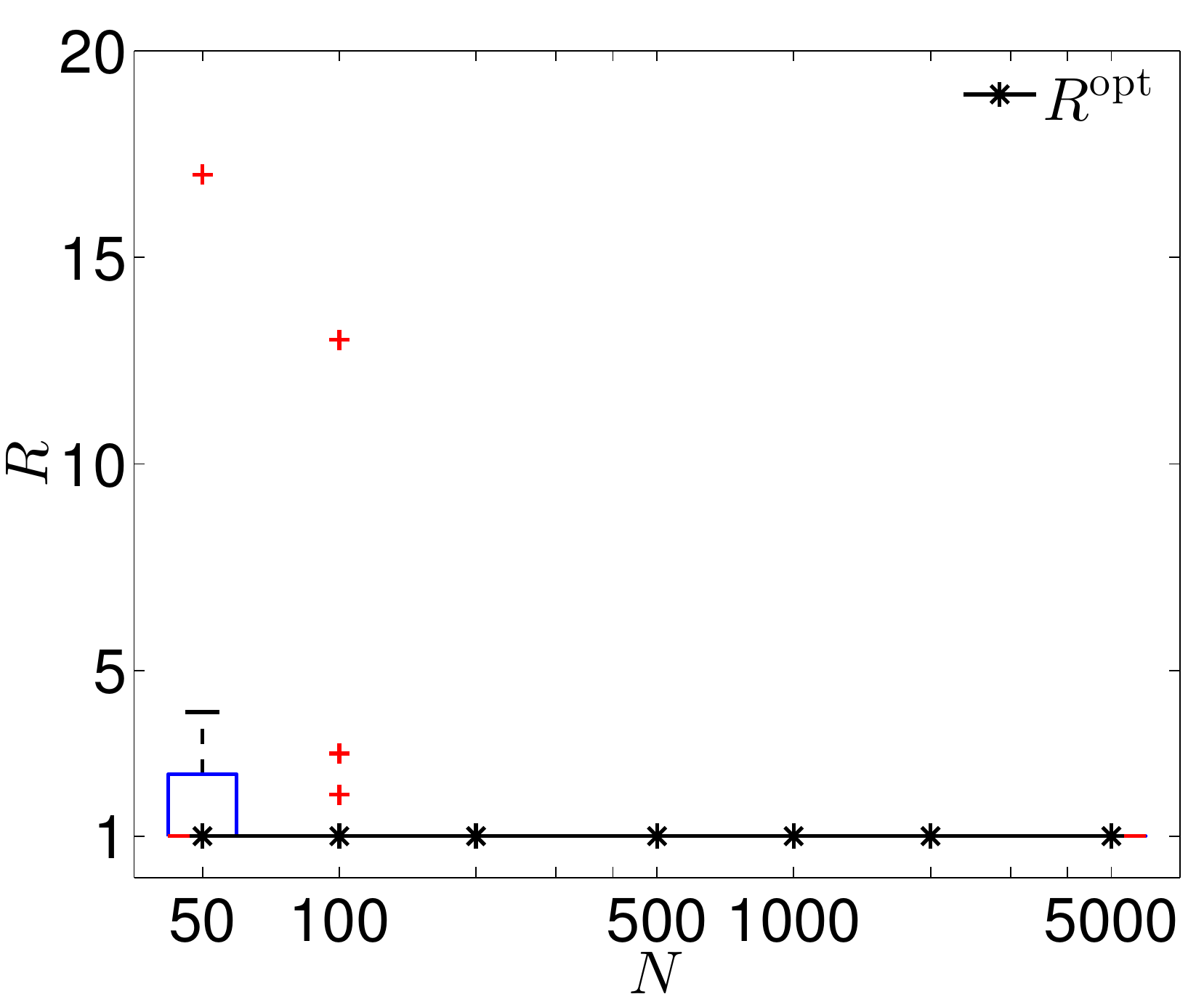}
	\includegraphics[width=0.48\textwidth]{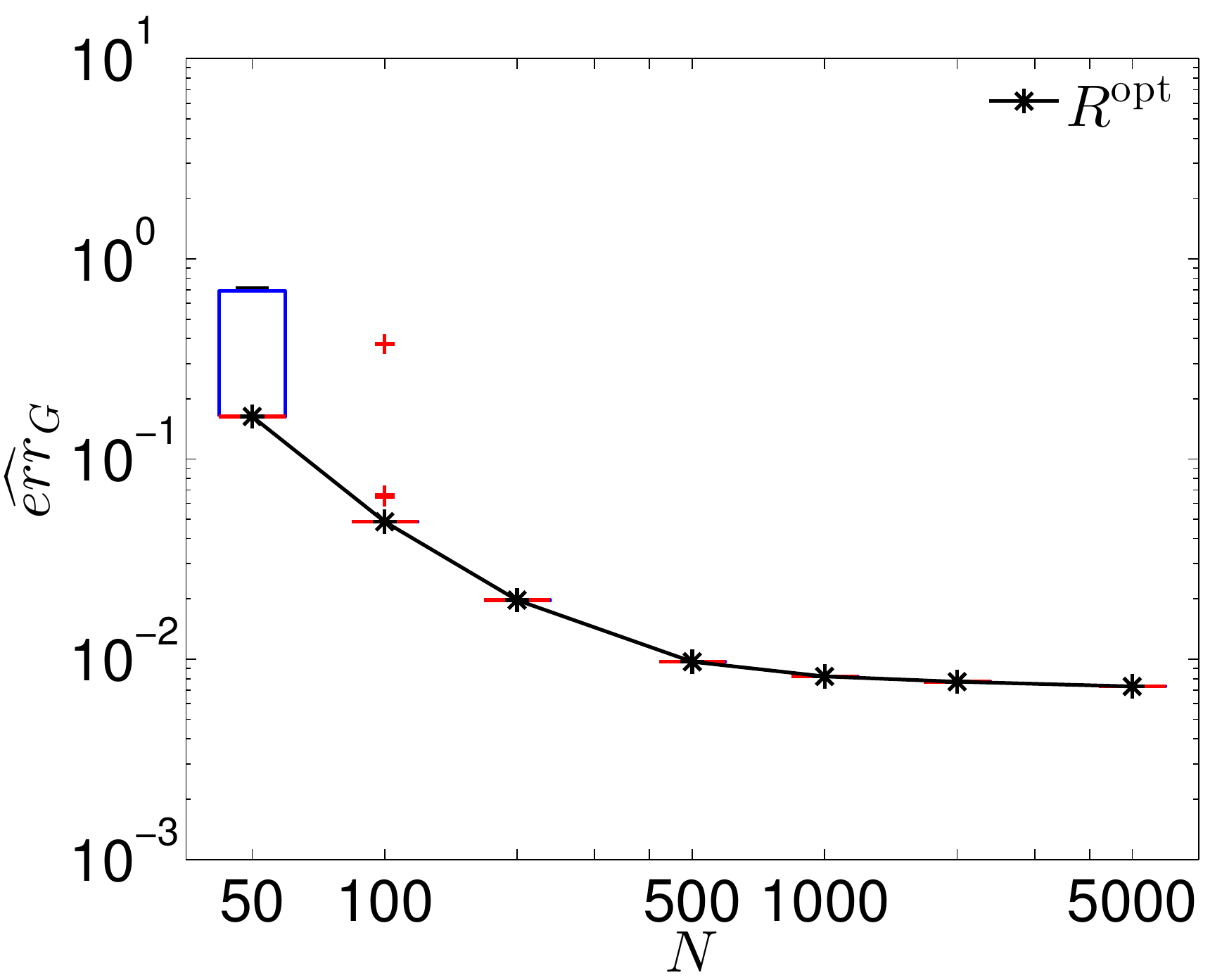}
	\caption{Heat-conduction problem: Comparison of ranks selected with 3-fold CV (20 replications) to optimal ranks (left) and corresponding relative generalization errors (right).}
	\label{fig:DiffusionRF_rank}
\end{figure}

The left graph of Figure~\ref{fig:DiffusionRF_rank} shows boxplots of the ranks selected with 3-fold CV for 20 random partitions of each ED. The optimal rank, also indicated in the graph, is equal to unity in all cases. For $N\geq 200$, the 3-fold CV approach consistently selects the rank-1 meta-models in all trials. For the two smaller EDs however, it may erroneously select higher ranks, with the effect on the LRA accuracy depicted in the right graph of the same figure. The LOO error criterion identifies the optimal rank $\RLOO=\Ropt=1$ for all considered EDs.

\begin{figure}[!ht]
	\centering
	\includegraphics[width=0.48\textwidth]{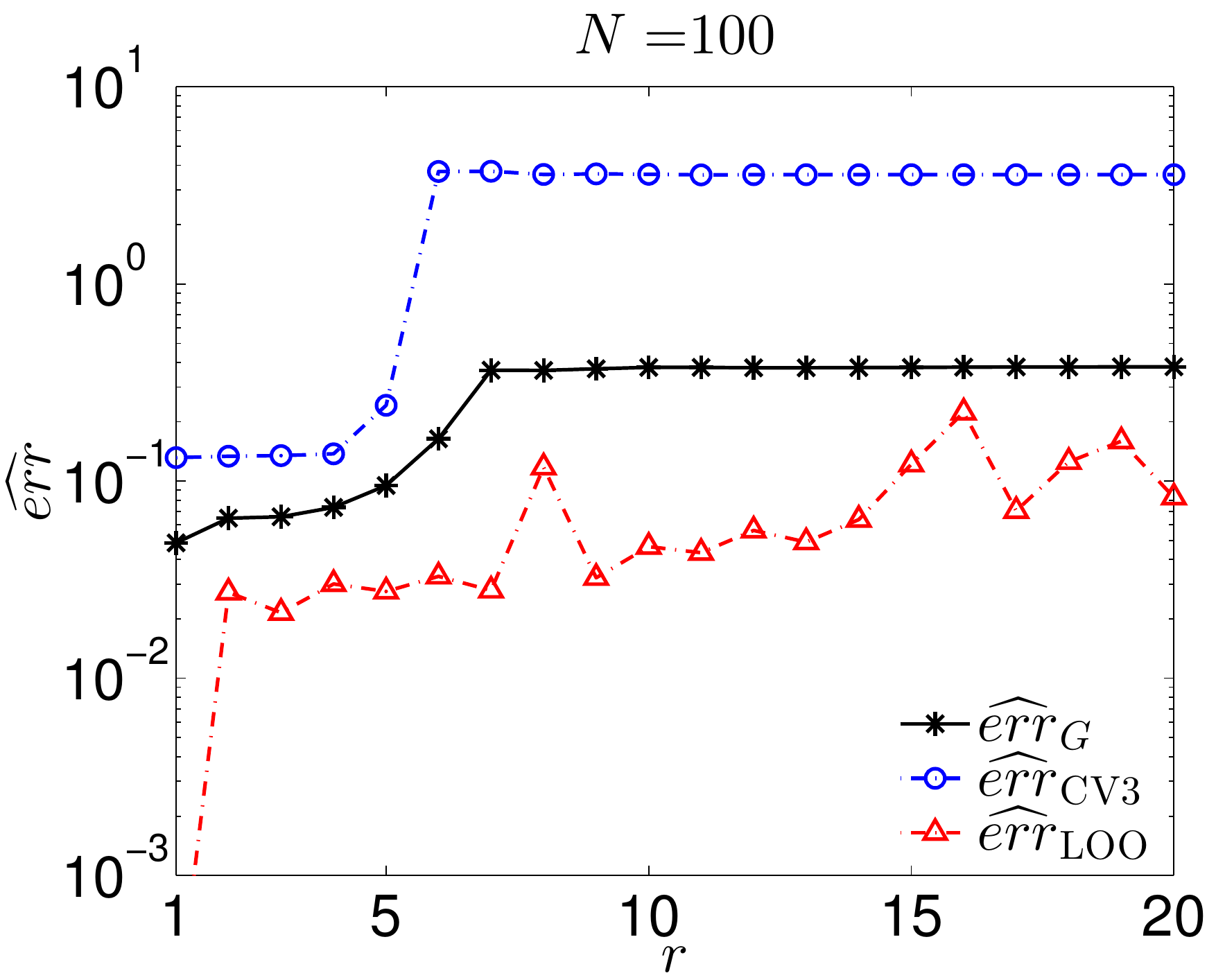}
	\includegraphics[width=0.48\textwidth]{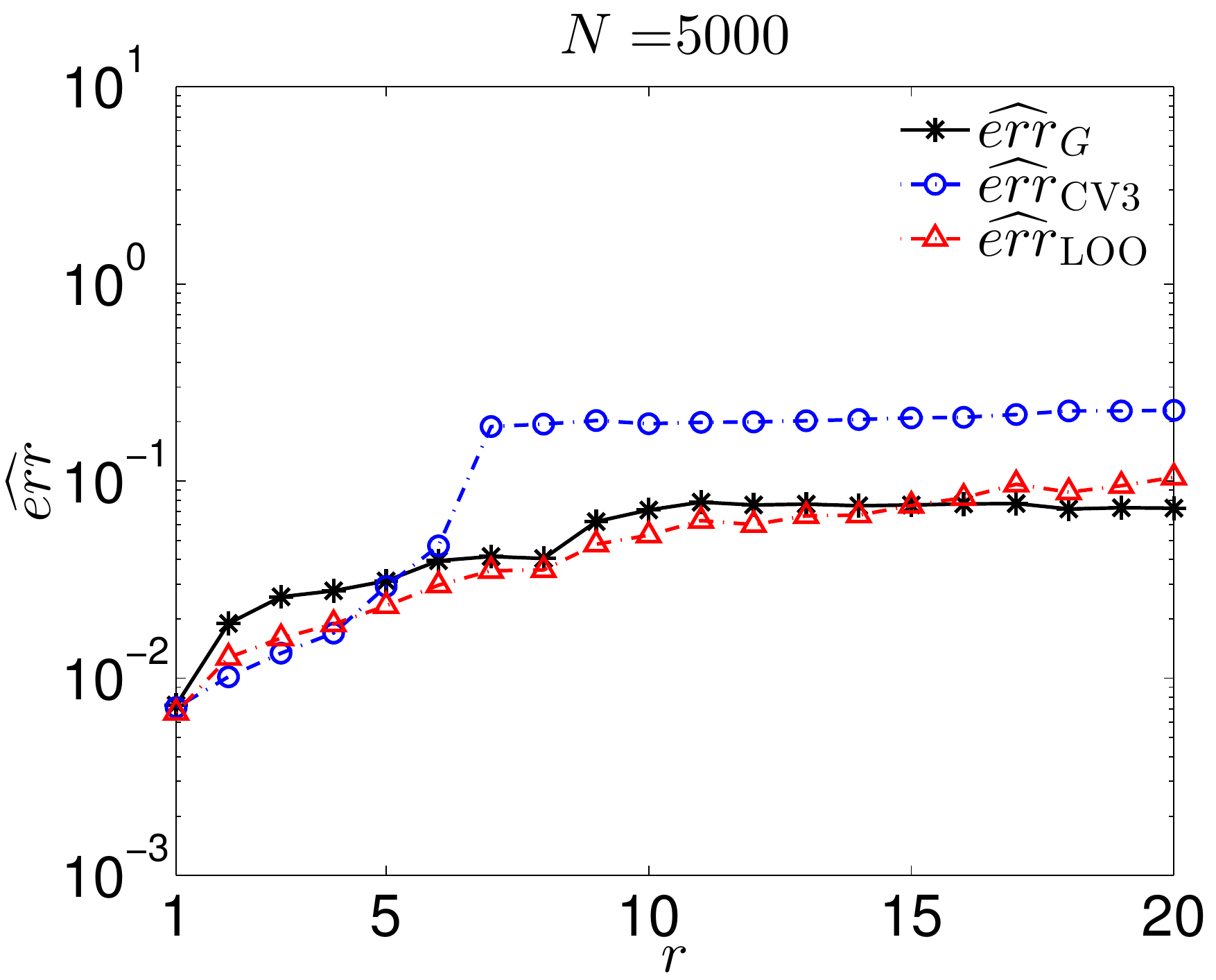}
	\caption{Heat-conduction problem: Comparison of error measures in rank-selection criteria with the relative generalization error based on the validation set.}
	\label{fig:DiffusionRF_errors}
\end{figure}

In Figure~\ref{fig:DiffusionRF_errors}, we compare $\errG$ with $\errCV$ and $\errLOO$ for all candidate ranks and two example EDs of size $N=100$ (left) and $N=5,000$ (right). The $\errCV$ errors for one example partition of the ED are shown.  Note in these graphs that both error criteria identify the unity rank as optimal. In this example, $\errLOO$ can provide reasonably good estimates of $\errG$ even for higher ranks; however $\errCV$ overall outperforms $\errLOO$ in the estimation of $\errG$ at $\Ropt=1$, which is herein evident in the case of $N=100$.

\subsubsection{Stopping criterion in the correction step}

\begin{figure}[!ht]
	\centering
	\includegraphics[width=0.48\textwidth]{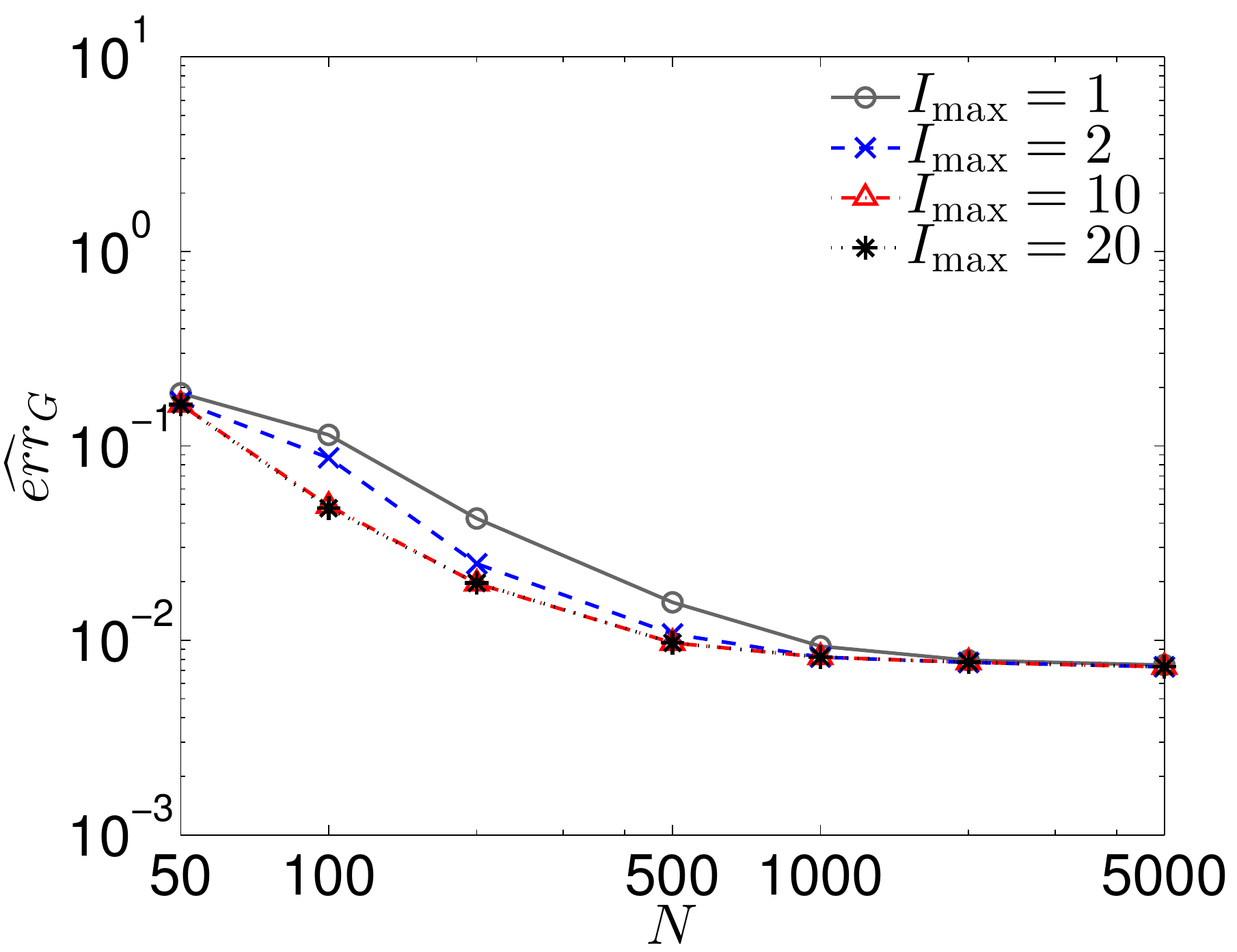}
	\includegraphics[width=0.48\textwidth]{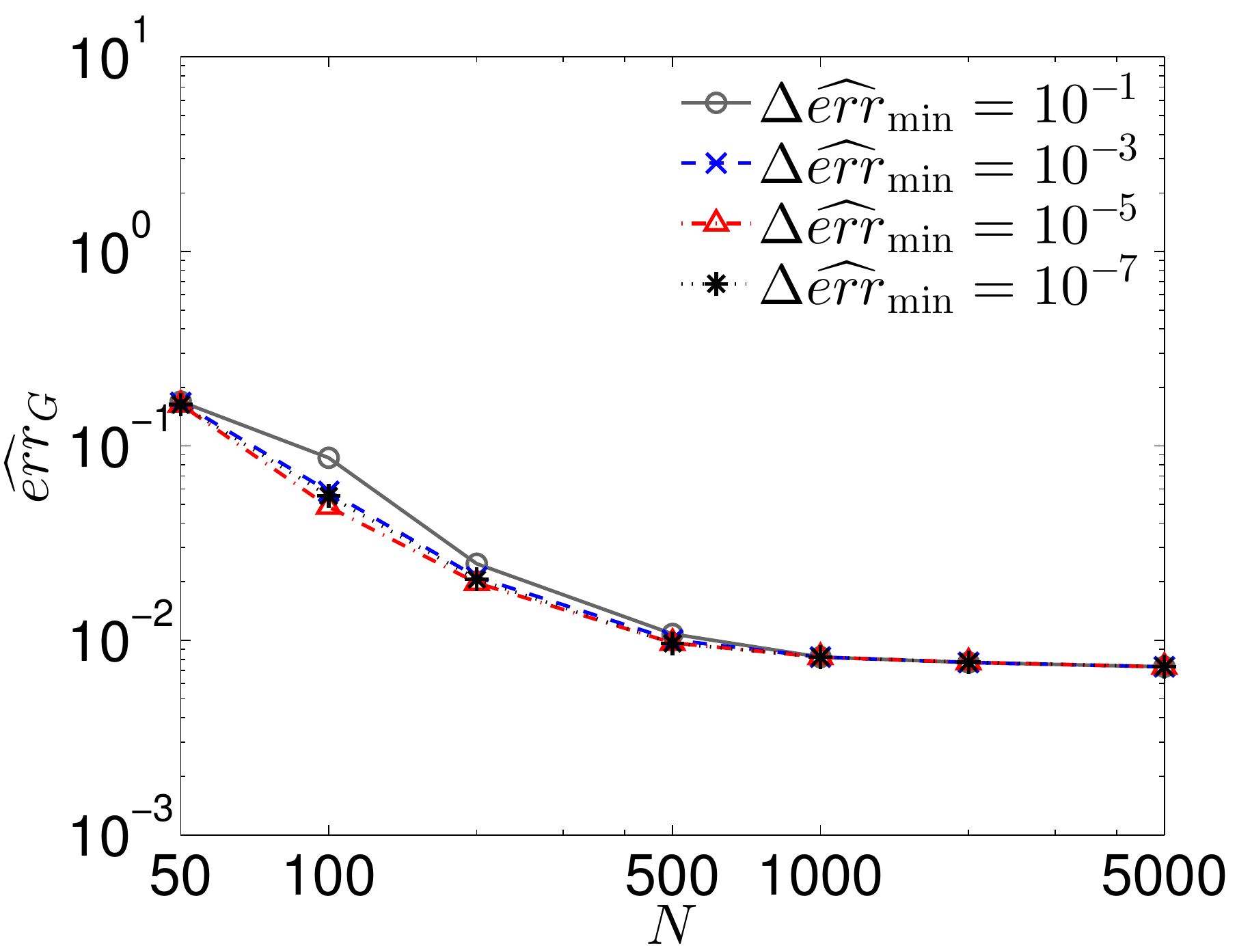}
	\caption{Heat-conduction problem: Relative generalization error for varying stopping criteria in the correction step versus size of experimental design.}
	\label{fig:Diffusion_Derr}
\end{figure}

Finally, we investigate effects of the parameters of the stopping criterion on the LRA accuracy. The left graph of Figure~\ref{fig:Diffusion_Derr} shows $\errG$ for LRA with rank $\Ropt=1$ versus $N$, while the maximum allowable number of iterations takes the values $\Imax=1,2,10,20$ and the differential error threshold is fixed to $\Derrmin=10^{-5}$. The right graph of the same figure shows $\errG$ for LRA with rank $\Ropt=1$ versus $N$, while the differential error threshold takes the values $\Derrmin=10^{-1},10^{-3},10^{-5},10^{-7}$ and the maximum allowable number of iterations is fixed to $I_{\rm max}=50$. Overall, both $I_{\rm max}$ and $\Derrmin$ have a relatively small effect on the LRA accuracy, which becomes negligible for $N>500$.

\subsection{Discussion on the construction of canonical low-rank approximations}
\label{sec:discussion}

In the above numerical investigations, the 3-fold CV error proves a more reliable criterion for rank selection as compared to the LOO error. The lack of orthogonality of the regressors may lead to excessive correction factors in the latter, particularly for higher ranks, whereas its uncorrected counterpart tends to be overly optimistic. On the other hand, the accuracy of the 3-fold CV error consistently improves with increasing size of the ED. For a certain ED, the median generalization error of LRA with rank selected via 3-fold CV is fairly close to the generalization error of the LRA with optimal rank. Divergence of the selected rank from the optimal one tends to have a smaller effect on the LRA accuracy as the ED size increases. In the examples involving finite-element models, EDs of size $N\geq10 \cdot M$ lead to LRA with generalization errors of the order of $10^{-3}$ or smaller, even in cases with non-optimal rank selection. Note that this level of accuracy is typically sufficient in several meta-modeling applications, including sensitivity analysis \citep{Konakli2016RESS}. For the analytical beam-deflection function with underlying rank-1 structure, the same condition on the ED size leads to even higher meta-model accuracy.

Optimal values for the parameters in the stopping criterion of the correction step appear to strongly depend on the specific application. Overall, setting a sufficiently low differential error threshold $\Derrmin$ tends to be more critical when small EDs are considered ($N\leq10 \cdot M$ in the examined applications). Based on the above results, we recommend the use of $\Derrmin$ values in the range $10^{-5}-10^{-6}$ and caution that lower values may lead to numerical instabilities in certain cases. By imposing a simultaneous constraint on the maximum number of iterations $\Imax$, excessive unnecessary iterations are avoided. In the above case studies, setting $\Imax=50$ allows the required number of iterations for achieving nearly the maximum or a sufficient meta-model accuracy. Note that depending on the application, the required number of iterations might be considerably smaller than that imposed by the recommended thresholds. Nevertheless, because these iterations involve minimization problems of only small size, they are relatively inexpensive from a computational viewpoint. We underline that in a typical realistic meta-modeling application, the main computational effort lies in the evaluation of the model responses at the ED.


\section{Canonical low-rank approximations versus sparse polynomial chaos expansions}
\label{sec:LRAvsPCE}

We next confront canonical LRA to sparse PCE in the same meta-modeling applications considered in Section~\ref{sec:Examples_LRA}. The focus of the comparison is set on the applications involving finite-element models. The beam-deflection problem is only briefly examined in order to confirm the superior performance of LRA when the original model has a low-rank structure. The two types of meta-models are built using polynomials from the same family. Consistently with Section~\ref{sec:Examples_LRA}, Hermite polynomials are herein used in all four applications. Following the discussion above, we build LRA with rank $\RCV$ (one random partition per ED is considered). A common polynomial degree is set in all dimensions, with its value $p \in \{1 \enum 20\}$ also selected via 3-fold CV (one random partition per ED is considered). The associated ED-based error estimate of the meta-model is denoted by $\CVLRA$. In building sparse PCE, we determine the candidate basis with a hyperbolic truncation scheme (see Section~\ref{sec:PCE_basis}) and compute the coefficients with hybrid LAR (see Section~\ref{sec:PCE_coef}). An optimal combination of the parameters $p^t \in \{1 \enum 20\}$ and $q \in \{0.25, 0.50, 0.75, 1\}$ is selected in terms of the minimum corrected LOO error, hereafter denoted by $\LOOPCE$. The PCE computations are performed with the software UQLab \citep{MarelliICVRAM2014, UQdoc_09_104}. Comparisons between LRA and PCE are based on EDs obtained with Sobol sequences as well as Latin Hypercube Sampling (LHS). Of interest are cases with small ED sizes, which are typically encountered in real-life meta-modeling applications. Similarly to Section~\ref{sec:Examples_LRA}, large MCS validation sets are used to determine the actual errors of the meta-models. Because however such validation sets are not available in real-life applications, we further examine the accuracy of the ED-based error measures $\CVLRA$ and $\LOOPCE$ in estimating the corresponding generalization errors $\errGLRA$ and $\errGPCE$.

\subsection{Beam deflection}

\begin{figure}[!ht]
	\centering
	\includegraphics[width=0.47\textwidth]{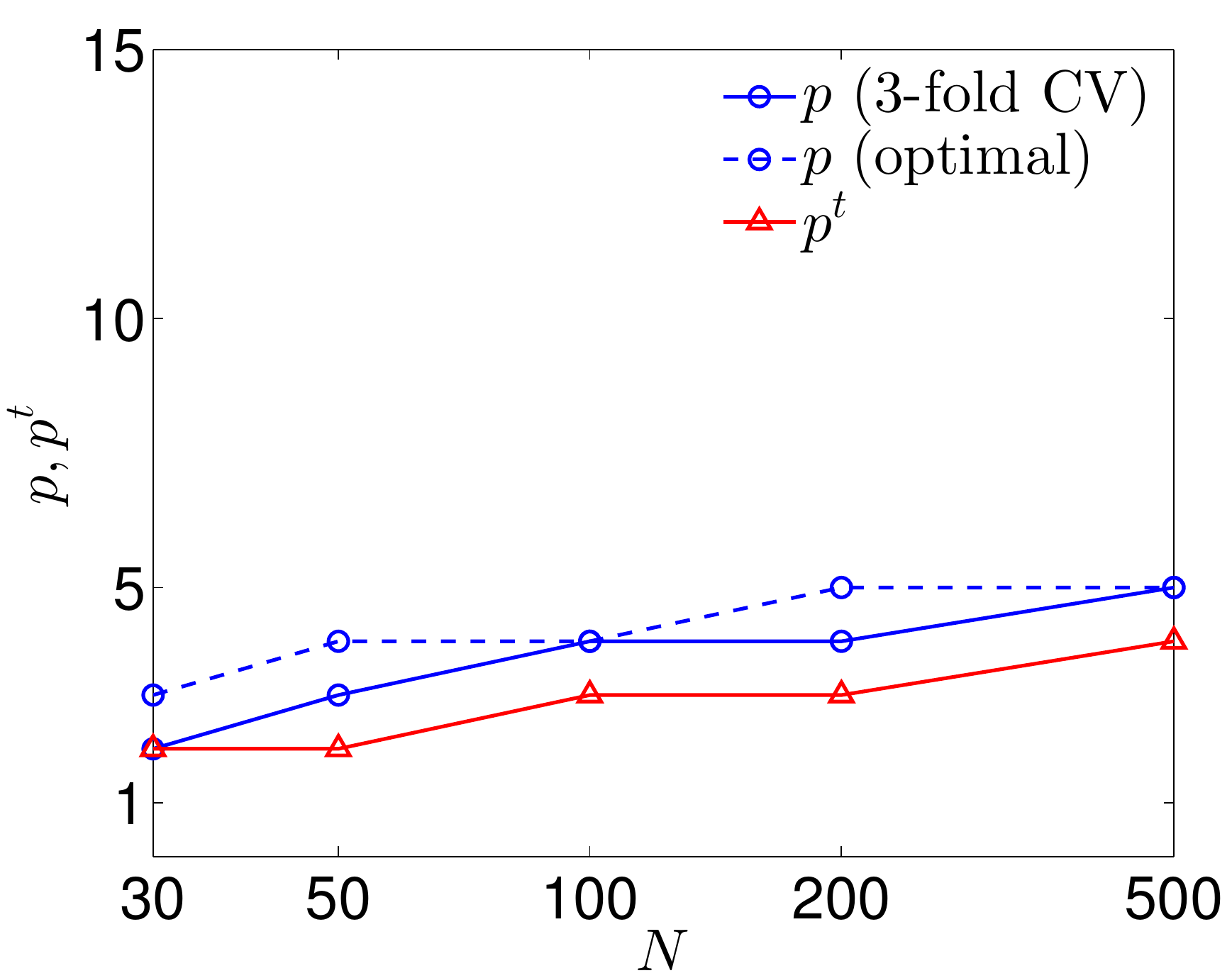}
	\includegraphics[width=0.48\textwidth]{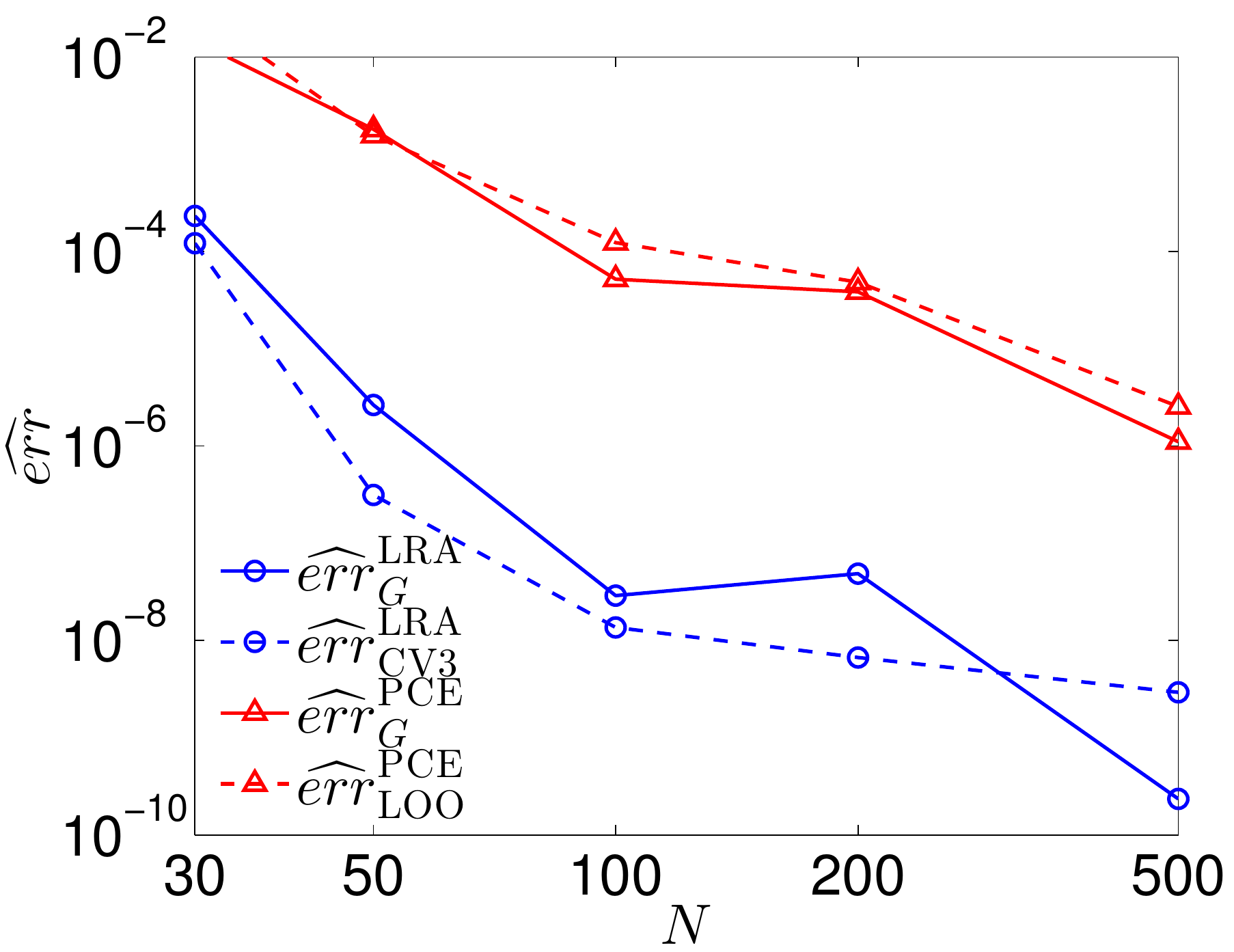}
	\caption{Beam deflection: Polynomial degrees (left) and corresponding error measures (right) of LRA and PCE meta-models based on Sobol sequences.}
	\label{fig:Beam_LRA_PCE}
\end{figure}

We herein develop LRA and PCE meta-models of the beam-deflection function described in Section~\ref{sec:Beam_LRA}, using EDs of size $30\leq N\leq500$ based on Sobol sequences. We assess the accuracy of the meta-models using a validation set of size  $\nval=10^6$. The left graph of Figure~\ref{fig:Beam_LRA_PCE} compares the maximum polynomial degree $p$ of the univariate polynomials in LRA  with the maximum total polynomial degree $p^t$  in PCE. In addition to the degree $p$ selected by means of 3-fold CV, the graph also shows the actual optimal degree based on the validation set; the two values either agree or differ by one. The ED-based and generalization errors of the two types of meta-models are shown in the right graph of the figure. Note that the depicted errors of LRA correspond to the degree selected by means of 3-fold CV, \ie we assess the accuracy of meta-models developed by using the information contained in the ED only. In the present case where the underlying model has a rank-1 structure, the LRA meta-models are 2-3 orders of magnitude more accurate than the PCE ones. It is noteworthy that EDs of size as small as $N=30-50\leq10\cdot M$ yield highly accurate LRA meta-models with errors of order $10^{-4}-10^{-6}$. The generalization errors of both LRA and PCE are approximated fairly well by the corresponding ED-based measures.

\subsection{Truss deflection}

\begin{figure}[!ht]
	\centering
	\includegraphics[width=0.47\textwidth]{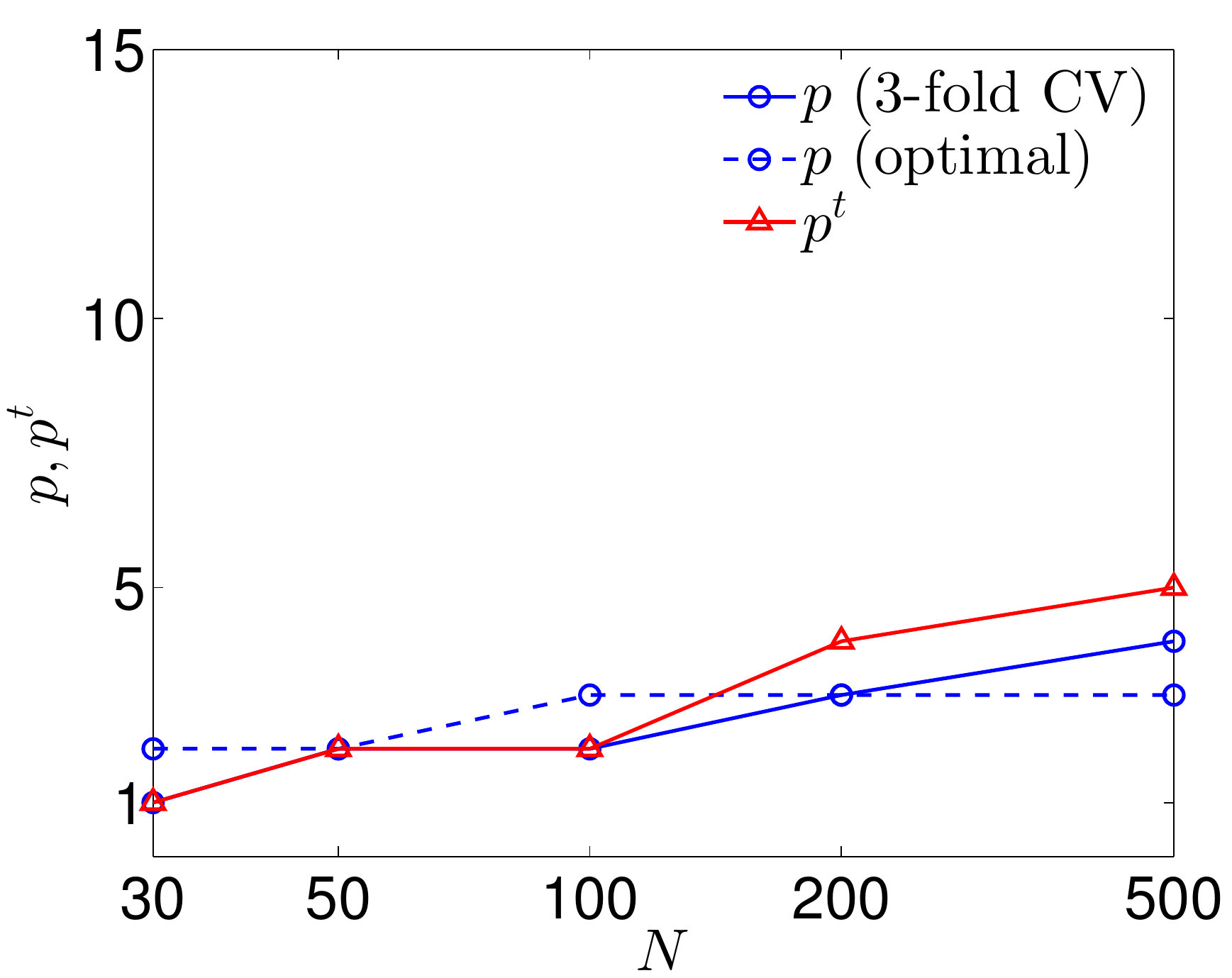}
	\includegraphics[width=0.48\textwidth]{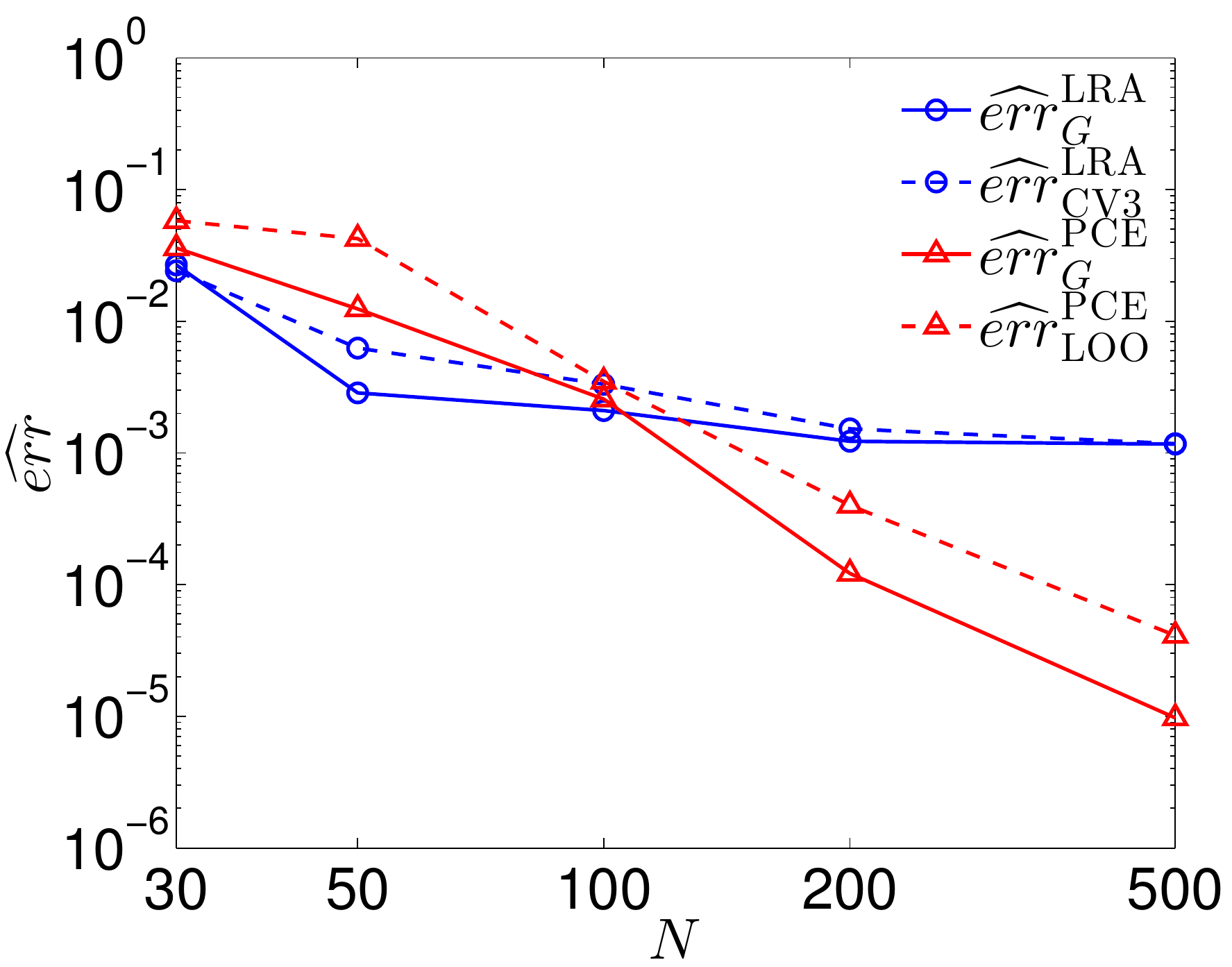}
	\caption{Truss deflection: Polynomial degrees (left) and corresponding error measures (right) of LRA and PCE meta-models based on Sobol sequences.}
	\label{fig:Truss_LRA_PCE}
\end{figure}

\begin{figure}[!ht]
	\centering
	\includegraphics[width=0.48\textwidth]{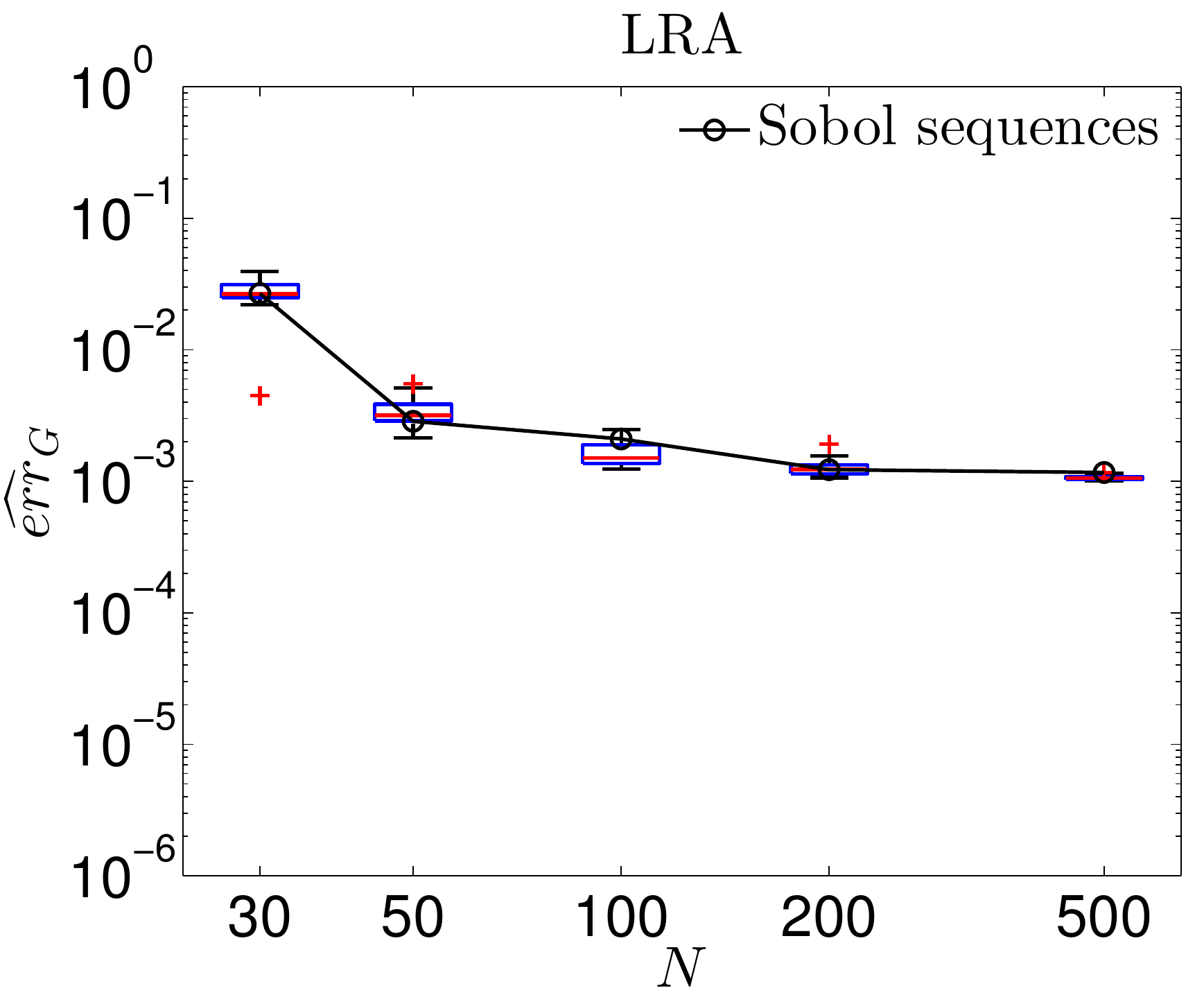}
	\includegraphics[width=0.48\textwidth]{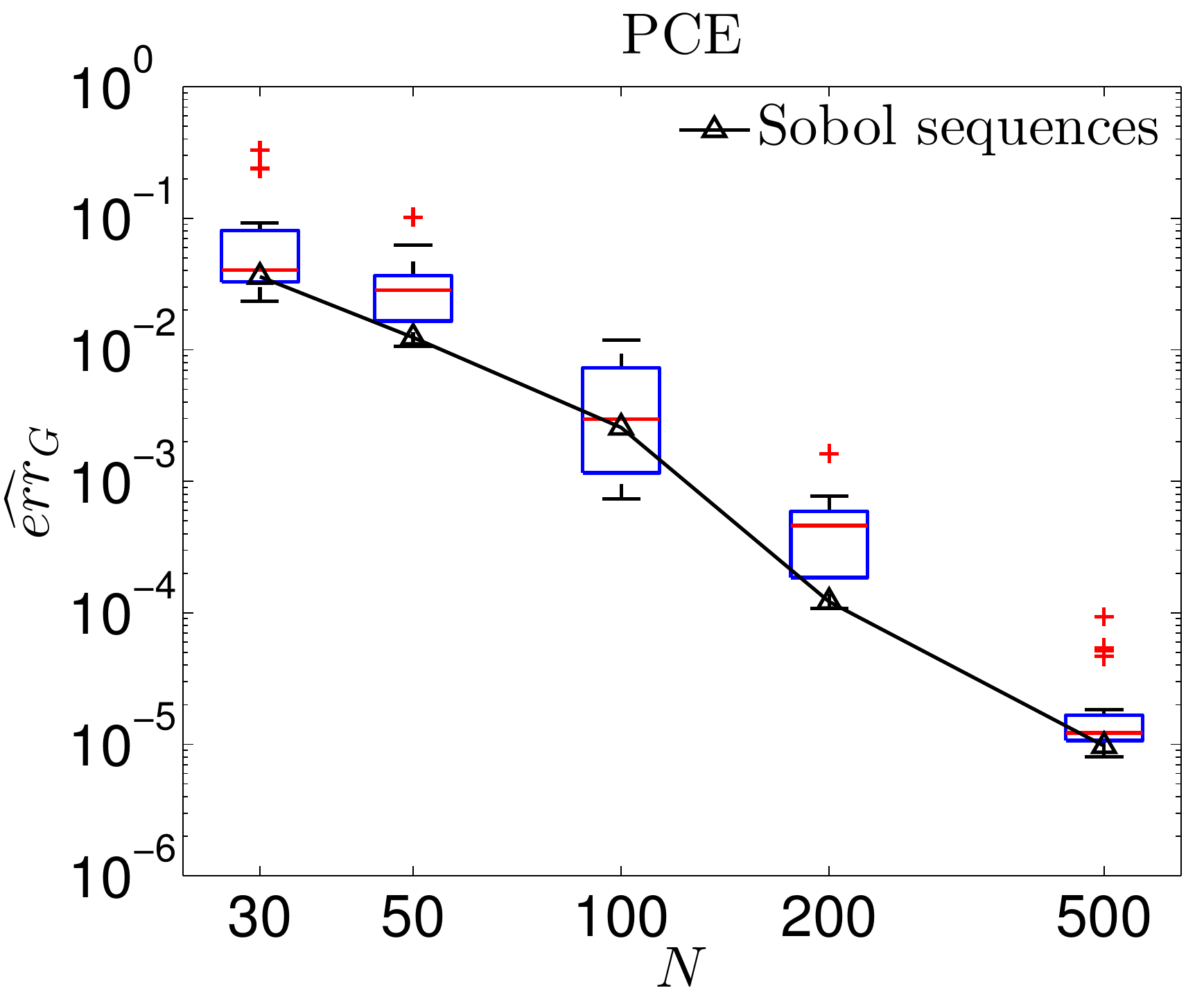}
	\caption{Truss deflection: Comparison of relative generalization errors of meta-models based on LHS (20 replications) to meta-models based on Sobol sequences.}
	\label{fig:Truss_LRA_PCE_LHS}
\end{figure}

We next assess the comparative accuracy of LRA and PCE in representing the truss-deflection model described in Section~\ref{sec:Truss_LRA}. To this end, we use EDs of size $30\leq N\leq500$, obtained with Sobol sequences as well as LHS, and a validation set of size $\nval=10^6$. We first examine the results based on Sobol sequences, presented in Figure~\ref{fig:Truss_LRA_PCE}. Similarly to the previous example, the left graph shows the polynomial degrees, while the right graph shows the corresponding ED-based and generalization errors. Again, the selected and optimal degrees $p$ of the LRA meta-models either coincide or differ by one. LRA are more accurate than PCE for the smaller EDs, but $\errGPCE$ decreases faster than $\errGLRA$ with increasing $N$, rendering PCE superior for the larger EDs. This trend is captured by the ED-based error measures, even though these are overall more pessimistic for PCE. The generalization errors of LRA and PCE meta-models built with LHS designs are shown in Figure~\ref{fig:Truss_LRA_PCE_LHS}. For each $N$, the corresponding generalization errors of the mera-models obtained with Sobol sequences are also shown for comparison reasons. The depicted boxplots correspond to 20 EDs with each representing the best among 5 random LHS designs, where the selection criterion is the maximum of the minimum distance between the points, so-called \emph{maximin} LHS designs. Clearly, the PCE errors exhibit a larger dispersion than the LRA ones. The median errors of the LHS-based LRA are very close to the corresponding errors of the LRA based on Sobol sequences. In the case of PCE, the median errors for LHS designs are either similar or larger than the corresponding errors for Sobol sequences. 

\begin{figure}[!ht]
	\centering
	\includegraphics[width=0.47\textwidth]{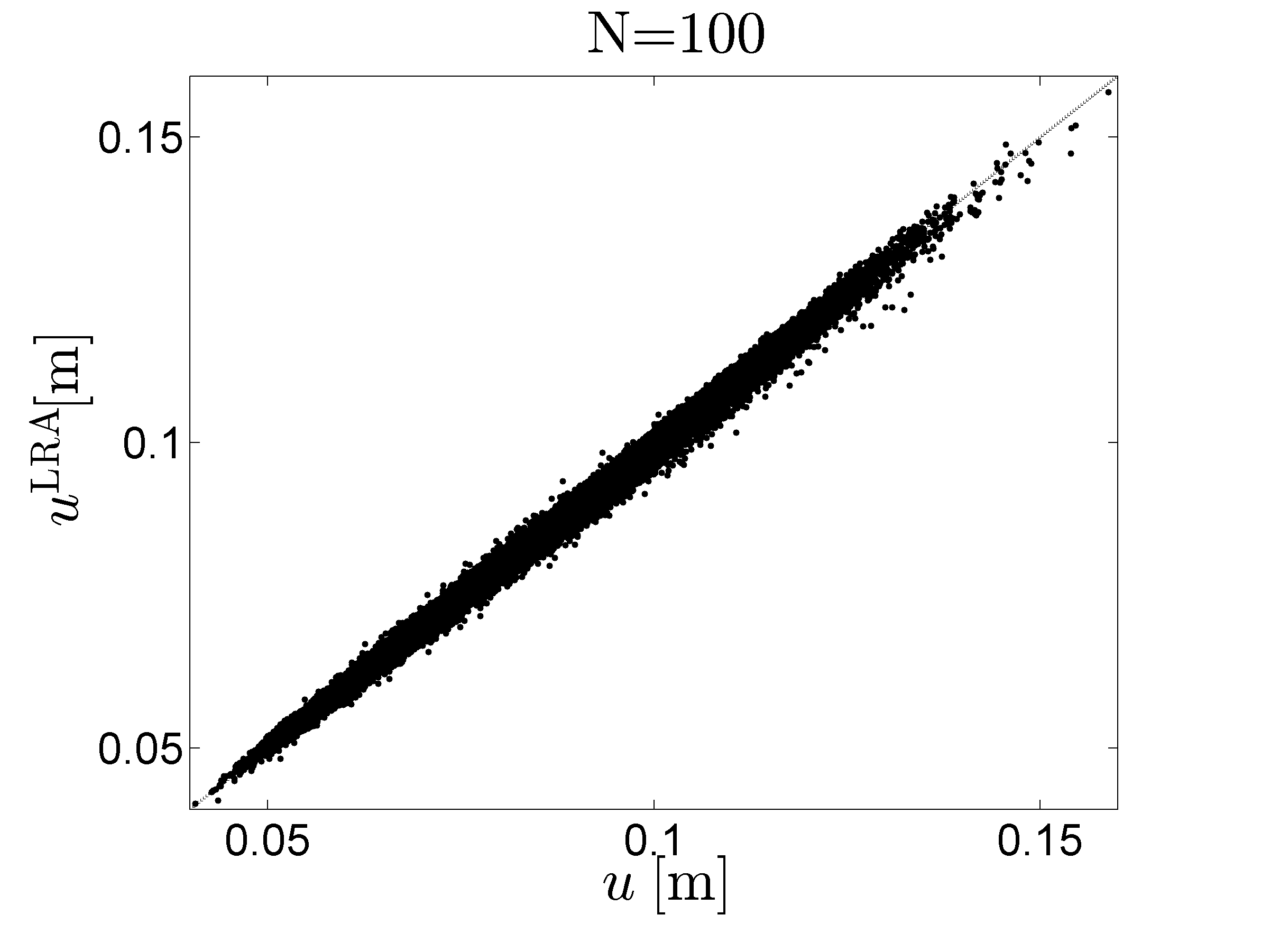}
	\includegraphics[width=0.47\textwidth]{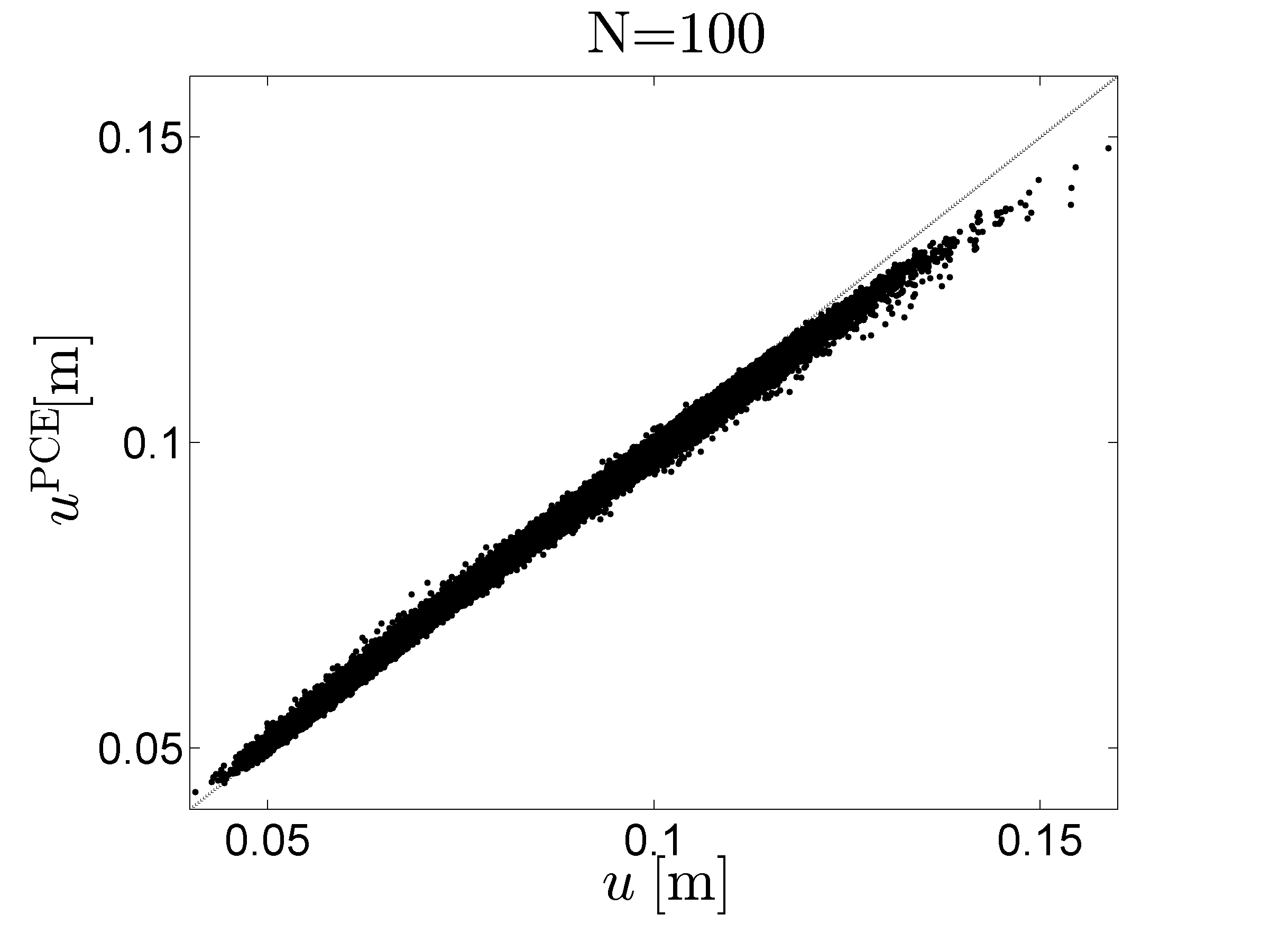}
	\caption{Truss deflection: Comparison of the exact model responses at the validation set with the respective responses of the LRA meta-model (left) and the PCE meta-model (right) for $N=100$.}
	\label{fig:Truss_LRA_PCE_N100}
\end{figure}

\begin{figure}[!ht]
	\centering
	\includegraphics[width=0.48\textwidth]{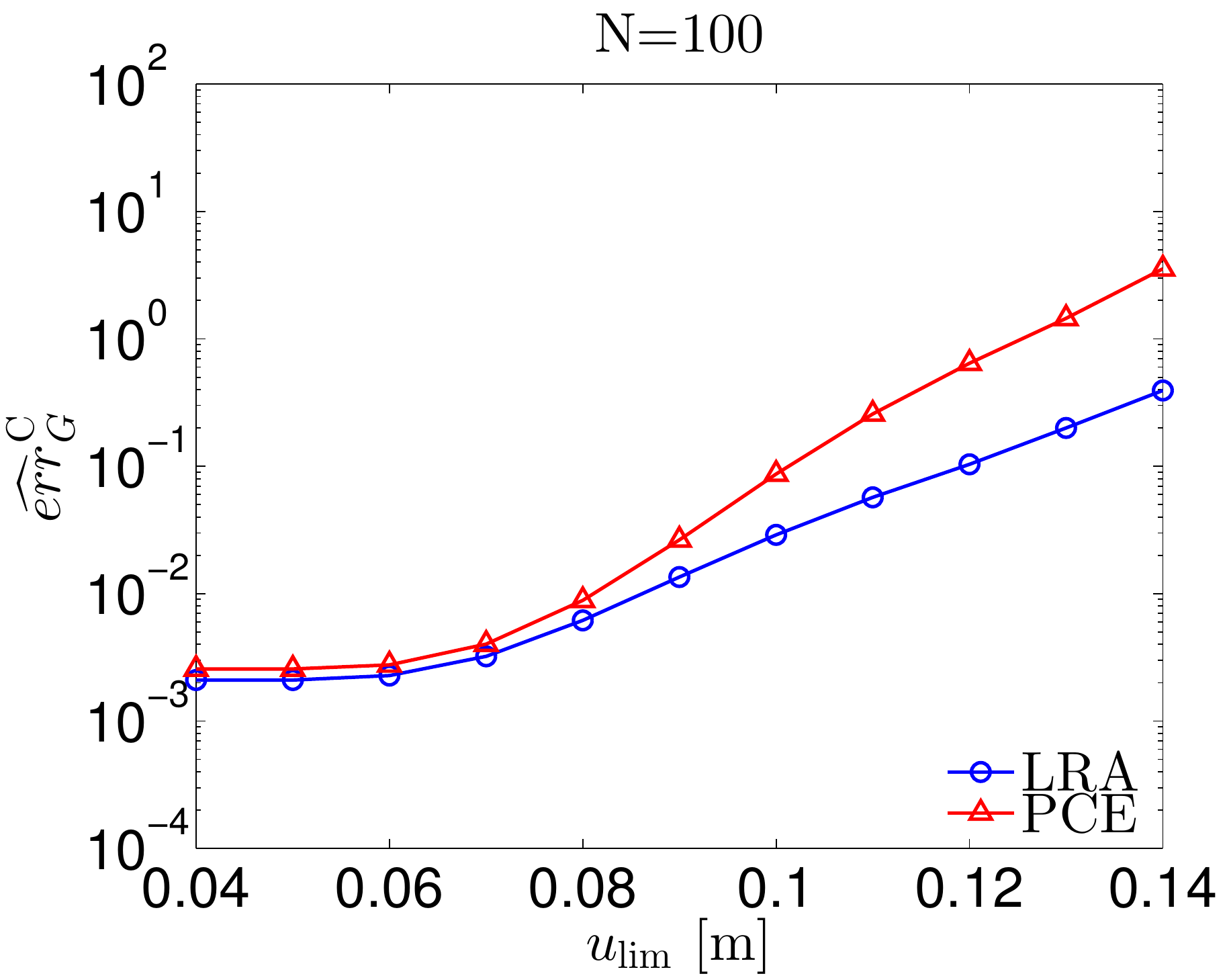}
	\includegraphics[width=0.47\textwidth]{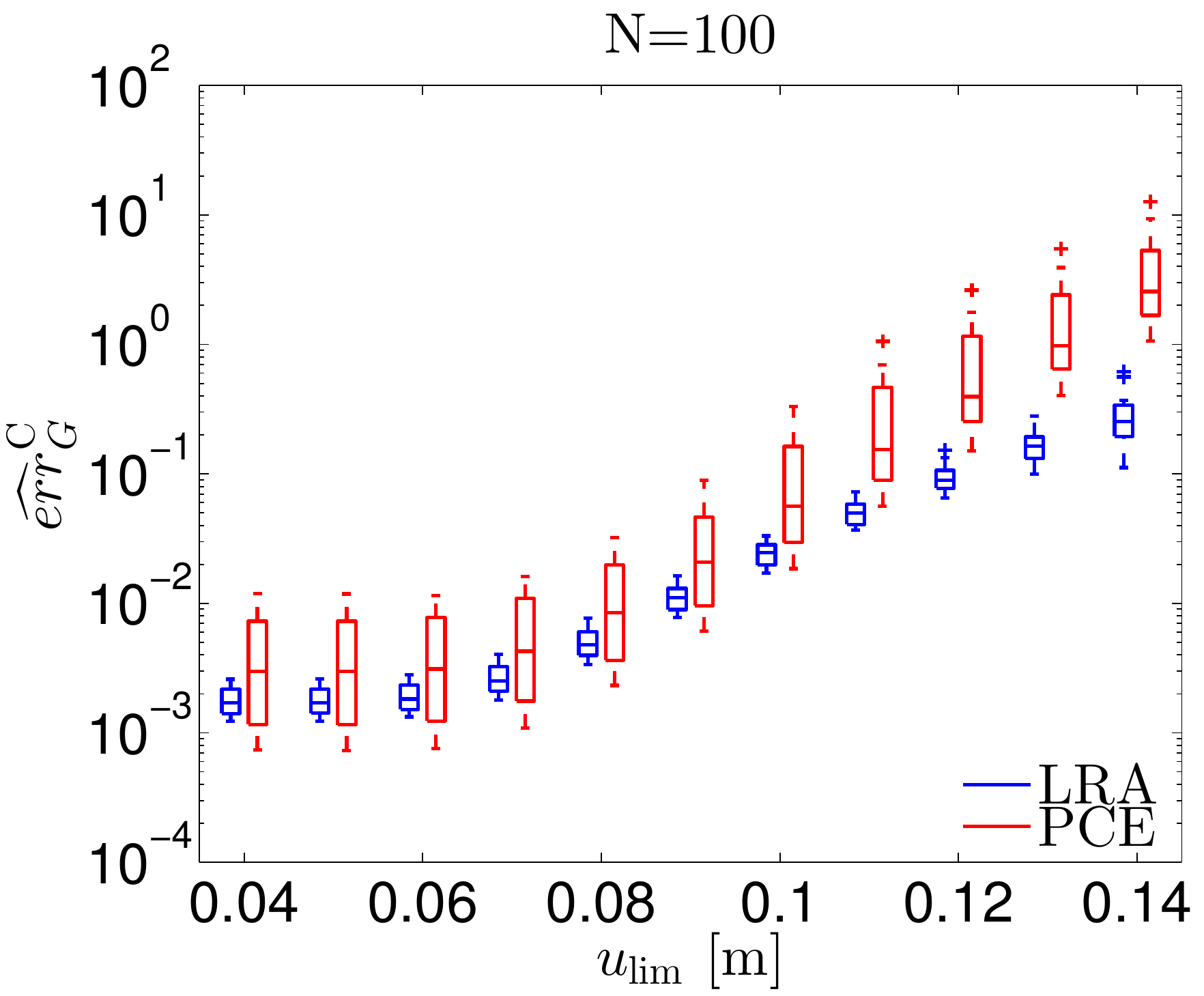}
	\caption{Truss deflection: Comparison between the LRA and PCE relative conditional generalization errors for EDs based on Sobol sequences (left) and LHS (right).}
	\label{fig:Truss_LRA_PCE_errGcond}
\end{figure}

To gain further insight into the behaviors of the two types of meta-models, we plot in Figure~\ref{fig:Truss_LRA_PCE_N100} the meta-model versus the actual-model responses at the points of the validation set for the case of Sobol sequences with $N=100$. Note in Figure~\ref{fig:Truss_LRA_PCE} that the corresponding LRA and PCE are characterized by similar generalization errors ($\errGLRA=2.10 \cdot10^{-3}, \errGPCE=2.56\cdot10^{-3}$). Obviously, LRA (left graph) provide better predictions of extreme responses, which are systematically underestimated by PCE (right graph). Because the extreme responses represent only a small fraction of the validation set, the observed differences have only a minor influence on the generalization error. In order to capture the meta-model performance in particular regions of interest, we introduce the \emph{conditional generalization error}:
\begin{equation}
\label{eq:hatErrGcond}
\widehat{Err}_G^{\rm C} =\left\|\cm-\widehat {\cm}\right\|_{\cx_{\rm val}^{\rm C}}^2.
\end{equation}
The above equation indicates that $\widehat{Err}_G^{\rm C}$  is computed similarly to $\widehat{Err}_G$  (see Eq.~(\ref{eq:ErrG_hat}), but by considering only a subset $\cx_{\rm val}^{\rm C}$ of the validation set $\cx_{\rm val}$, defined by an appropriate condition. The subsets $\cx_{\rm val}^{\rm C}$ of interest in the prediction of extreme responses are defined by:
\begin{equation}
\label{eq:Xcond}
\cx_{\rm val}^{\rm C}=\{\ve x \in \cx_{\rm val}: Y=\cm(\ve x)\geq y_{\rm lim}\}.
\end{equation}
The corresponding relative error $\widehat{err}_G^{\rm C}$ is obtained after normalization with the empirical variance of $\cy_{\rm val}^{\rm C}$, the latter denoting the set of model responses at $\cx_{\rm val}^{\rm C}$. The left graph of Figure~\ref{fig:Truss_LRA_PCE_errGcond} shows the evolution of $\widehat{err}_G^{\rm C}$ versus the response threshold $u_{\rm lim}$ for the same meta-models considered in Figure~\ref{fig:Truss_LRA_PCE_errGcond}. The conditional error increases faster for PCE than for LRA with the former being about an order of magnitude larger at the highest response threshold. The right graph of Figure~\ref{fig:Truss_LRA_PCE_errGcond} compares the LRA and PCE boxplots of $\widehat{err}_G^{\rm C}$ for the 20 LHS designs of size $N=100$. Note again the faster increase of the median PCE error, which is about an order of magnitude larger than the median LRA error at the highest response threshold.

\subsection{Frame displacement}

\begin{figure}[!ht]
	\centering
	\includegraphics[width=0.47\textwidth]{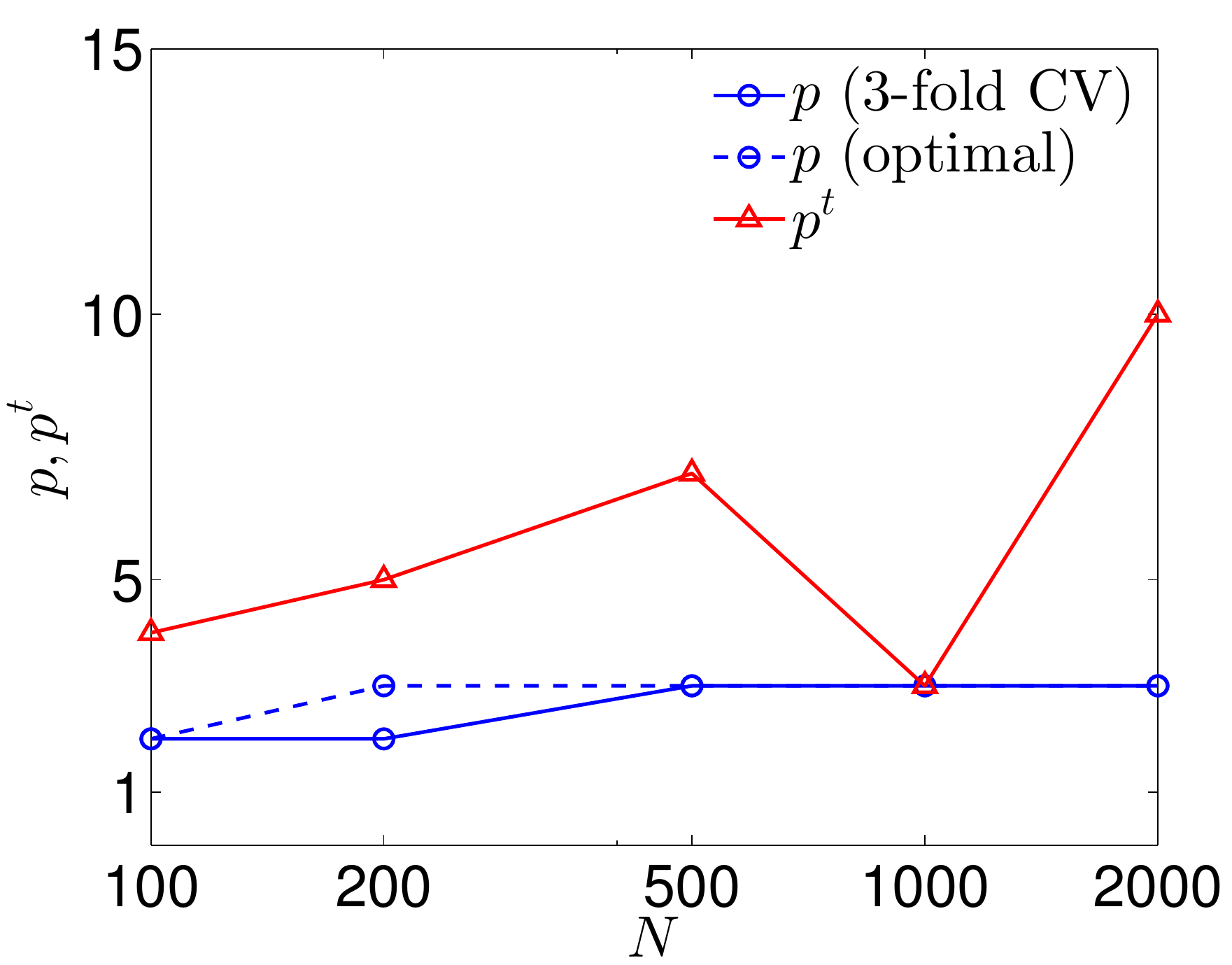}
	\includegraphics[width=0.48\textwidth]{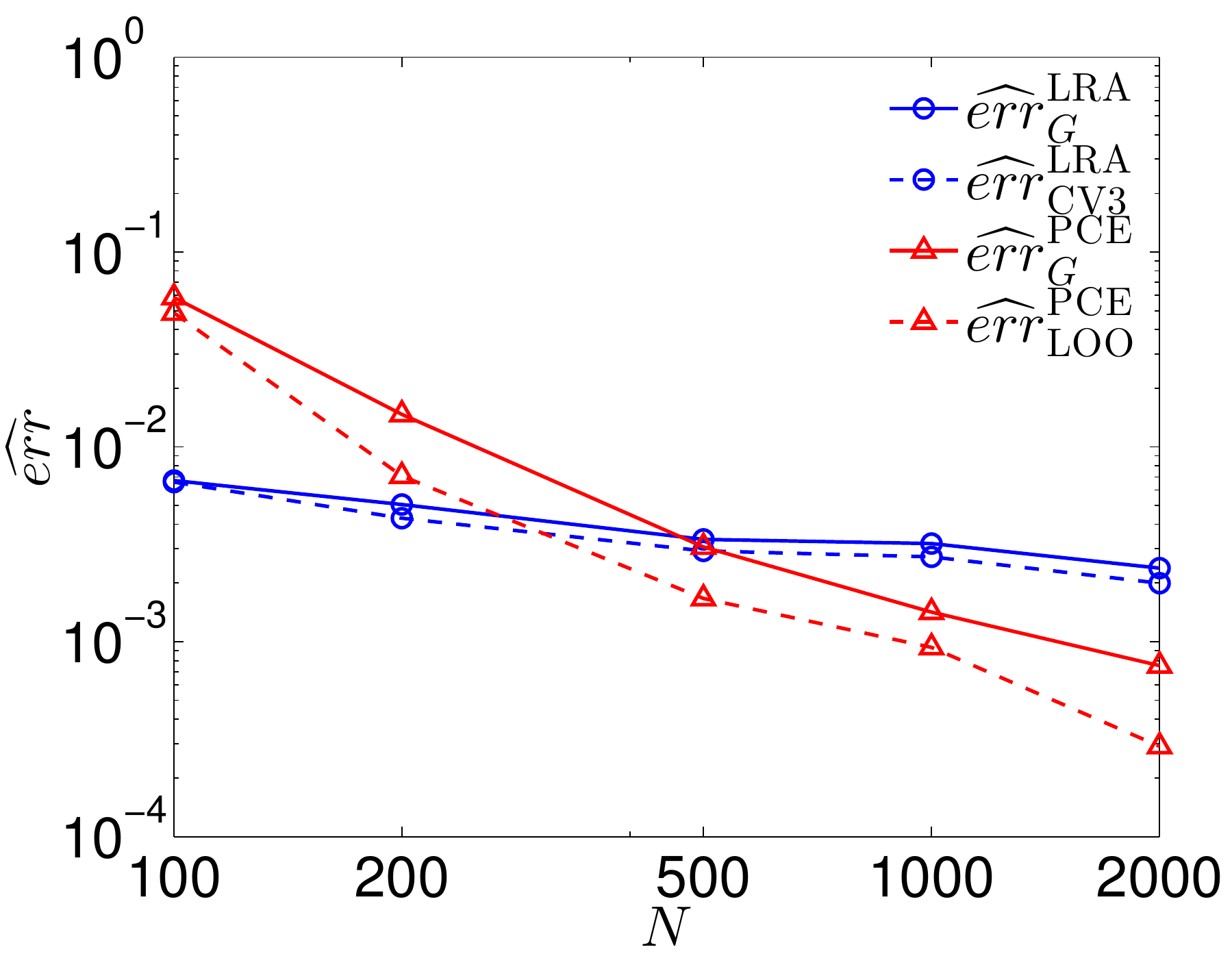}
	\caption{Frame displacement: Polynomial degrees (left) and corresponding error measures (right) of LRA and PCE meta-models based on Sobol sequences.}
	\label{fig:Frame_LRA_PCE}
\end{figure}

\begin{figure}[!ht]
	\centering
	\includegraphics[width=0.48\textwidth]{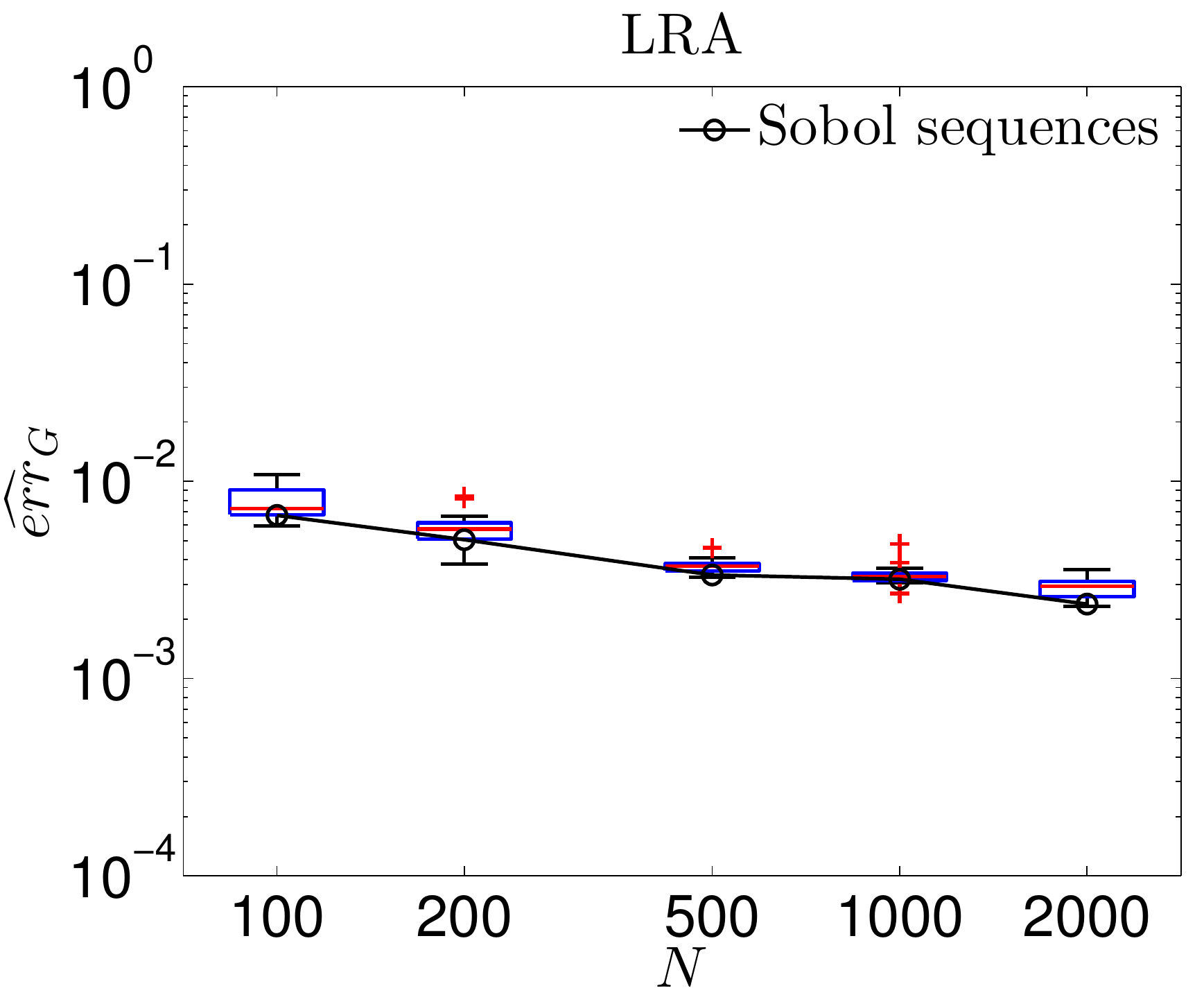}
	\includegraphics[width=0.48\textwidth]{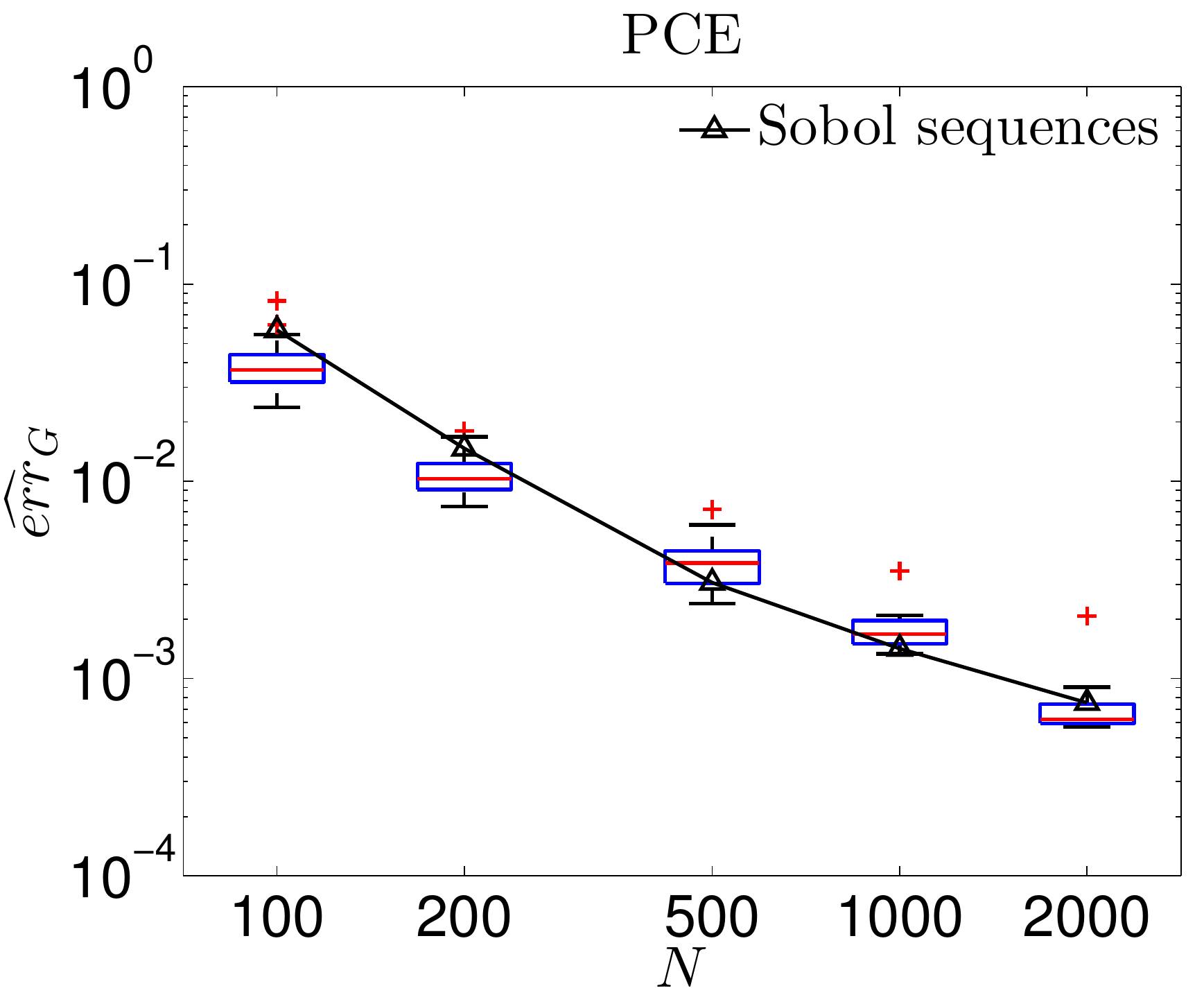}
	\caption{Frame displacement: Comparison of relative generalization errors of meta-models based on LHS (20 replications) to meta-models based on Sobol sequences.}
	\label{fig:Frame_LRA_PCE_LHS}
\end{figure}

For the frame-displacement problem, described in Section~\ref{sec:Frame_LRA}, we again examine the comparative accuracy of LRA and PCE by drawing EDs with Sobol sequences and LHS designs. The EDs are of size $100\leq N\leq2,000$, while the validation set comprises $\nval=10^6$ points. Figures~\ref{fig:Frame_LRA_PCE} and~\ref{fig:Frame_LRA_PCE_LHS} respectively show results for EDs drawn with Sobol sequences and 20 maximin LHS designs (each being the best among 5 random ones), in a manner similar to Figures~\ref{fig:Truss_LRA_PCE} and~\ref{fig:Truss_LRA_PCE_LHS} in the previous example. The polynomial degrees for the case of Sobol sequences are shown in the left graph of Figure~\ref{fig:Frame_LRA_PCE}. The selected degrees for LRA match the optimal degrees except for $N=200$. The increasing trend in the PCE degree with increasing ED size is interrupted at $N=1,000$ due to the change of the optimal truncation parameter $q$ from $0.50$ to $1.0$. Note in Figure~\ref{fig:Frame_LRA_PCE} that similarly to the truss-deflection problem, LRA exhibit smaller errors than PCE for the smaller EDs, but the PCE errors decrease faster with increasing ED size. The LRA generalization error is rather accurately estimated with the 3-fold CV approach, whereas the PCE generalization error is slightly underestimated by the LOO error. Figure~\ref{fig:Frame_LRA_PCE_LHS} shows that the LRA errors for the LHS designs exhibit a slightly smaller dispersion than the PCE ones. In all cases, the median errors for the LHS designs are fairly close to the corresponding errors for the Sobol sequences.

\begin{figure}[!ht]
	\centering
	\includegraphics[width=0.47\textwidth]{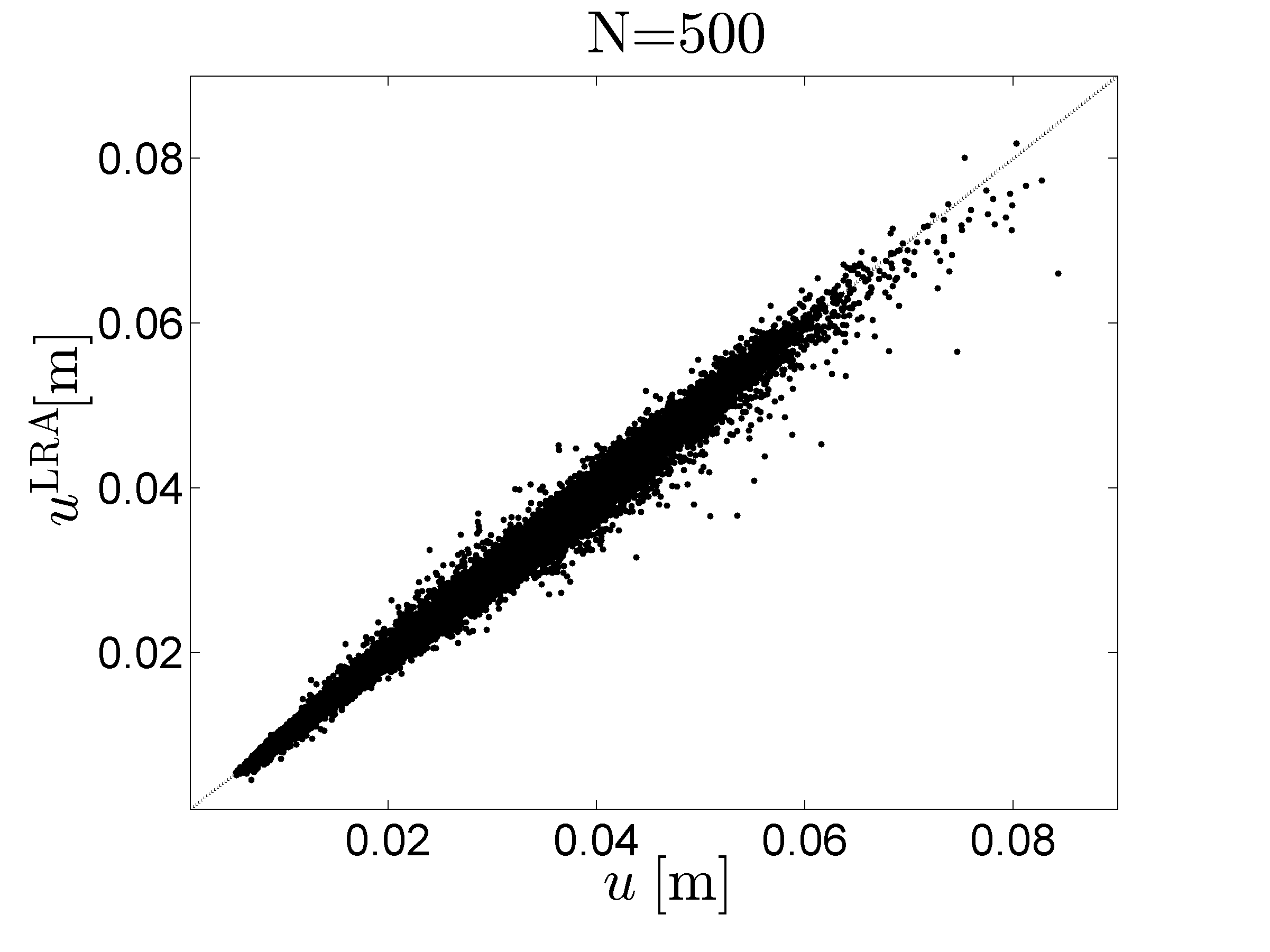}
	\includegraphics[width=0.47\textwidth]{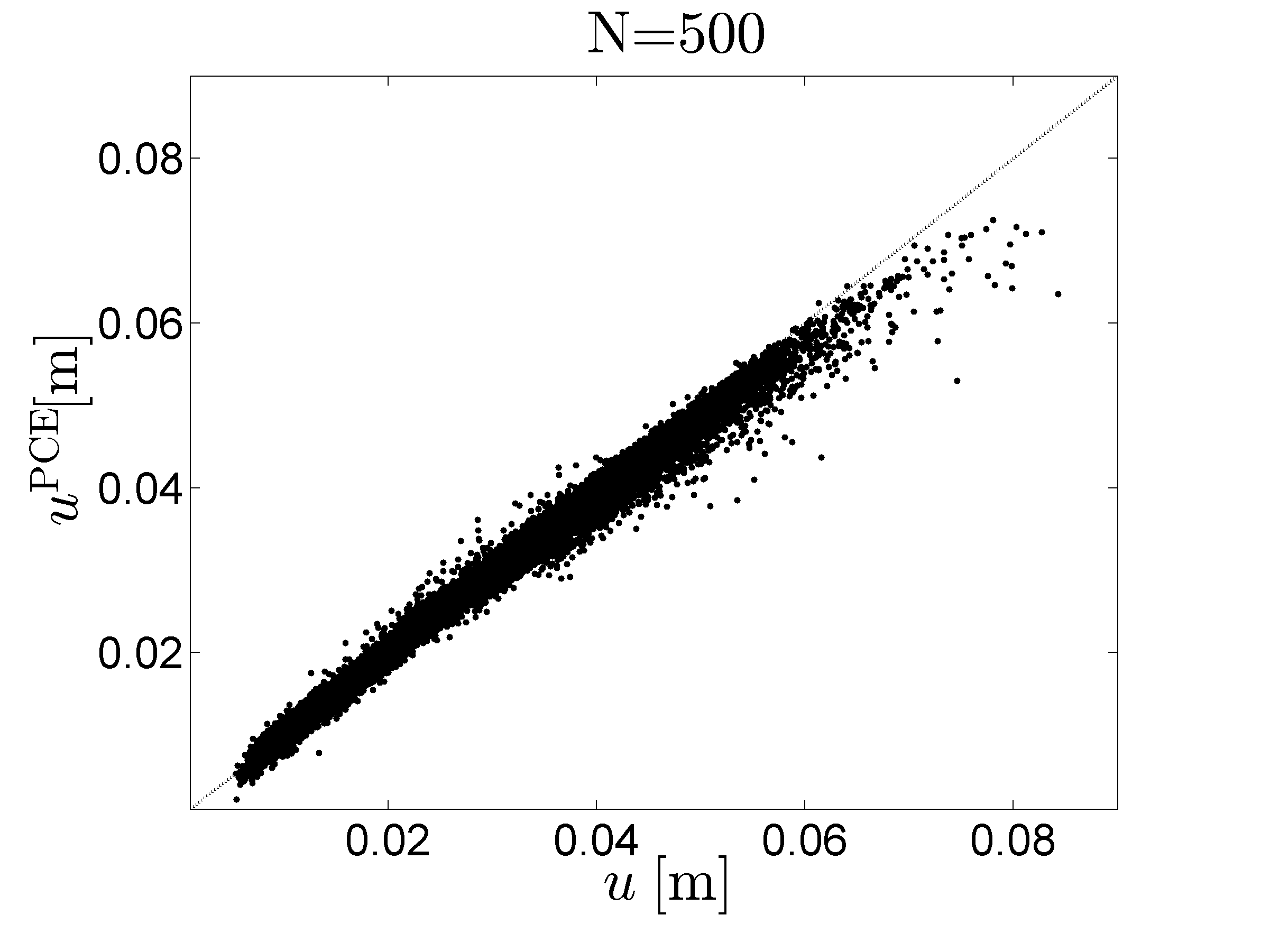}
	\caption{Frame displacement: Comparison of the exact model responses at the validation set with the respective responses of the LRA meta-model (left) and the PCE meta-model (right) for $N=500$.}
	\label{fig:Frame_LRA_PCE_N500}
\end{figure}

\begin{figure}[!ht]
	\centering
	\includegraphics[width=0.48\textwidth]{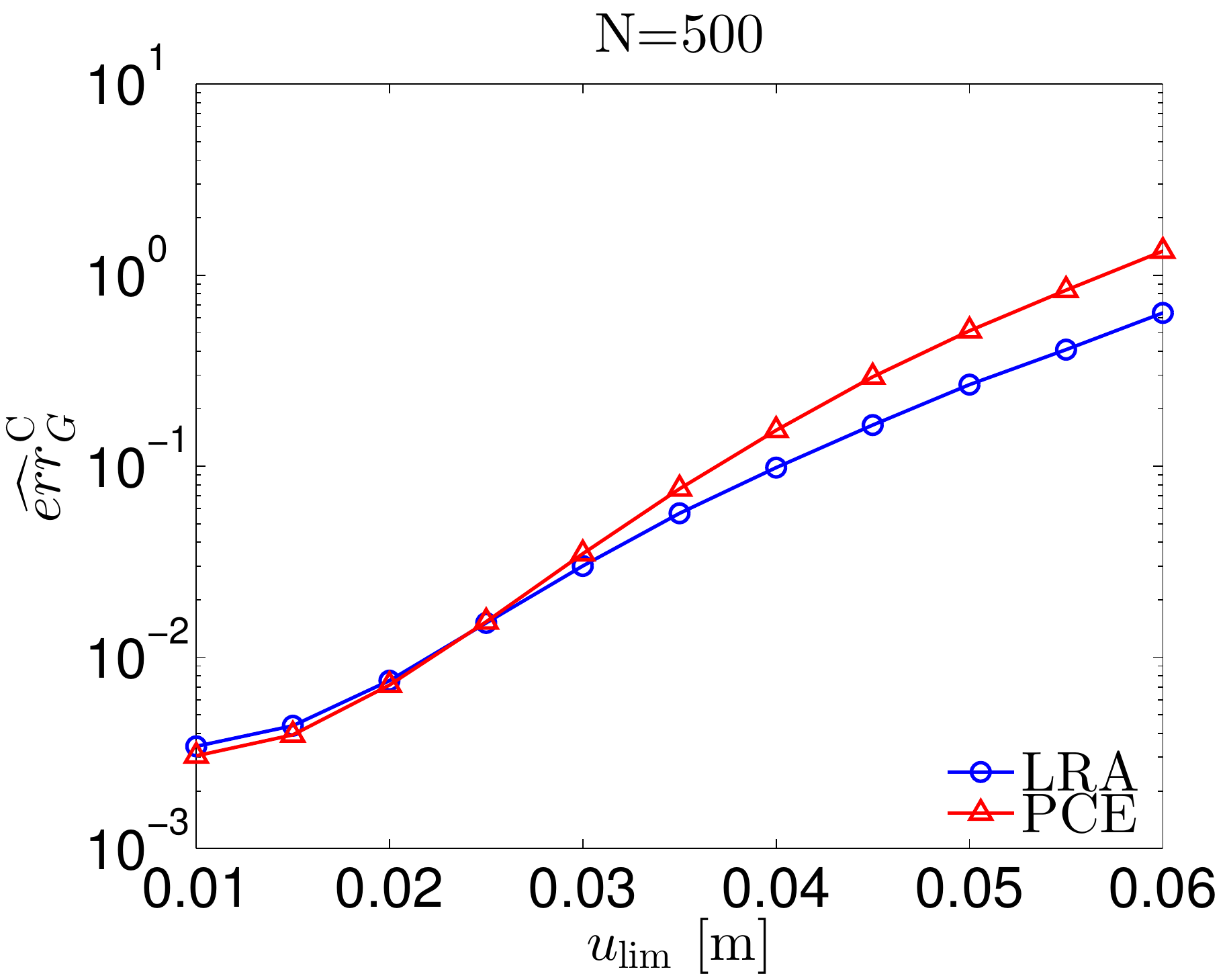}
	\includegraphics[width=0.47\textwidth]{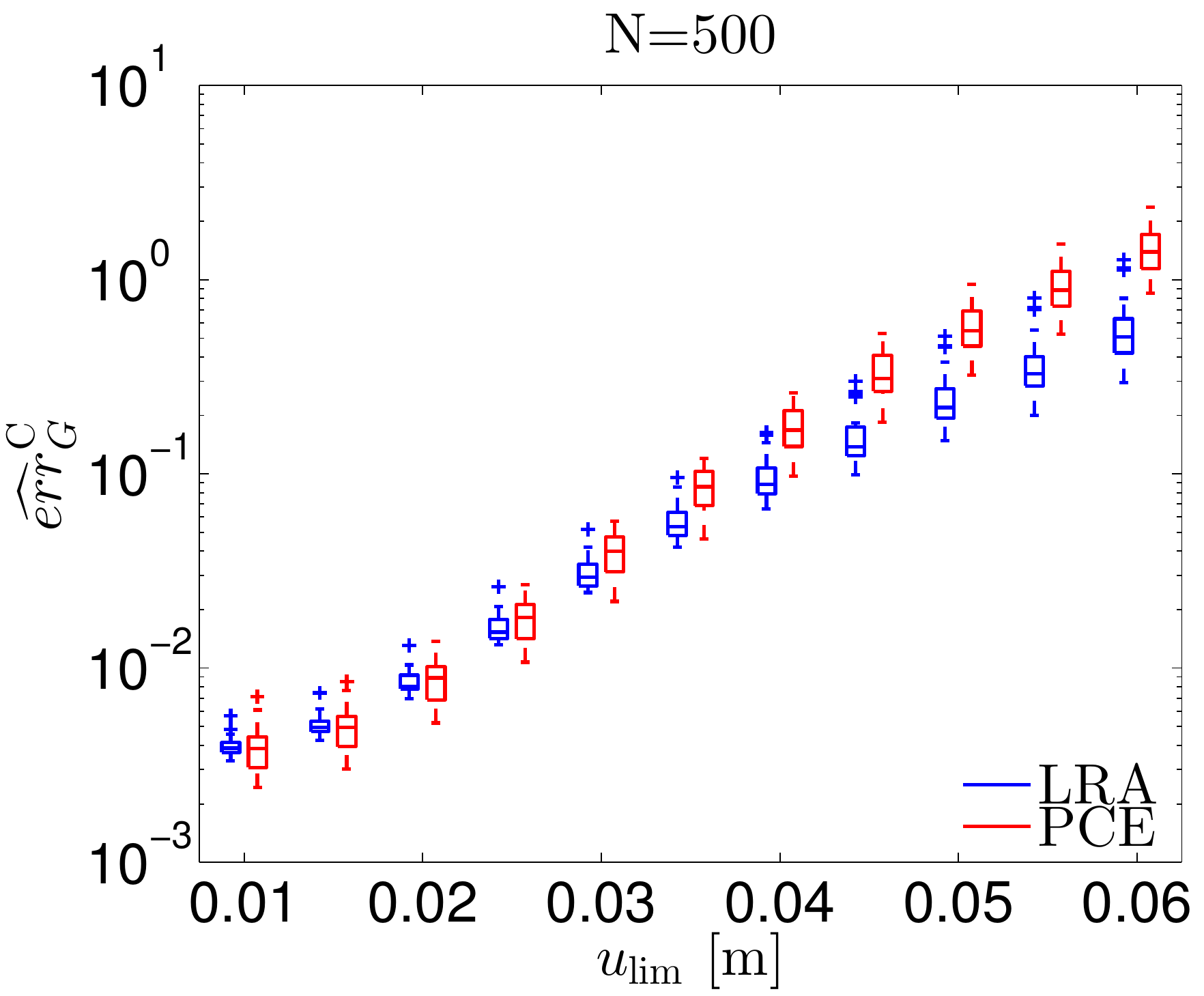}
	\caption{Frame displacement: Comparison between the LRA and PCE relative conditional generalization errors for EDs based on Sobol sequences (left) and LHS (right).}
	\label{fig:Frame_LRA_PCE_errGcond}
\end{figure}

Note in Figures~\ref{fig:Frame_LRA_PCE} and~\ref{fig:Frame_LRA_PCE_LHS} that the two types of meta-models exhibit similar generalization errors at $N=500$. However, as in the truss-deflection problem, LRA outperform PCE in the prediction of extreme responses. This can be observed in Figure~\ref{fig:Frame_LRA_PCE_N500}, which shows the meta-model versus the actual-model responses at the validation set for an example LHS design of size $N=500$. The superior performance of LRA at the upper tail of the response distribution is not reflected on the generalization errors ($\errGLRA=3.68 \cdot10^{-3}, \errGPCE=3.47\cdot10^{-3}$), but can be captured by the conditional generalization error (Eq.~(\ref{eq:hatErrGcond})-(\ref{eq:Xcond})). For the meta-models developed with EDs of size $N=500$, Figure~\ref{fig:Frame_LRA_PCE_errGcond} depicts the evolution of $\widehat{err}_G^{\rm C}$ with increasing response threshold; the left graph of the figure shows this error for PCE and LRA based on Sobol sequences, while the right graph shows respective boxplots for the 20 LHS designs. In both graphs, we observe that the conditional errors of LRA and PCE are similar for the lower responses thresholds, but the PCE errors are larger for the higher ones.

\subsection{Heat conduction}

\begin{figure}[!ht]
	\centering
	\includegraphics[width=0.47\textwidth]{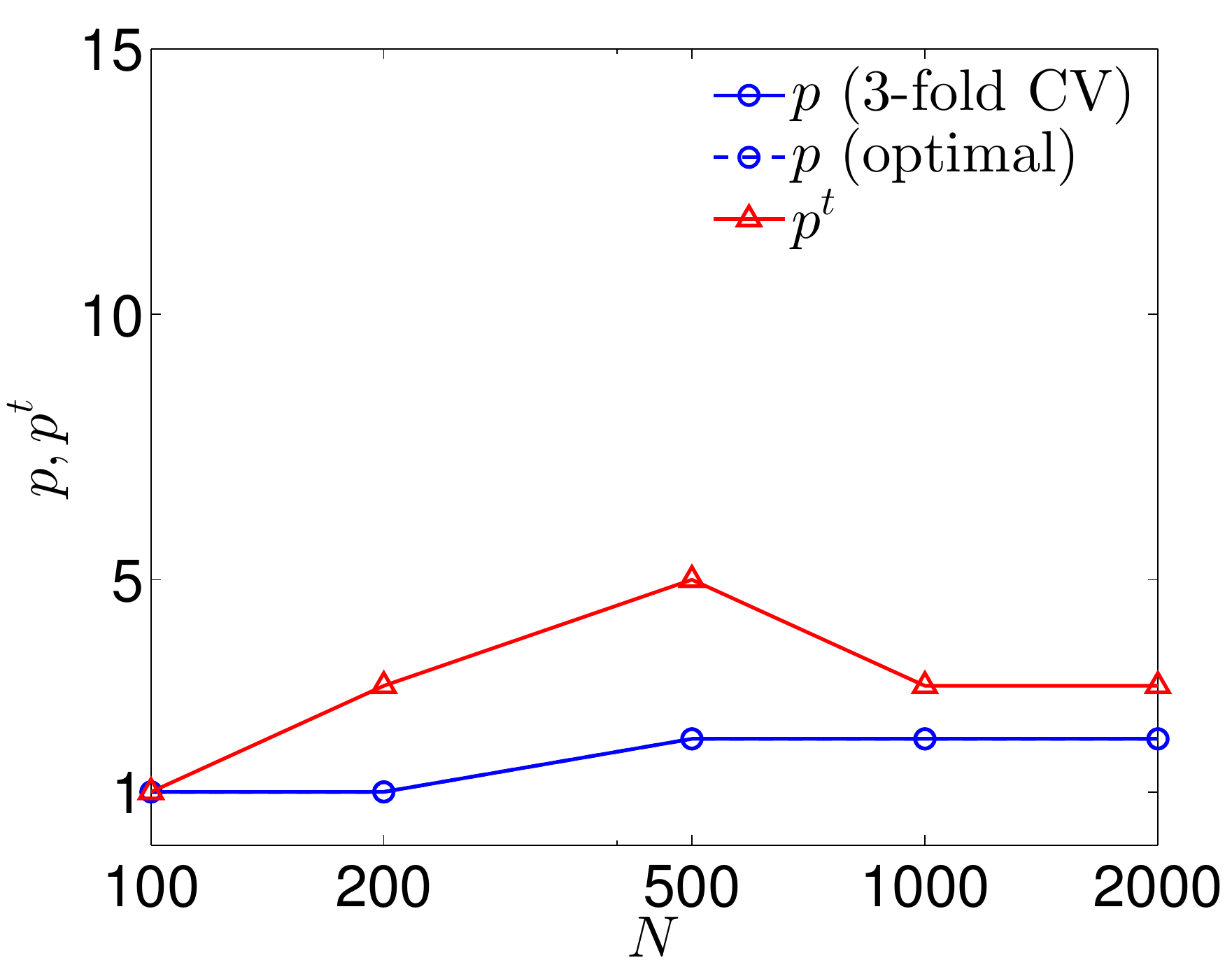}
	\includegraphics[width=0.48\textwidth]{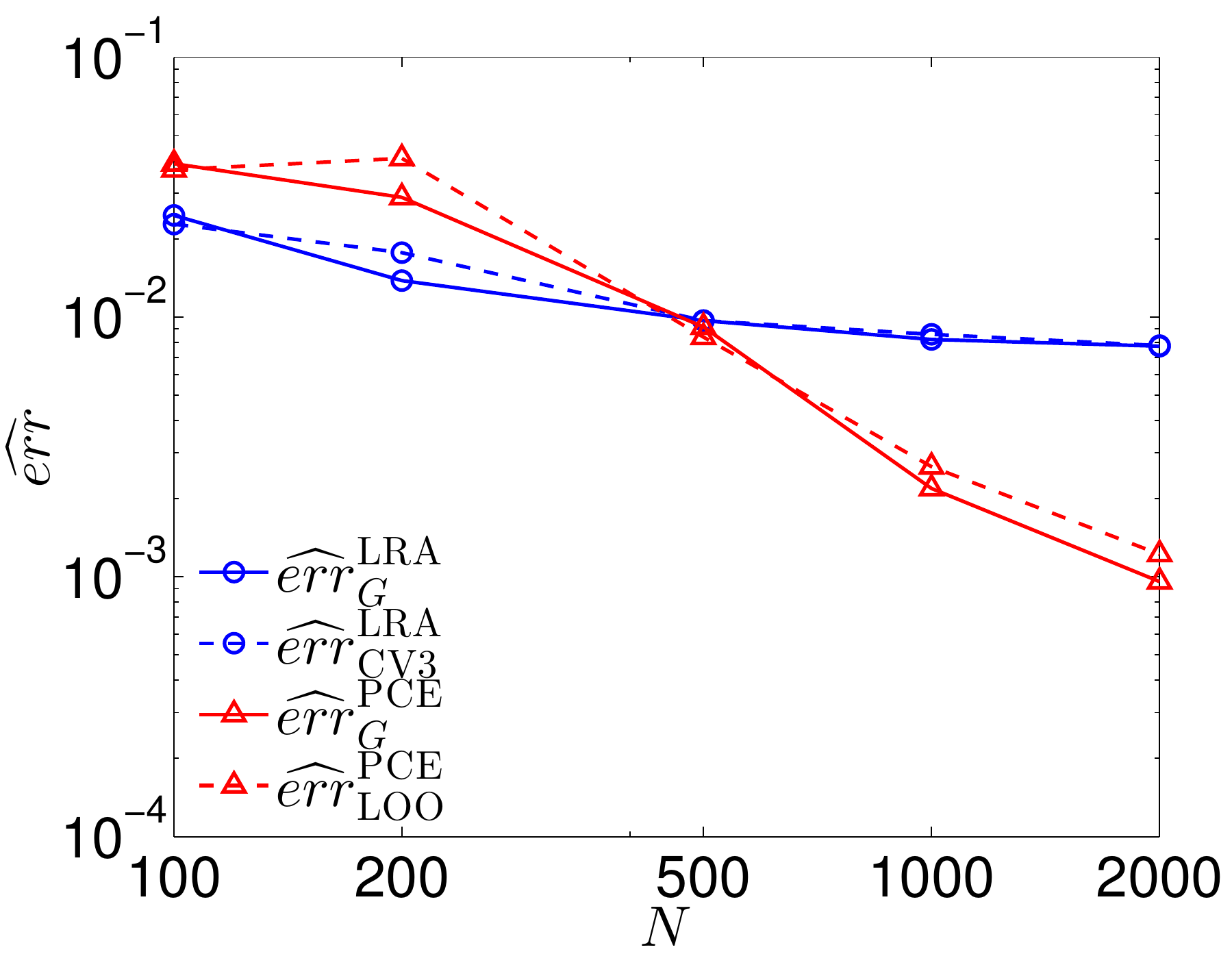}
	\caption{Heat conduction: Polynomial degrees (left) and corresponding error measures (right) of LRA and PCE meta-models based on Sobol sequences.}
	\label{fig:DiffusionRF_LRA_PCE}
\end{figure}

In this last example, we consider again the heat-conduction problem investigated in Section~\ref{sec:DiffusionRF_LRA}. We develop LRA and PCE meta-models for the average temperature $\widetilde{T}$ in domain B (see Figure~\ref{fig:RF_domain}) using EDs of size  $100 \leq N\leq 2,000$ based on Sobol sequences. To assess the comparative accuracy of the meta-models, we use a validation set of size $n_{\rm val}=10^4$. The left and right graphs of Figure~\ref{fig:DiffusionRF_LRA_PCE} respectively show the polynomial degrees and error measures for the two types the meta-models. Note that the optimal LRA degree $p$ is accurately identified by 3-fold CV for all considered EDs. The increasing trend in the PCE degree $p^t$ with increasing $N$ is interrupted at $N=1,000$ due to the change of the optimal truncation parameter $q$ from $0.50$ to $0.75$. The errors demonstrate similar trends with those observed for the previously examined finite-element models, \ie $\errGLRA$ is smaller than $\errGPCE$ for the smaller EDs, but the latter decreases faster with increasing $N$. For both types of meta-models, the generalization errors are estimated fairly well by the corresponding ED-based measures.

\begin{figure}[!ht]
	\centering
	\includegraphics[width=0.47\textwidth]{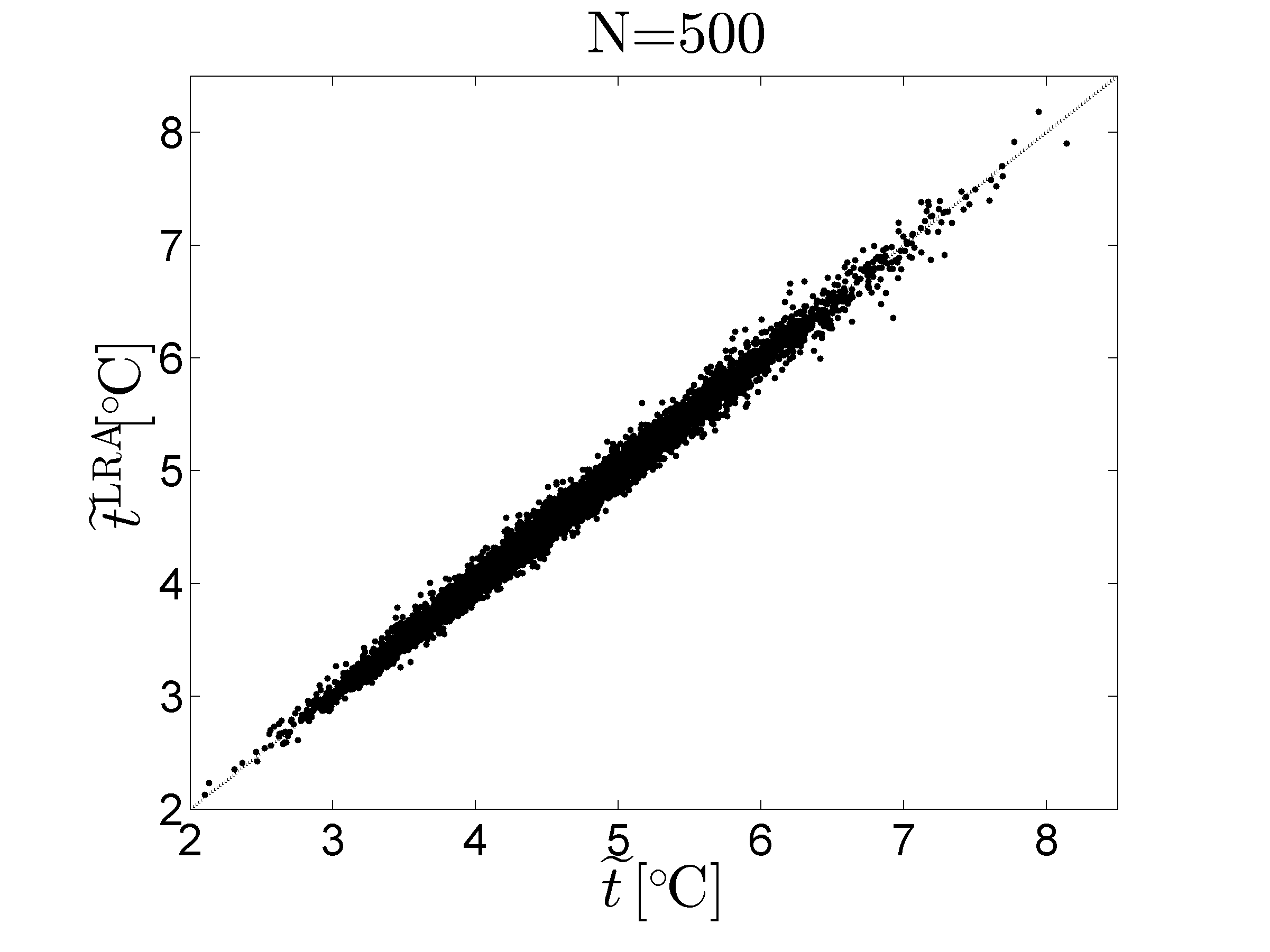}
	\includegraphics[width=0.47\textwidth]{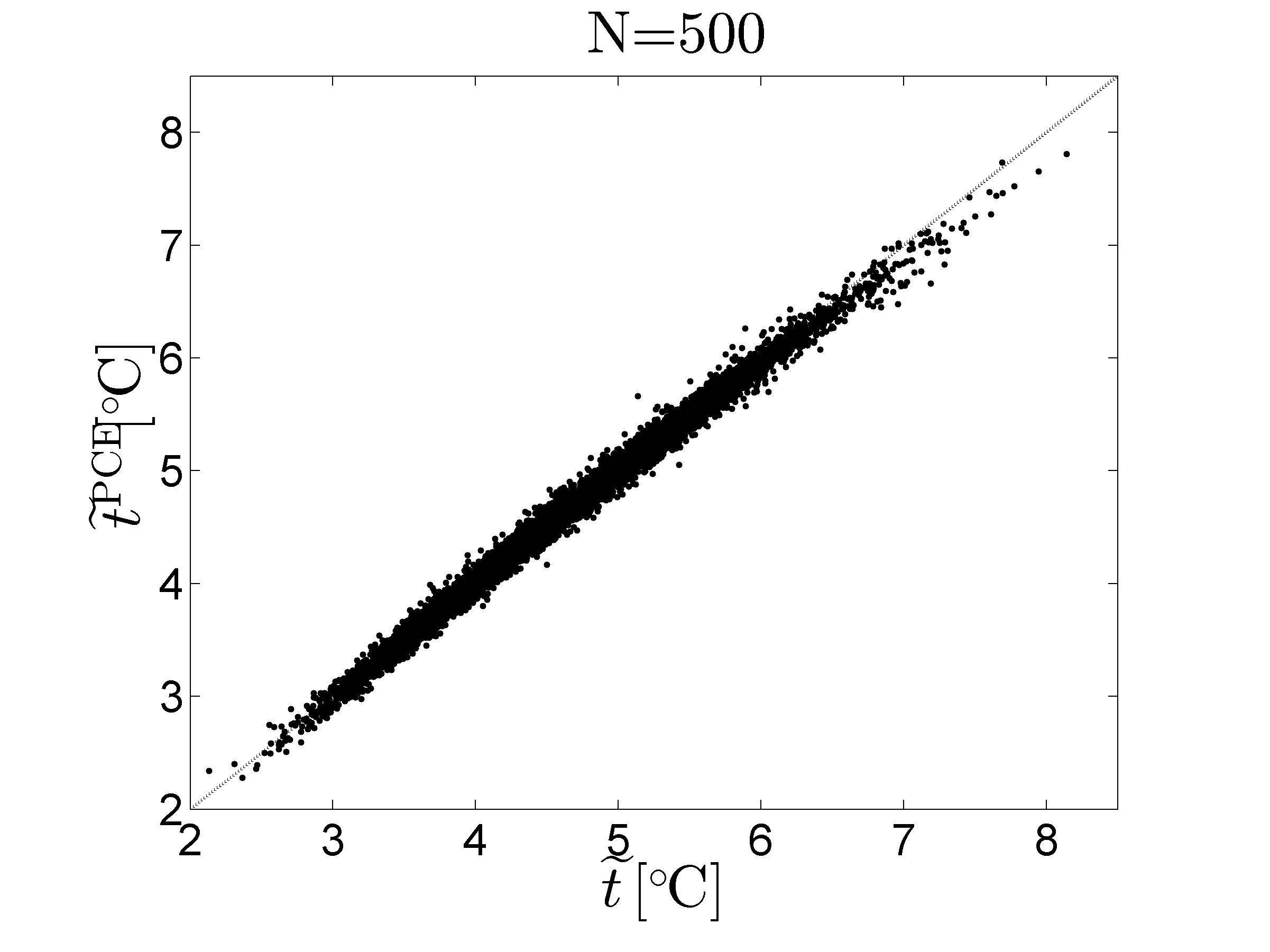}
	\caption{Heat conduction: Comparison of the exact model responses at the validation set with the respective responses of the LRA meta-model (left) and the PCE meta-model (right) for $N=500$.}
	\label{fig:DiffusionRF_LRA_PCE_N500}
\end{figure}

\begin{figure}[!ht]
	\centering
	\includegraphics[width=0.48\textwidth]{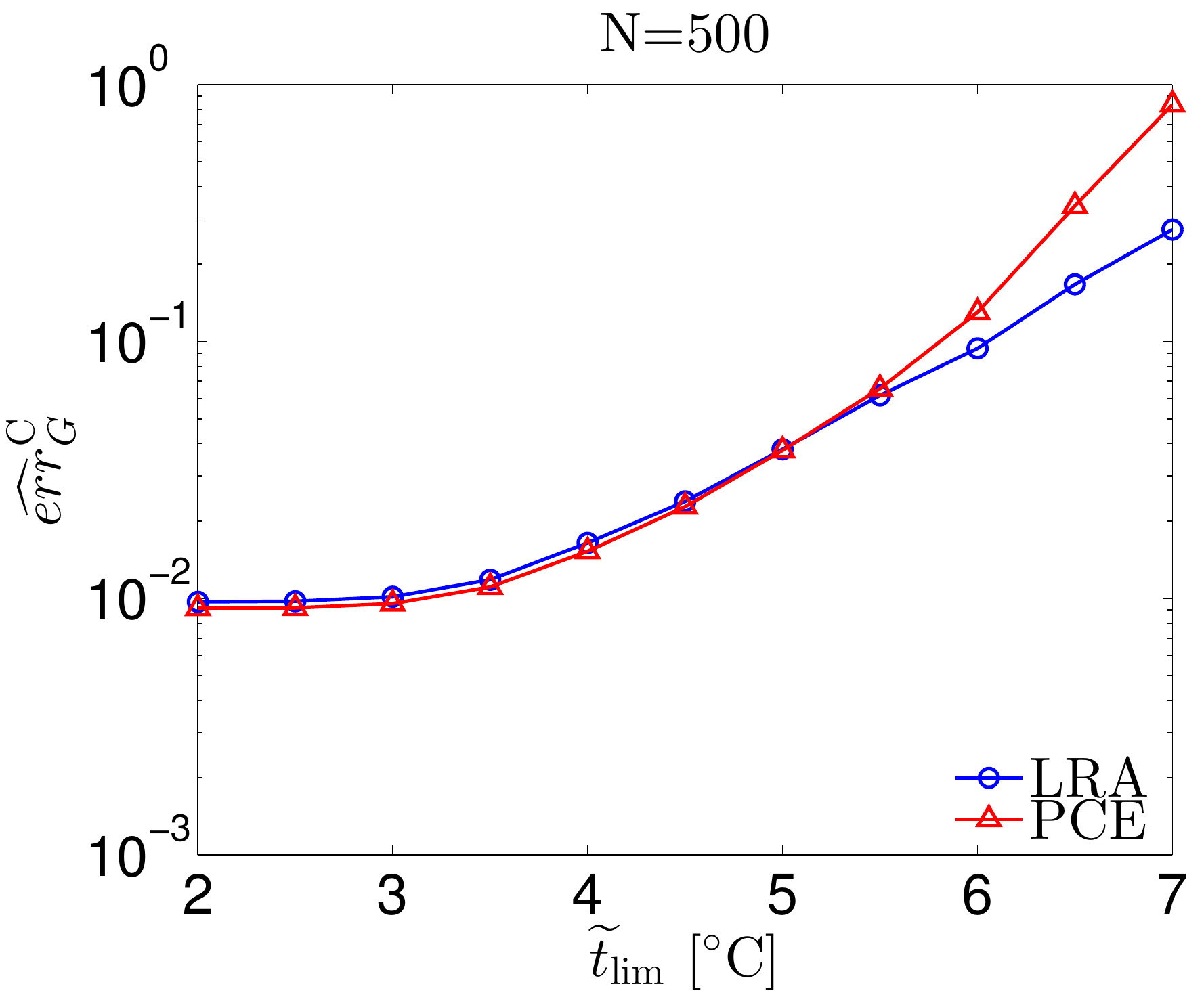}
	\caption{Heat conduction: Comparison between the LRA and PCE relative conditional generalization errors.}
	\label{fig:DiffusionRF_LRA_PCE_errGcond}
\end{figure}

For the ED of size $N=500$, the generalization errors of LRA and PCE are nearly equal ($\errGLRA=9.68 \cdot10^{-3}, \errGPCE=9.15\cdot10^{-3}$); however, as in the previous examples, LRA provide better predictions of the extreme responses. The left and right graphs of Figure~\ref{fig:DiffusionRF_LRA_PCE_N500} respectively show the LRA and PCE responses versus the actual model responses at the validation set. The PCE predictions of the extreme responses exhibit a negative bias, which does not appear in the LRA predictions. The lower accuracy of the PCE meta-model at the tail of the response distribution is captured by the conditional generalization error, which is plotted in Figure~\ref{fig:DiffusionRF_LRA_PCE_errGcond} for varying response thresholds.

\section{Conclusions}

In this study, a newly emerged class of meta-models based on canonical low-rank approximations (LRA) was confronted to the popular sparse polynomial chaos expansions (PCE). Both meta-modeling approaches hold strong promise against the curse of dimensionality that often poses a major challenge in real-life uncertainty propagation problems. We examined LRA and PCE meta-models developed in an non-intrusive manner by relying on the same polynomial families to build the basis functions; the polynomial coefficients in both meta-model types were obtained by solving error-minimization problems. In the considered sparse PCE approach, a single large-size minimization problem is efficiently solved by retaining only the significant basis terms, as identified with the least angle regression algorithm. In the canonical LRA approach, the tensor-product form of the basis is retained, leading to a series of small-size minimization problems that were herein solved using ordinary least-squares (OLS).

We first shed light on issues pertinent to the construction of LRA through extensive numerical investigations. To this end, we considered a particular greedy algorithm comprising a series of pairs of a correction and an updating step. In a correction step, a rank-one component is built based on the sequential updating of the polynomial coefficients in different dimensions, whereas in an updating step, the coefficients of the new set of rank-one components are determined. Two approaches for rank selection were examined; the first relies on error estimation with 3-fold cross-validation (CV), whereas the second, based on a simpler computation, uses a corrected version of the leave-one-out (LOO) error in the updating step. The lack of orthogonality of the regressors led to excessive correction factors of the LOO error in certain cases, particularly for higher ranks. On the other hand, the 3-fold CV error was found an overall reliable estimator of the generalization error, also appropriate for the selection of the polynomial degree. Finally, a stopping criterion for the correction step was examined, combining a threshold on the differential empirical error in two successive iterations with a maximum allowable iteration number. Effects of the two thresholds on the LRA accuracy were investigated and appropriate values were proposed.

The comparative accuracy of the particular LRA and PCE approaches was investigated in three problems involving finite-element models of varying dimensionality $M$. Experimental designs (EDs) of varying sizes were drawn with Sobol sequences and Latin hypercube sampling (LHS). The errors obtained by using Sobol sequences as well as the median errors obtained by using LHS designs were lower for LRA when small EDs (approximately up to size $10 M-20 M$) were considered. However, the PCE errors decreased faster with increasing ED size. In cases when the two types of meta-models exhibited similar generalization errors, LRA provided better estimates of the extreme responses. The superiority of LRA in predicting  extreme responses was quantified by means of the introduced conditional generalization error. This finding renders LRA particularly promising for problems where the accuracy in estimating the tail of the response distribution is important, such as reliability applications \citep{Konakli2015ICASP,Konakli2015ESREL}. It is further emphasized that by relying on a series of OLS minimizations, the construction of LRA involved simpler computations compared to PCE. Because of the small size of these minimization problems, the LRA construction also required far less computer memory.

\bibliographystyle{chicago}
\bibliography{bib}

\end{document}